\documentclass[10pt]{article}
\usepackage{amssymb}
\usepackage{latexsym}

\usepackage[a4paper,top=5.35cm,bottom=5.35cm,left=4.75cm,right=4.75cm,%
            bindingoffset=0mm]{geometry}

\usepackage{setspace}
\makeatletter
\let\@xp=\expandafter

\let\@nx=\noexpand
\newcommand{\numberwithin}[3][\arabic]{%
  \@ifundefined{c@#2}{\@nocounterr{#2}}{%
    \@ifundefined{c@#3}{\@nocnterr{#3}}{%
      \@addtoreset{#2}{#3}%
      \@xp\xdef\csname the#2\endcsname{%
        \@xp\@nx\csname the#3\endcsname .\@nx#1{#2}}}}%
}
\makeatother
\numberwithin{equation}{section}
\newtheorem{thm}[equation]{Theorem}
\newtheorem{prop}[equation]{Proposition}
\newtheorem{defn}[equation]{Definition}
\newtheorem{rem}[equation]{Remark}
\newtheorem{lem}[equation]{Lemma}
\newtheorem{corol}[equation]{Corollary}
\title{Regularizing properties\\
 of the double layer potential
 \\
of  second order elliptic differential operators
}

\author{Francesco Dondi \& Massimo Lanza de Cristoforis}
\date\ 

\begin{document}

\maketitle

\noindent
{\bf Abstract:} We prove the validity of regularizing properties of a  double layer potential associated to the fundamental solution of a {\em nonhomogeneous} second order elliptic differential operator with constant coefficients in Schauder spaces by exploiting an 
explicit formula for the tangential derivatives of the double layer potential itself. We also introduce ad hoc norms for kernels of integral operators in order to prove continuity results of  integral operators upon variation of the kernel, which we apply to layer potentials.
 
\vspace{\baselineskip}


\section{Introduction.} In this paper, we consider the double layer potential associated to the fundamental solution of a second order differential operator with constant coefficients. Throughout the paper, we assume that
\[
n\in {\mathbb{N}}\setminus\{0,1\}\,,
\]
where ${\mathbb{N}}$ denotes the set of natural numbers including $0$. Let $\alpha\in]0,1[$, $m\in {\mathbb{N}}\setminus\{0\}$. Let $\Omega$ be a bounded open subset of ${\mathbb{R}}^{n}$ of class $C^{m,\alpha}$. Let $\nu \equiv (\nu_{l})_{l=1,\dots,n}$ denote the external unit normal to $\partial\Omega$. Let $N_{2}$ denote the number of multi-indexes $\gamma\in {\mathbb{N}}^{n}$ with $|\gamma|\leq 2$. For each 
\begin{equation}
\label{introd0}
{\mathbf{a}}\equiv (a_{\gamma})_{|\gamma|\leq 2}\in {\mathbb{C}}^{N_{2}}\,, 
\end{equation}
we set 
\[
a^{(2)}\equiv (a_{lj} )_{l,j=1,\dots,n}\qquad
a^{(1)}\equiv (a_{j})_{j=1,\dots,n}\qquad
a\equiv a_{0}\,.
\]
with $a_{lj} \equiv 2^{-1}a_{e_{l}+e_{j}}$ for $j\neq l$, $a_{jj} \equiv
 a_{e_{j}+e_{j}}$,
and $a_{j}\equiv a_{e_{j}}$, where $\{e_{j}:\,j=1,\dots,n\}$  is the canonical basis of ${\mathbb{R}}^{n}$. We note that the matrix $a^{(2)}$ is symmetric. 
Then we assume that 
  ${\mathbf{a}}\in  {\mathbb{C}}^{N_{2}}$ satisfies the following ellipticity assumption
\begin{equation}
\label{ellip}
\inf_{
\xi\in {\mathbb{R}}^{n}, |\xi|=1
}{\mathrm{Re}}\,\left\{
 \sum_{|\gamma|=2}a_{\gamma}\xi^{\gamma}\right\} >0\,,
\end{equation}
and we consider  the case in which
\begin{equation}
\label{symr}
a_{lj} \in {\mathbb{R}}\qquad\forall  l,j=1,\dots,n\,.
\end{equation}
Then we introduce the operators
\begin{eqnarray*}
P[{\mathbf{a}},D]u&\equiv&\sum_{l,j=1}^{n}\partial_{x_{l}}(a_{lj}\partial_{x_{j}}u)
+
\sum_{l=1}^{n}a_{l}\partial_{x_{l}}u+au\,,
\\
B_{\Omega}^{*}v&\equiv&\sum_{l,j=1}^{n} \overline{a}_{jl}\nu_{l}\partial_{x_{j}}v
-\sum_{l=1}^{n}\nu_{l}\overline{a}_{l}v\,,
\end{eqnarray*}
for all $u,v\in C^{2}(\overline{\Omega})$, and a fundamental solution $S_{{\mathbf{a}} }$ of $P[{\mathbf{a}},D]$, and the   double layer potential 
\begin{eqnarray}
\label{introd3}
\lefteqn{
w[\partial\Omega       ,{\mathbf{a}},S_{{\mathbf{a}}}   ,\mu](x) \equiv 
\int_{\partial\Omega}\mu (y)\overline{B^{*}_{\Omega,y}}\left(S_{{\mathbf{a}}}(x-y)\right)
\,d\sigma_{y}
}
\\  \nonumber
&&
\qquad
=-\int_{\partial\Omega}\mu(y)\sum_{l,j=1}^{n} a_{jl}\nu_{l}(y)\frac{\partial S_{ {\mathbf{a}} } }{\partial x_{j}}(x-y)\,d\sigma_{y}
\\  \nonumber
&&
\qquad\quad
-\int_{\partial\Omega}\mu(y)\sum_{l=1}^{n}\nu_{l}(y)a_{l}
S_{ {\mathbf{a}} }(x-y)\,d\sigma_{y} \qquad\forall x\in {\mathbb{R}}^{n}\,,
\end{eqnarray}
where the density or moment $\mu$ is a function  from $\partial\Omega$ to ${\mathbb{C}}$. Here the subscript $y$ of $\overline{B^{*}_{\Omega,y}}$ means that we are taking $y$ as variable of the differential operator $\overline{B^{*}_{\Omega,y}}$. 
The role of the double layer potential in the solution of boundary value problems for the operator $P[{\mathbf{a}},D]$ is well known (cf.~\textit{e.g.}, 
G\"{u}nter~\cite{Gu67}, Kupradze,  Gegelia,  Basheleishvili and 
 Burchuladze~\cite{KuGeBaBu79}, Mikhlin \cite{Mik70}.) 
 
 The analysis of the continuity and compactness properties of the integral operator  associated to the double layer potential is  a classical topic. In particular, it has long been known that if $\mu$ is of class $C^{m,\alpha}$, then the restriction of the double layer potential to the sets
 \[
 \Omega^{+}\equiv \Omega\,,
 \qquad
 \Omega^{-}\equiv {\mathbb{R}}^{n}\setminus{\mathrm{cl}}\Omega\,,
 \]
can be extended to a function of $C^{m,\alpha}({\mathrm{cl}}\Omega^{+})$
and to a function of $C^{m,\alpha}_{{\mathrm{loc}} }({\mathrm{cl}}\Omega^{-})$, respectively (cf.~\textit{e.g.}, Miranda~\cite{Mi65},  
Wiegner~\cite{Wi93}, Dalla Riva~\cite{Da13}, Dalla Riva, Morais and Musolino~\cite{DaMoMu13}.)

In case $n=3$ and $\Omega$ is of class $C^{1,\alpha}$ and 
$S_{{\mathbf{a}}}$ is the fundamental solution of the Laplace operator,
it has long been known that  $w[\partial\Omega, {\mathbf{a}},S_{{\mathbf{a}}},\cdot]_{|\partial\Omega}$
is  a linear and compact operator in $C^{1,\alpha}(\partial\Omega)$ and is linear and continuous from $C^{0 }(\partial\Omega)$ to $C^{0,\alpha}(\partial\Omega)$ (cf.~
Schauder~\cite{Sc31}, \cite{Sc32}, Miranda~\cite{Mi65}.)

In case  $n=3$, $m\geq 1$ and $\Omega$ is of class $C^{m+1}$ and if  $P[{\mathbf{a}},D]$ is the  Laplace operator,  G\"{u}nter~\cite[Ch.~II, \S\ 21, Thm.~3]{Gu67} has 
proved that   $w[\partial\Omega ,{\mathbf{a}},S_{{\mathbf{a}}},\cdot]_{|\partial\Omega}$ is bounded from $C^{m-1,\alpha'}(\partial\Omega)$ to $C^{m,\alpha}(\partial\Omega)$ for $\alpha'\in]\alpha,1[$ and that accordigly it
 is  compact in $C^{m,\alpha}(\partial\Omega)$.

Fabes,  Jodeit and Rivi\`{e}re~\cite{FaJoRi78} have proved that if $\Omega$ is of class $C^{1}$
and  if $P[{\mathbf{a}},D]$ is the  Laplace operator, then $w[\partial\Omega, {\mathbf{a}},S_{{\mathbf{a}}},\cdot]_{|\partial\Omega}$ is compact in $L^{p}(\partial\Omega)$ for $p\in]1,+\infty[$. Later Hofmann,  M.~Mitrea and Taylor~\cite{HoMitTa10} have proved the same compactness result under more general conditions on $\partial\Omega$.

In case $n=2$ and $\Omega$ is of class $C^{2,\alpha}$ and if $P[{\mathbf{a}},D]$ is the  Laplace operator, Schippers~\cite{Schi82} has proved that $w[\partial\Omega, {\mathbf{a}},S_{{\mathbf{a}}},\cdot]_{|\partial\Omega}$ is continuous from $C^{0 }(\partial\Omega)$ to $C^{1,\alpha}(\partial\Omega)$.   

In case $n=3$  and $\Omega$ is of class $C^{2}$   and if  $P[{\mathbf{a}},D]$ is the  Helmholtz operator, Colton and Kress~\cite{CoKr83} have 
developed previous work of G\"{u}nter~\cite{Gu67} and Mikhlin~\cite{Mik70}  and proved that the operator $w[\partial\Omega ,{\mathbf{a}},S_{{\mathbf{a}}},\cdot]_{|\partial\Omega}$
is bounded from $C^{0,\alpha}(\partial\Omega)$ to $C^{1,\alpha}(\partial\Omega)$ and that accordigly it
 is  compact in $C^{1,\alpha}(\partial\Omega)$.

Wiegner~\cite{Wi93} has proved that if $\gamma\in {\mathbb{N}}^{n}$ has odd length and $\Omega$ is of class $C^{m,\alpha}$, then the operator with kernel 
$ (x-y)^{\gamma}  |x-y|^{-(n-1)- |\gamma|}  
$
is continuous from  $C^{m-1,\alpha}(\partial\Omega)$ to $C^{m-1,\alpha}({\mathrm{cl}}\Omega)$ (and a corresponding result for the exterior of $\Omega$). 

 In case  $n=3$, $m\geq 2$ and $\Omega$ is of class $C^{m,\alpha}$  and if  $P[{\mathbf{a}},D]$ is the  Helmholtz operator, Kirsch~\cite{Ki89} has proved that  the operator  $w[\partial\Omega ,{\mathbf{a}},S_{{\mathbf{a}}},\cdot]_{|\partial\Omega}$ is  bounded from $C^{m-1,\alpha}(\partial\Omega)$ to $C^{m,\alpha}(\partial\Omega)$ and that accordigly it
 is  compact in $C^{m,\alpha}(\partial\Omega)$. 

von Wahl~\cite{vo90} has considered the case of Sobolev spaces and has proved that if $\Omega$ is of class $C^{\infty}$ and if 
 $S_{{\mathbf{a}}}$ is the fundamental solution of the Laplace operator,  then  the double layer improves the regularity of one unit on the boundary.
 
 Then Heinemann~\cite{He92} has developed the ideas of  von Wahl in the frame of Schauder spaces and has proved that
 if $\Omega$ is of class $C^{m+5}$ and if 
 $S_{{\mathbf{a}}}$ is the fundamental solution of the Laplace operator,  then  the double layer improves the regularity of one unit on the boundary, \textit{i.e.}, 
 $w[\partial\Omega, {\mathbf{a}},S_{{\mathbf{a}}},\cdot]_{|\partial\Omega}$ is linear and continuous from $C^{m,\alpha}(\partial\Omega)$ to $C^{m+1,\alpha}(\partial\Omega)$.

 Maz'ya and  Shaposhnikova~\cite{MaSh05}
 have proved that $w[\partial\Omega ,{\mathbf{a}},S_{{\mathbf{a}}},\cdot]_{|\partial\Omega}$ is continuous in fractional Sobolev spaces under sharp   regularity assumptions on the boundary and if  $P[{\mathbf{a}},D]$ is the  Laplace operator. 

Mitrea~\cite{Mit14}  has proved that the double layer of second order equations and systems is compact in  $C^{0,\beta}(\partial\Omega)$   for $\beta\in]0, \alpha[$ and bounded  in $C^{0,\alpha}(\partial\Omega)$
  under the assumption that $\Omega$ is of class $C^{1,\alpha}$. Then by exploiting a formula for the tangential derivatives such results have been extended to  compactness  and boundedness results in  $C^{1,\beta}(\partial\Omega)$ and $C^{1,\alpha}(\partial\Omega)$, respectively.

Mitrea,  Mitrea  and Verdera~\cite{MitMitVe14} have proved that if $q$ is a homogeneous polynomial of odd order, then the operator with kernel
$ q(x-y) |x-y|^{-(n-1)-{\mathrm{deg}} (q)}$  maps $C^{0,\alpha}(\partial\Omega)$ to $C^{1,\alpha}({\mathrm{cl}}\Omega)$. 

 In this paper we are interested into the regularizing properties of the operator $w[\partial\Omega, {\mathbf{a}},S_{{\mathbf{a}}},\cdot]_{|\partial\Omega}$
  in Schauder spaces under the assumption that $\Omega$ is   of class $C^{m,\alpha}$.  We prove our statements by exploiting 
tangential derivatives and an inductive argument to reduce the problem to the case of the action of  $w[\partial\Omega, {\mathbf{a}},S_{{\mathbf{a}}},\cdot]_{|\partial\Omega}$ on  $C^{0,\alpha}(\partial\Omega)$
 instead of flattening the boundary with   parametrization functions as done by other authors. We mention that the idea of exploiting an inductive argument together with a  formula for the tangential gradient in order to prove continuity and compactness properties of the double layer potential  has been exploited by 
 Kirsch~\cite[Thm.~3.2]{Ki89} in case $n=3$ and $P[{\mathbf{a}},D]$  equals the Helmholtz operator and  $S_{{\mathbf{a}}}$ is the fundamental solution satisfying the radiation condition.  The tangential derivatives of $f\in C^{1}(\partial\Omega)$ are defined by the equality 
 \[
M_{lr}[f]\equiv \nu_{l}\frac{\partial\tilde{f}}{\partial x_{r}}-
\nu_{r}\frac{\partial\tilde{f}}{\partial x_{l}}\qquad {\mathrm{on}}\ \partial\Omega\,,
\]
for all   $l,r\in\{1,\dots,n\}$. Here $\tilde{f} $ denotes an extension of $f$ to an open neighbourhood of $\partial\Omega$, and one can easily verify that $M_{lr}[f]$ is independent of the specific choice of the extension $\tilde{f}$ of $f$. Then we prove an explicit  formula for
\begin{equation}
\label{introd4}
M_{lr}[  w[\partial\Omega ,{\mathbf{a}},S_{{\mathbf{a}}},\mu] ](x)
-
w[\partial\Omega ,{\mathbf{a}},S_{{\mathbf{a}}},M_{lr}[\mu]](x)\qquad\forall x\in \partial\Omega\,,
\end{equation}
for all $\mu\in C^{1}(\partial\Omega)$  and $l,r\in\{1,\dots,n\}$ (see formula (\ref{wtg1}).)  

We note that   G\"{u}nter~\cite[Ch.~II, \S\ 10, (42)]{Gu67} contains a formula for the partial derivatives of the double layer with respect to the variables in ${\mathbb{R}}^{n}$ in case $n=3$ and $P[{\mathbf{a}},D]$  equals the  Laplace operator (see (\ref{dlay0}) in case of the Laplace operator.)    A similar formula can be found in 
Kupradze,  Gegelia,  Basheleishvili and 
 Burchuladze~\cite[Ch.~V, \S\ 6, (6.11)]{KuGeBaBu79} for the elastic double layer potential in case $n=3$.      Schwab and Wendland~\cite{ScWe99} have proved that   the difference in (\ref{introd4}) can be written in terms of  pseudodifferential operators of order $-1$.   Dindo\v{s} and   Mitrea  have proved a number of properties of the double layer potential. In particular, \cite[Prop.~3.2]{DiMi04} proves the existence of integral operators such that  the gradient of the double layer potential corresponding to the Stokes system  can be written as a sum of such integral operators applied to the gradient of the moment of the double layer.    Duduchava, Mitrea, and Mitrea~\cite{DuMiMi06} analyze  various properties of the tangential deriatives. Duduchava~\cite{Du10} investigates partial differential equations on hypersurfaces and Bessel potential operators. In   particular   \cite[point B of the proof of Lem.~2.1]{Du10} analyzes the properties of a commutator of a Bessel potential operator and of a tangential derivative.
 Hofmann, Mitrea and Taylor~\cite[(6.2.6)]{HoMitTa10} prove a general formula for the tangential derivatives of the double layer potential corresponding to second order elliptic  {\em homogeneous} equations and systems in explicit terms.  
 
The formula  (\ref{wtg1}) we compute here extends a formula of  \cite{La08} for the Laplace operator,  which has been computed with arguments akin to those of G\"{u}nter~\cite[Ch.~II, \S\ 10, (42)]{Gu67}, and a formula of \cite{Do12}  for the Helmholtz operator, and can be considered as a variant of the formula of Hofmann, Mitrea and Taylor~\cite[(6.2.6)]{HoMitTa10} for the second order    {\em nonhomogeneous} elliptic differential operator $P[{\mathbf{a}},D]$.

Formula  (\ref{wtg1}) involves auxiliary operators, which we analyze in section \ref{aio}. We have based our analysis of the auxiliary operators involved in formula (\ref{wtg1}) on the introduction of boundary norms for weakly singular kernels and on a result of joint continuity of weakly singular integrals both on the kernel of the integral and on the functional variable of the corresponding integral operator (see section \ref{bonoke}.) For fixed choices of the kernel and for some choices of the parameters, such lemmas are known (cf.~\textit{e.g.}, Kirsch and Hettlich~\cite[Thm.~3.17, p.~121]{KiHe15}.) The authors believe that the methods of section \ref{bonoke} may be applied  to simplify also  the exposition of other classical proofs of properties of layer potentials. 

By exploiting formula (\ref{wtg1}), we can prove that $w[\partial\Omega , {\mathbf{a}}, S_{{\mathbf{a}}},\cdot,]_{|\partial\Omega}$ induces a linear and continuous operator from $C^{m}(\partial\Omega)$ to the generalized Schauder space
$C^{m,\omega_{\alpha}}(\partial\Omega)$ of functions with $m$-th order derivatives which satisfy a generalized $\omega_{\alpha}$-H\"{o}lder condition with
\[
\omega_{\alpha}(r)\sim r^{\alpha}|\ln r| 
\qquad{\mathrm{as}}\ r\to 0,
\]
and that $w[\partial\Omega ,{\mathbf{a}},S_{{\mathbf{a}}},\cdot]_{|\partial\Omega}$ induces a linear and continuous operator from $C^{m,\beta}(\partial\Omega)$
to $C^{m,\alpha}(\partial\Omega)$ for all $\beta\in]0,\alpha]$. In particular, the double layer potential has a regularizing effect on the boundary if $\Omega$ is of class $C^{m,\alpha}$.
As a consequence of our result,  $w[\partial\Omega ,{\mathbf{a}},S_{{\mathbf{a}}},\cdot]_{|\partial\Omega}$ induces a compact operator from $C^{m }(\partial\Omega)$ to itself, and from $C^{m,\omega_{\alpha}(\cdot)}(\partial\Omega)$ to itself, and from 
$C^{m,\alpha}(\partial\Omega)$ to itself when $\Omega$ is of class $C^{m,\alpha}$.

 \section{Notation}\label{nota}

We  denote the norm on 
a   normed space ${\mathcal X}$ by $\|\cdot\|_{{\mathcal X}}$. Let 
${\mathcal X}$ and ${\mathcal Y}$ be normed spaces. We endow the  
space ${\mathcal X}\times {\mathcal Y}$ with the norm defined by 
$\|(x,y)\|_{{\mathcal X}\times {\mathcal Y}}\equiv \|x\|_{{\mathcal X}}+
\|y\|_{{\mathcal Y}}$ for all $(x,y)\in  {\mathcal X}\times {\mathcal 
Y}$, while we use the Euclidean norm for ${\mathbb{R}}^{n}$.
 For 
standard definitions of Calculus in normed spaces, we refer to 
Deimling~\cite{De85}.  If $A$ is a matrix with real or complex entries, then $A^{t}$ denotes the transpose matrix of $A$. 
The set $M_{n}({\mathbb{R}})$ denotes the set of $n\times n$ matrices with real entries.  Let 
${\mathbb{D}}\subseteq {\mathbb {R}}^{n}$. Then $\mathrm{cl}{\mathbb{D}}$ 
denotes the 
closure of ${\mathbb{D}}$, and $\partial{\mathbb{D}}$ denotes the boundary of ${\mathbb{D}}$, 
and ${\mathrm{diam}} ({\mathbb{D}})$ denotes the diameter of   ${\mathbb{D}}$. The symbol
$| \cdot|$ denotes the Euclidean modulus   in
${\mathbb{R}}^{n}$ or in ${\mathbb{C}}$. For all $R\in]0,+\infty[$, $ x\in{\mathbb{R}}^{n}$, 
$x_{j}$ denotes the $j$-th coordinate of $x$, and  
 ${\mathbb{B}}_{n}( x,R)$ denotes the ball $\{
y\in{\mathbb{R}}^{n}:\, | x- y|<R\}$. 
Let $\Omega$ be an open 
subset of ${\mathbb{R}}^{n}$. The space of $m$ times continuously 
differentiable complex-valued functions on $\Omega$ is denoted by 
$C^{m}(\Omega,{\mathbb{C}})$, or more simply by $C^{m}(\Omega)$. 
Let $s\in {\mathbb{N}}\setminus\{0\}$, $f\in \left(C^{m}(\Omega)\right)^{s} $. Then   $Df$ denotes the Jacobian matrix of $f$. 
Let  $\eta\equiv
(\eta_{1},\dots ,\eta_{n})\in{\mathbb{N}}^{n}$, $|\eta |\equiv
\eta_{1}+\dots +\eta_{n}  $. Then $D^{\eta} f$ denotes
$\frac{\partial^{|\eta|}f}{\partial
x_{1}^{\eta_{1}}\dots\partial x_{n}^{\eta_{n}}}$.    The
subspace of $C^{m}(\Omega )$ of those functions $f$ whose derivatives $D^{\eta }f$ of
order $|\eta |\leq m$ can be extended with continuity to 
$\mathrm{cl}\Omega$  is  denoted $C^{m}(
\mathrm{cl}\Omega )$.\par 

The
subspace of $C^{m}(\mathrm{cl}\,\Omega ) $  whose derivatives up to order $m$ are bounded is denoted 
$C^{m}_{b}(\mathrm{cl}\,\Omega ) $. Then $C^{m}_{b}(\mathrm{cl}\,\Omega )$ 
endowed with the norm
$\|f\|_{C^{m}_{b}(\mathrm{cl}\,\Omega )}\equiv \sum
_{|\eta |\leq m}\sup_{\mathrm{cl}\,\Omega}|D^\eta f|$ is a Banach space.\par 

Now let $\omega$ be a function of $]0,+\infty[$ to itself such that
\begin{equation}
\label{om}
\omega\ {\mathrm{is\   increasing\  and}}\ \lim_{r\to 0^{+}}\omega(r)=0\,,
\end{equation}
and that
\begin{equation}
\label{om1}
\sup_{r\in ]0,1[} \omega^{-1}(r)r<\infty\,.
\end{equation}
If $f$ is a function from a subset ${\mathbb{D}}$ of ${\mathbb{R}}^{n}$   to ${\mathbb{C}}$, we set
\[
|f:{\mathbb{D}}|_{\omega (\cdot)
}
\equiv
\sup\left\{
\frac{|f( x )-f( y)|}{\omega(| x- y|)
}: x, y\in {\mathbb{D}} ,  x\neq
 y\right\}\,.
\]
If $|f:{\mathbb{D}}|_{\omega(\cdot)}<\infty$, we say that the function $f$ is $\omega(\cdot)$-H\"{o}lder continuous. Sometimes, we simply write $|f|_{\omega(\cdot)}$ instead of $|f:{\mathbb{D}}|_{\omega(\cdot)}$.
If $\omega(r)=r$, and if $|f:{\mathbb{D}}|_{\omega(\cdot)}<\infty$, then we say that    $f$ is Lipschitz continuous and we set ${\mathrm{Lip}}(f)\equiv 
|f:{\mathbb{D}}|_{\omega(\cdot)}$.  The
subspace of $C^{0}({\mathbb{D}} ) $  whose
functions  are
$\omega(\cdot)$-H\"{o}lder continuous    is denoted $C^{0,\omega(\cdot)} ({\mathbb{D}})$, and the subspace of $C^{0}({\mathbb{D}} ) $  whose
functions  are Lipschitz continuous is denoted ${\mathrm{Lip}}({\mathbb{D}})$.

Let $\Omega$ be an open subset of ${\mathbb{R}}^{n}$. The
subspace of $C^{m}(\mathrm{cl}\,\Omega ) $  whose
functions have $m$-th order derivatives that are
$\omega(\cdot)$-H\"{o}lder continuous    is denoted $C^{m,\omega(\cdot)} (\mathrm{cl}\,\Omega )$. Then we set
\[
C^{m,\omega(\cdot)}_{b} (\mathrm{cl}\,\Omega )
\equiv C^{m,\omega(\cdot)} (\mathrm{cl}\,\Omega )
\cap
C^{m}_{b}(\mathrm{cl}\,\Omega )\,.
\]
The space $C^{m,\omega(\cdot)}_{b}({\mathrm{cl}}\,
\Omega )$, equipped with its usual norm 
\[
\|f\|_{C^{m,\omega(\cdot)}_{b}(\mathrm{cl}\,\Omega )}=
\|f\|_{C^{m}_{b}(\mathrm{cl}\,\Omega )}+\sum_{|\eta |=m}
|D^{\eta }f:\,\Omega|_{\omega(\cdot)}\,, 
\]
is well-known to be a Banach space. \par

Obviously, $C^{m,\omega(\cdot)}_{b}(\mathrm{cl}\,\Omega )=C^{m,\omega(\cdot)}(\mathrm{cl}\,\Omega )$ if $\Omega $ is bounded (and in this case, we shall always drop the subscript $b$.) The subspace of $C^{m}(\mathrm{cl}\, \Omega  ) $ of those functions $f$ such that $f_{|{\mathrm{cl}}(\Omega \cap{\mathbb{B}}_{n}(0,R))}\in
C^{m,\omega(\cdot)}({\mathrm{cl}}(\Omega \cap{\mathbb{B}}_{n}(0,R)))$ for all $R\in]0,+\infty[$ is denoted $C^{m,\omega(\cdot)}_{{\mathrm{loc}}}(\mathrm{cl}\,\Omega ) $.
 Clearly, 
$C^{m,\omega(\cdot)}_{{\mathrm{loc}}}(\mathrm{cl}\,\Omega ) = C^{m,\omega(\cdot)}(\mathrm{cl}\, \Omega  )$ if $\Omega $ is bounded.

Particularly important is the case in which $\omega(\cdot)$ is the function 
$r^{\alpha}$ for some fixed $\alpha\in]0,1]$. In this case, we simply write $|\cdot:{\mathrm{cl}}\,\Omega|_{\alpha}$ instead of
$|\cdot:{\mathrm{cl}}\,\Omega|_{r^{\alpha}}$, and
$C^{m,\alpha}(\mathrm{cl}\,\Omega )$
instead of
$C^{m,r^{\alpha}}(\mathrm{cl}\,\Omega )$, and $C^{m,\alpha}_{b}(\mathrm{cl}\,\Omega )$
instead of
$C^{m,r^{\alpha}}_{b}(\mathrm{cl}\,\Omega )$.
We observe that property (\ref{om1}) implies that
\[
C^{m,1}_{b}(\mathrm{cl}\,\Omega )\subseteq
C^{m,\omega(\cdot)}_{b}(\mathrm{cl}\,\Omega )\,.
\]
For the definition of a bounded open Lipschitz subset of ${\mathbb{R}}^{n}$, we refer for example to Ne\v{c}as~\cite[\S 1.3]{Ne67}.  Let $m\in {\mathbb{N}}\setminus\{0\}$. 
We say that a bounded open subset $\Omega $ of ${\mathbb{R}}^{n}$ is of class 
 $C^{m,\alpha}$, if  for every 
$P\in \partial \Omega $, there exist  an open neighborhood $W$ of $P$ in ${\mathbb{R}}^{n}$, 
and a diffeomorphism  $\psi\in C^{m,\alpha}
\left({\mathrm{cl}}\,{\mathbb{B}}_{n},
{\mathbb{R}}^{n}\right)$ of ${\mathbb{B}}_{n}\equiv\{x\in {\mathbb{R}}^{n}:\,|x|<1\}$ onto
 $W$ such that $\psi (0)=P$, $\psi(
 \{x\in {\mathbb{B}}_{n}:\, x_{n}=0\}
 )=W\cap\partial \Omega $,  $\psi(
 \{x\in {\mathbb{B}}_{n}:\, x_{n}<0\}
 )=W\cap  \Omega $ ($\psi $ is said to be a parametrization of $\partial \Omega $ around $P$.)  
 Now let $\Omega $ be bounded and of class $C^{m,\alpha}$.  By compactness of $\partial \Omega $ and by definition of set of class $C^{m,\alpha}$,   there exist 
$P_{1}$,\dots,$P_{r}\in \partial \Omega $, and parametrizations 
$\{\psi_{i}\}_{i=1,\dots,r}$, with 
$\psi_{i}\in C^{m,\alpha}\left({\mathrm{cl}}\,{\mathbb{B}}_{n},
{\mathbb{R}}^{n}\right)$  such that 
$\cup_{i=1}^{r}\psi_{i}( \{x\in {\mathbb{B}}_{n}:\, x_{n}=0\}
)=\partial \Omega $.  Let $h\in\{1,\dots,m\}$. Let $\omega$ be as in (\ref{om}), (\ref{om1}). Let
\begin{equation}
\label{om3}
\sup_{r\in ]0,1[} \omega^{-1}(r)r^{\alpha}<\infty\,.
\end{equation} 
We denote by $C^{h,\omega(\cdot)}\left(\partial \Omega \right)$ the linear space of functions $f$ of $\partial \Omega $ to 
${\mathbb{R}}$ such that 
$f\circ\psi_{i}(\cdot,0)\in C^{h,\omega(\cdot)}
\left({\mathrm{cl}}\,{\mathbb{B}}_{n-1}\right)$ for all $i=1,\dots,r$, and we set
\[
\|f\|_{C^{h,\omega(\cdot)}\left(\partial \Omega \right)}
\equiv
\sup_{i=1,\dots,r}\|f\circ\psi_{i}(\cdot,0)\|_{C^{h,\omega(\cdot)}
\left({\mathrm{cl}}\,{\mathbb{B}}_{n-1}\right)}\qquad
\forall f\in C^{h,\omega(\cdot)}\left(\partial \Omega \right). 
\]
It is well known that by choosing a different finite family of 
parametrizations as $\{\psi_{i}\}_{i=1,\dots,r}$, we would obtain an 
equivalent norm. In case $\omega(\cdot)$ is the function 
$r^{\alpha}$, we have the spaces $C^{h,\alpha}\left(\partial \Omega \right)$.

It is known that $(C^{h,\omega(\cdot)}\left(\partial \Omega\right),\|\cdot \|_{C^{h,\omega(\cdot)}\left(\partial \Omega\right)})$ is complete. Moreover condition (\ref{om3}) implies that  the restriction operator is linear and continuous from 
$C^{h,\omega(\cdot)}\left({\mathrm{cl}} \Omega\right)$ to
$C^{h,\omega(\cdot)}\left(\partial \Omega\right)$.

\par
We denote by $d\sigma$ the area element of a manifold imbedded in ${\mathbb{R}}^{n}$. We retain the standard notation for the Lebesgue spaces. \par

\begin{rem}
\label{h=m}
Let $m\in{\mathbb{N}}\setminus\{0\}$, $\alpha\in]0,1[$. Let $\Omega$ be a  bounded open subset of ${\mathbb{R}}^{n}$ of class $C^{m,\alpha}$.  

Let $\omega$ be as in (\ref{om}), (\ref{om1}). If $h\in\{1,\dots,m\}$, $h<m$, then $m-1\geq 1$ and $\Omega$ is of class $C^{m-1,1}$ and  condition (\ref{om1}) implies the validity of condition  (\ref{om3})  with $\alpha$ replaced by $1$, and thus we can consider the space $C^{h,\omega (\cdot)}(\partial\Omega)$ even if we do not assume condition (\ref{om3}). If instead  $h=m$, the definition we gave requires (\ref{om3}). 
\end{rem}

\begin{rem}
\label{om4}
Let $\omega$ be as in (\ref{om}). 
Let ${\mathbb{D}}$ be a   subset of ${\mathbb{R}}^{n}$. Let $f$ be a bounded function from $ {\mathbb{D}}$ to ${\mathbb{C}}$, $a\in]0,+\infty[$.  Then,
\[
\sup_{x,y\in {\mathbb{D}},\ |x-y|\geq a}\frac{|f(x)-f(y)|}{\omega(|x-y|)}
\leq \frac{2}{\omega(a)} \sup_{{\mathbb{D}}}|f|\,.
\]
\end{rem}
Thus the difficulty of estimating the H\"{o}lder quotient $\frac{|f(x)-f(y)|}{\omega(|x-y|)}$ of a bounded function $f$ lies entirely in case $0<|x-y|<a$.  
Then we have the following well known extension result. For a proof, we refer to Troianiello \cite[Thm.~1.3, Lem.~1.5]{Tr87}. 

\begin{lem}
\label{extsch}
Let $m\in {\mathbb{N}}\setminus\{0\}$, $\alpha\in ]0,1[$, $j\in\{0,\dots,m\}$. Let $\Omega$ be a bounded open subset of ${\mathbb{R}}^{n}$ of class $C^{m,\alpha}$. Let $R\in]0,+\infty[$ be such that ${\mathrm{cl}}\Omega\subseteq{\mathbb{B}}_{n}(0,R)$. 
Then there exists a linear and continuous extension operator ` $\tilde{\ }$ ' of $C^{j,\alpha}(\partial\Omega)$ to $C^{j,\alpha}({\mathrm{cl}}
{\mathbb{B}}_{n}(0,R))$, which takes $\mu \in C^{j,\alpha}(\partial\Omega)$ to a map $\tilde{\mu}\in C^{j,\alpha}({\mathrm{cl}}
{\mathbb{B}}_{n}(0,R))$ such that $\tilde{\mu}_{|\partial\Omega}=\mu$ and such that the support of $\mu$ is compact and contained in ${\mathbb{B}}_{n}(0,R)$. 
The same statement holds by replacing $C^{m,\alpha}$ by $C^{m}$ and $C^{j,\alpha}$ by $C^{j}$.
\end{lem}
Let $\Omega$ be a bounded open subset of ${\mathbb{R}}^{n}$ of class $C^{1}$.
The tangential gradient  $D_{\partial\Omega} f$ of $f\in C^{1}(\partial\Omega)$ is defined as
\[
D_{\partial\Omega}f\equiv D\tilde{f}-(\nu\cdot D\tilde{f})\nu
\qquad{\mathrm{on}}\ \partial\Omega\,, 
\]
where $\tilde{f}$ is an extension  of $f$ of class $C^{1}$ in an open neighborhood of $\partial\Omega$, 
and we have 
\[
\frac{\partial\tilde{f}}{\partial x_{r}}-
(\nu\cdot D\tilde{f})\nu_{r}=
\sum_{l=1}^{n}M_{lr}[f]\nu_{l}
\qquad{\mathrm{on}}\ \partial\Omega\,,
\]
for all $r\in\{1,\dots,n\}$. If ${\mathbf{a}}$ is as in (\ref{introd0}), (\ref{ellip}), then we also set
\[
D_{ {\mathbf{a}} }f\equiv (D_{ {\mathbf{a}} ,r}f)_{r=1,\dots,n}\equiv
D\tilde{f}
-
\frac{D\tilde{f}a^{(2)}\nu }{\nu^{t}
a^{(2)}\nu 
}\nu \qquad{\mathrm{on}}\ \partial\Omega\,.
\]
Since
\begin{equation}
\label{dam}
D_{ {\mathbf{a}} ,r}f=\frac{\partial\tilde{f}}{\partial x_{r}}-
\frac{ D\tilde{f} a^{(2)} \nu  }{ \nu^{t}
a^{(2)}\nu }\nu_{r}
=\sum_{l=1}^{r}M_{lr}[f]\left(
\frac{\sum_{h=1}^{n}a_{lh}\nu_{h}}
{
\nu^{t}a^{(2)}\nu
}
\right)  \qquad{\mathrm{on}}\ \partial\Omega\,,
\end{equation}
for all $r\in\{1,\dots,n\}$,  $D_{ {\mathbf{a}} }f$ is independent of the specific choice of the extension $\tilde{f}$ of $f$. 
We also need the following well known consequence of the Divergence Theorem.
\begin{lem}
\label{gagre}
Let $\Omega$ be a bounded open subset of ${\mathbb{R}}^{n}$ of class $C^{1}$. If $\varphi$, $\psi\in C^{1}(\partial\Omega)$, then
\[
\int_{\partial\Omega}M_{lj}[\varphi]\psi\,d\sigma=-
\int_{\partial\Omega}\varphi M_{lj}[\psi] \,d\sigma
\]
for all $l,j\in \{1,\dots,n\}$.
\end{lem}

Next we introduce the following auxiliary Lemmas, whose proof is based on the definition of norm in a Schauder space.

\begin{lem}
\label{tanco}
Let $m\in {\mathbb{N}}\setminus\{0\}$, $\alpha\in]0,1]$. Let $\omega$ be as in (\ref{om}), (\ref{om1}), (\ref{om3}). Let $\Omega$ be a bounded open connected subset of ${\mathbb{R}}^{n}$ of class $C^{m,\alpha}$. Then the following statements hold.
\begin{enumerate}
\item[(i)] A function $f\in C^{1}(\partial\Omega)$ belongs to $C^{m,\omega(\cdot)}(\partial\Omega)$ if and only if $M_{lr}[f]\in C^{m-1,\omega(\cdot)}(\partial\Omega)$ for all $l,r\in\{1,\dots,n\}$.
\item[(ii)] The norm $\|\cdot\|_{C^{m,\omega(\cdot)}(\partial\Omega)}$ is equivalent to the norm on $C^{m,\omega(\cdot)}(\partial\Omega)$ defined by
\[
\|f\|_{C^{0}(\partial\Omega)}+
\sum_{l,r=1}^{n}\|M_{lr}[f]\|_{C^{m-1,\omega(\cdot)}(\partial\Omega)}\qquad\forall
f\in C^{m,\omega(\cdot)}(\partial\Omega)\,.
\]
\end{enumerate}
\end{lem}
Then we have the following (see also Remark \ref{h=m}.)
\begin{lem}
\label{thgs}
Let $m\in {\mathbb{N}}\setminus\{0\}$, $\alpha\in ]0,1]$. Let $\Omega$ be a bounded open connected subset of ${\mathbb{R}}^{n}$ of class $C^{m,\alpha}$. Let $h\in\{1,\dots,m\}$. Then the following statements hold.
\begin{enumerate}
\item[(i)]  Let  $h<m$. Let $\omega$  be as in (\ref{om}), (\ref{om1}). Then $M_{lj}$ is linear and continuous from
$C^{h,\omega(\cdot)}(\partial\Omega)$ to $C^{h-1,\omega(\cdot)}(\partial\Omega)$ for all $l,j\in\{1,\dots,n\}$. 
If we further assume that $\omega$ satisfies condition (\ref{om3}), then the same statement holds 
also for $h=m$.
\item[(ii)]  Let  $h<m$. Let $\omega$  be as in (\ref{om}), (\ref{om1}). Let ${\mathbf{a}}$ be as in (\ref{introd0}), (\ref{ellip}). Then the function from 
$C^{h,\omega(\cdot)}(\partial\Omega)$ to $C^{h-1,\omega(\cdot)}(\partial\Omega, {\mathbb{R}}^{n})$, which takes $f$ to $D_{ {\mathbf{a}} }f$
is linear and continuous. If we further assume that $\omega$ satisfies condition (\ref{om3}), then the same statement holds 
also for $h=m$.
\item[(iii)] Let $h<m$. 
Let $\omega$  be as in (\ref{om}), (\ref{om1}). Then the space $C^{h,\omega(\cdot)}(\partial\Omega)$ is continuously imbedded into 
$C^{h-1,1}(\partial\Omega)$. If we further assume that $\omega$ satisfies condition (\ref{om3}), then the same statement holds 
also for $h=m$.
\item[(iv)]  Let $h<m$.  Let $\psi_{1}$, $\psi_{2}$  be as in (\ref{om}), (\ref{om1}). Let  condition  
\[
\sup_{r\in ]0,1[} \psi_{2}^{-1}(r)\psi_{1}(r)<\infty
\]
 hold.
Then $C^{h,\psi_{1}(\cdot)}(\partial\Omega)$ is continuously imbedded into
$C^{h,\psi_{2}(\cdot)}(\partial\Omega)$.
If we further assume that $\psi_{j}$ satisfies condition (\ref{om3}) for $j\in\{1,2\}$, then the same statement holds also for $h=m$.
\item[(v)]  Let   $h<m$.  Let $\psi_{1}$, $\psi_{2}$, $\psi_{3}$ be as in (\ref{om}), (\ref{om1}). Let conditions
$
\sup_{j=1,2}\sup_{r\in ]0,1[} \psi_{j}(r)\psi_{3}^{-1}(r)<\infty$ 
hold. Then the pointwise product is bilinear and continuous from
$C^{h,\psi_{1}(\cdot)}(\partial\Omega)\times C^{h,\psi_{2}(\cdot)}(\partial\Omega)$ to $C^{h,\psi_{3}(\cdot)}(\partial\Omega)$.
If we further assume that $\psi_{j}$ satisfies condition (\ref{om3}) for $j\in\{1, 2,3\}$,
then the same statement holds also for $h=m$.
\end{enumerate}
\end{lem}
 
\begin{lem}
Let $\Omega$ be a bounded open Lipschitz subset of ${\mathbb{R}}^{n}$. 
Let $\psi_{1}$, $\psi_{2}$, $\psi_{3}$ be as in (\ref{om}), (\ref{om1}). Let conditions
$
\sup_{j=1,2}\sup_{r\in ]0,1[} \psi_{j}(r)\psi_{3}^{-1}(r)<\infty$ 
hold. Then the pointwise product is bilinear and continuous from the space 
$C^{0,\psi_{1}(\cdot)}(\partial\Omega)\times C^{0,\psi_{2}(\cdot)}(\partial\Omega)$ to $C^{0,\psi_{3}(\cdot)}(\partial\Omega)$.
\end{lem}

\section{Preliminary inequalities}
\label{preineq}
We first introduce the following elementary lemma on matrices.
\begin{lem}
\label{t-1}
Let $\Lambda\in M_{n}({\mathbb{R}})$ be invertible. Let $|\Lambda|\equiv
\sup_{|x|=1}|\Lambda x|$. Then the following statements hold. 
\begin{enumerate}
\item[(i)] Let $\tau_{\Lambda}\equiv\max\{|\Lambda|,|\Lambda^{-1}|\}$. Then 
\[
\tau_{\Lambda}^{-1}|x|\leq |\Lambda x|\leq\tau_{\Lambda}|x|\qquad\forall x\in {\mathbb{R}}^{n}\,.
\]
\item[(ii)] Let $r\in ]0,+\infty[$. Then
\[
|\Lambda^{-1}x|^{-r}\leq |\Lambda|^{r}|x|^{-r}\qquad\forall x\in {\mathbb{R}}^{n}\setminus\{0\}\,.
\]
\end{enumerate}
\end{lem}
{\bf Proof.} Statement (i) is well known. We now consider statement (ii). Let $x\in {\mathbb{R}}^{n}\setminus\{0\}$. Then we have
\[
|x|=|\Lambda(\Lambda^{-1}x)|\leq |\Lambda|\,|\Lambda^{-1}x|\,.
\]
Hence, $|\Lambda^{-1}x|\geq |\Lambda|^{-1}|x|$ and the statement follows. \hfill  $\Box$ 

\vspace{\baselineskip}

Then we introduce  the following elementary lemma, which collects either
known inequalities or variants of known inequalities, which we need in the sequel.
\begin{lem}
\label{rec}
Let $\gamma\in {\mathbb{R}}$. Let $\Lambda\in M_{n}({\mathbb{R}})$ be invertible. The following statements hold.
\begin{enumerate}
\item[(i)]
\begin{eqnarray*}
&&
\frac{1}{2 }|x'-y|\leq |x''-y|\leq 2|x'-y|\,,
\\ \nonumber
&&
\frac{1}{2\tau_{\Lambda}^{2} }|\Lambda x'-\Lambda y|\leq |\Lambda x''-\Lambda y|\leq 2\tau_{\Lambda }^{2}  |\Lambda x'-\Lambda y| \,,
\end{eqnarray*}
for all  $x',x''\in {\mathbb{R}}^{n}$, $x'\neq x''$, $y\in {\mathbb{R}}^{n}
\setminus {\mathbb{B}}_{n}(x',2|x'-x''|)$.
\item[(ii)] 
\begin{eqnarray*}
\lefteqn{
|x'-y|^{\gamma}\leq 2^{|\gamma|}|x''-y|^{\gamma}\,,
\quad
|x''-y|^{\gamma}\leq 2^{|\gamma|}|x'-y|^{\gamma}
\,,
}
\\ \nonumber
\lefteqn{
|\Lambda x'-\Lambda y|^{\gamma}\leq (2\tau_{\Lambda }^{2})^{|\gamma|}|\Lambda x''-\Lambda y|^{\gamma}\,,
}
\\ \nonumber
&&\qquad\qquad\qquad\qquad\qquad\qquad
|\Lambda x''-\Lambda y|^{\gamma}\leq (2\tau_{\Lambda }^{2})^{|\gamma|}|\Lambda x'-\Lambda y|^{\gamma}
\,,
\end{eqnarray*}
for all  $x',x''\in {\mathbb{R}}^{n}$, $x'\neq x''$, $y\in {\mathbb{R}}^{n}
\setminus {\mathbb{B}}_{n}(x',2|x'-x''|)$.
\item[(iii)] 
\[
|
|x'-y|^{\gamma}-|x''-y|^{\gamma}
|
\leq  (2^{|\gamma|}-1)|x'-y|^{\gamma}
\quad\forall y\in {\mathbb{R}}^{n}
\setminus {\mathbb{B}}_{n}(x',2|x'-x''|)\,,
\]
for all $x',x''\in {\mathbb{R}}^{n}$, $x'\neq x''$.
\item[(iv)] There exist  $m_{\gamma}$, $m_{\gamma}(\Lambda )\in ]0,+\infty[$ such that
\begin{eqnarray*}
&&
|
|x'-y|^{\gamma}-|x''-y|^{\gamma}
|
\leq m_{\gamma}|x'-x''| \, |x'-y|^{\gamma-1}
\\ \nonumber
&&
|
|\Lambda x'-\Lambda y|^{\gamma}-|\Lambda x''-\Lambda y|^{\gamma}
|
\leq m_{\gamma}(\Lambda )|\Lambda x'-\Lambda x''| \, |\Lambda x'-\Lambda y|^{\gamma-1}
\end{eqnarray*}
 for all  $x',x''\in {\mathbb{R}}^{n}$, $x'\neq x''$, $y\in {\mathbb{R}}^{n}
\setminus {\mathbb{B}}_{n}(x',2|x'-x''|)$.
\item[(v)]  
\[
|
\ln |x'-y| -\ln |x''-y| 
|
\leq 2|x'-x''| \, |x'-y|^{-1}
\ \
\forall y\in {\mathbb{R}}^{n}
\setminus {\mathbb{B}}_{n}(x',2|x'-x''|)\,,
\]
for all $x',x''\in {\mathbb{R}}^{n}$, $x'\neq x''$.
\end{enumerate}
\end{lem}
{\bf Proof.}   The first two inequalities of statement (i) follows by the triangular inequality. Then we have
\[
|\Lambda x'-\Lambda y|\leq \tau_{\Lambda }|x'-y|\leq \tau_{\Lambda }2|x''-y|\leq 2\tau_{\Lambda }^{2}|\Lambda x''-\Lambda y|\,,
\]
and thus the first of the second two inequalities of statement (i) holds true.
The second of the second two inequalities of statement (i) can be proved by interchanging the roles of $x'$ and $x''$. 

We now prove only the second inequalities in statements  (ii), (iv). Indeed the first inequalities follow by the second ones and by the equality $\tau_{\Lambda }=1$ when $\Lambda $ is the identity matrix.  The first of the second inequalities  in (ii) for $\gamma\geq 0$ follows by raising the inequality
$|\Lambda x'-\Lambda y|\leq (2\tau_{\Lambda }^{2}) |\Lambda x''-\Lambda y|$ of statement (i) to the power $\gamma$. Instead for $\gamma<0$ the same inequality follows by raising the inequality
$|\Lambda x''-\Lambda y|\leq (2\tau_{\Lambda }^{2})   |\Lambda x'-\Lambda y|$ of statement (i) to the power $\gamma$. The  second of the second inequalities of (ii) can be proved by interchanging the roles of $x'$ and $x''$. 

Statement (iii) follows by a direct application of (ii). To prove (iv) and (v), we follow
Cialdea~\cite[\S 8]{Ci00}. We first consider (iv), and we assume that
$|\Lambda x'-\Lambda y|\leq |\Lambda x''-\Lambda y|$.   By the Lagrange Theorem,  there exists $\zeta\in [|\Lambda x'-\Lambda y|,|\Lambda x''-\Lambda y|]$ such that
\[
||\Lambda x'-\Lambda y|^{\gamma}-|\Lambda x''-\Lambda y|^{\gamma}|\leq |\gamma|\zeta^{\gamma-1}\,|\,|\Lambda x'-\Lambda y|-|\Lambda x''-\Lambda y|\,|\,.
\]
If $\gamma\geq 1$, the inequality $\zeta\leq |\Lambda x''-\Lambda y|$ and (i) imply that
\[
\zeta^{\gamma-1}\leq |\Lambda x''-\Lambda y|^{\gamma-1}\leq (2\tau_{\Lambda }^{2})^{|\gamma-1|}|\Lambda x'-\Lambda y|^{\gamma-1}\,.
\]
If $\gamma<1$, then inequalities $\zeta\geq |\Lambda x'-\Lambda y|$ and $\tau_{\Lambda }\geq 1$   imply that
\[
\zeta^{\gamma-1}\leq |\Lambda x'-\Lambda y|^{\gamma-1}\leq (2\tau_{\Lambda }^{2})^{|\gamma-1|}
|\Lambda x'-\Lambda y|^{\gamma-1}\,.
\]
Then we have 
\begin{eqnarray}
\label{rec1}
\lefteqn{
||\Lambda x'-\Lambda y|^{\gamma}-|\Lambda x''-\Lambda y|^{\gamma}|
}
\\ \nonumber
&&\qquad\qquad
\leq
|\gamma|(2\tau_{\Lambda }^{2})^{|\gamma-1|}
|\,|\Lambda x'-\Lambda y| -|\Lambda x''-\Lambda y|\, |\,|\Lambda x'-\Lambda y|^{\gamma-1}\,,
\end{eqnarray}
which implies the validity of (iv). Similarly, in case $|\Lambda x'-\Lambda y|> |\Lambda x''-\Lambda y|$, we can prove that (\ref{rec1}) holds with $x'$ and $x''$ interchanged. Then (ii)  implies the validity of (iv). 

We now consider statement (v) and we assume that $|x'-y|\leq |x''-y|$. By the Lagrange Theorem,  there exists $\zeta\in [|x'-y|,|x''-y|]$ such that
\begin{equation}
\label{rec2}
|\ln  |x'-y| -\ln  |x''-y| |\leq \zeta^{-1}|\,|x'-y|-|x''-y|\,|\leq \zeta^{-1}|x'-x''| \,.
\end{equation}
By the above assumption, $\zeta^{-1}\leq|x'-y|^{-1}$, and thus statement (v) follows. Similarly, if $|x'-y|> |x''-y|$, we can prove that (\ref{rec2}) holds with $x'$ and $x''$ interchanged and  (i) implies that $\zeta^{-1}\leq|x''-y|^{-1}\leq 2|x'-y|^{-1}$, which implies the validity of (v). \hfill  $\Box$ 

\vspace{\baselineskip}

\begin{lem}
\label{fanes}
Let $G$ be a nonempty bounded subset of ${\mathbb{R}}^{n}$. Then the following statements hold.  
\begin{enumerate}
\item[(i)] Let $F\in {\mathrm{Lip}}(\partial{\mathbb{B}}_{n}\times [0,{\mathrm{diam}}\,(G)])$ with
\begin{eqnarray*}
\lefteqn{
{\mathrm{Lip}}(F)
\equiv\biggl\{\biggr.
\frac{|F(\theta',r')-F(\theta'',r'')|}{ |\theta'-\theta''|+|r'-r''| }:\,
}
\\ \nonumber
&&\qquad\quad
(\theta',r'),(\theta'',r'')\in \partial{\mathbb{B}}_{n}\times [0,{\mathrm{diam}}\,(G)],\ (\theta',r')\neq (\theta'',r'')
\biggl.\biggr\}\,.
\end{eqnarray*}
Then 
\begin{eqnarray}
\label{fanes1}
\lefteqn{
\left|
F\left(
\frac{x'-y}{|x'-y|},|x'-y|
\right)
-
F\left(
\frac{x''-y}{|x''-y|},|x''-y|
\right)
\right|
}
\\ \nonumber
&&\quad
\leq
{\mathrm{Lip}}(F) (2+   {\mathrm{diam}}\,(G))
\frac{|x'-x''|}{|x'-y|}\,
\quad\forall y\in G
\setminus {\mathbb{B}}_{n}(x',2|x'-x''|)\,,
\end{eqnarray}
for all $x',x''\in G$, $x'\neq x''$. In particular, if $f\in C^{1}(\partial{\mathbb{B}}_{n}\times{\mathbb{R}},{\mathbb{C}})$, then
\begin{eqnarray*}
\lefteqn{
M_{f,G}\equiv\sup\biggl\{\biggr.
\left|
f\left(
\frac{x'-y}{|x'-y|},|x'-y|
\right)
-
f\left(
\frac{x''-y}{|x''-y|},|x''-y|
\right)
\right|\frac{|x'-y|}{|x'-x''|}
}
\\ \nonumber
&& \qquad\qquad\quad \  
:\,x',x''\in G, x'\neq x'',  y\in G   
\setminus {\mathbb{B}}_{n}(x',2|x'-x''|)
\biggl.\biggr\}<\infty\,.
\end{eqnarray*}
\item[(ii)] Let $W$ be an open neighbourhood of ${\mathrm{cl}}(G-G)$. Let $f\in C^{1}(W,{\mathbb{C}})$. Then 
\begin{eqnarray*}
\lefteqn{
\tilde{M}_{f,G}\equiv
 \sup\biggl\{\biggr.
|
f(x'-y)-f(x''-y)|\,|x'-x''|^{-1}
}
\\ \nonumber
&&\qquad\qquad\qquad\qquad\qquad\qquad\quad
:\,x',x''\in G, x'\neq x'',  y\in G 
\biggl.\biggr\}<\infty\,.
\end{eqnarray*}
Here $G-G\equiv\{y_{1}-y_{2}:\ y_{1}, y_{2}\in G\}$.
\end{enumerate}
\end{lem}
{\bf Proof.} We first consider statement (i). The Lipschitz continuity of $F$ 
implies that the left hand side of (\ref{fanes1})  is less or equal to
\begin{eqnarray*}
\lefteqn{
{\mathrm{Lip}}(F) \biggl\{\biggr.
\biggl|\biggr.
\frac{x'-y}{|x'-y|} -\frac{x''-y}{|x''-y|}
\biggl.\biggr|
+
\biggl|\biggr.
|x'-y|-|x''-y|
\biggl.\biggr|
\biggl.\biggr\}
}
\\ \nonumber
&&\qquad
\leq 
{\mathrm{Lip}}(F) \biggl\{\biggr.
|x''-y|
\biggl|\biggr.
\frac{1}{|x''-y|}-\frac{1}{|x'-y|}
\biggl.\biggr|
\\ \nonumber
&&\qquad\qquad
+\frac{1}{|x'-y|} |\, |x''-y|- |x'-y|\, |
+|x'-x''|\biggl.\biggr\}
\\ \nonumber
&&\qquad
\leq 
{\mathrm{Lip}}(F) \biggl\{\biggr.
|x''-y|\frac{  |x'-x''|  }{|x''-y|\,|x'-y|}
+\frac{|x'-x''|}{|x'-y|}+|x'-x''|\biggl.\biggr\}
\\ \nonumber
&&\qquad
\leq 
{\mathrm{Lip}}(F) |x'-x''|
\biggl\{\biggr.
\frac{2+|x'-y| }{|x'-y|}
\biggl.\biggr\}\,,
\end{eqnarray*}
and thus inequality (\ref{fanes1}) holds true. 

Since $\partial {\mathbb{B}}_{n}\times{\mathbb{R}}$ is a manifold of class $C^{\infty}$ imbedded into ${\mathbb{R}}^{n+1}$, there exists $F\in C^{1}({\mathbb{R}}^{n+1}) $ which extends $f$. Since $\partial{\mathbb{B}}_{n}\times [0,{\mathrm{diam}}\,(G)]$ is a compact subset of ${\mathbb{R}}^{n+1}$, $F$ is Lipschitz continuous on $\partial{\mathbb{B}}_{n}\times [0,{\mathrm{diam}}\,(G)]$
and the second part of statement (i) follows by inequality (\ref{fanes1}). 

We now consider statement (ii). Since $f\in C^{1}(W,{\mathbb{C}})$, $f$ is Lipschitz continuous on the compact set ${\mathrm{cl}}(G-G)$, and statement (ii) follows.\hfill  $\Box$ 

\vspace{\baselineskip}

Then we have the following well known statement.
\begin{lem}
\label{com}
Let $\alpha\in ]0,1]$. Let $\Omega$ be a bounded open connected subset of ${\mathbb{R}}^{n}$ of class $C^{1,\alpha}$. Then there exists $c_{\Omega,\alpha}>0$ such that
\[
|\nu (y)\cdot (x-y)|\leq c_{\Omega,\alpha} |x-y|^{1+\alpha}
\qquad
\forall x,y\in \partial\Omega\,.
\]
\end{lem}
Next we introduce a list of classical inequalities which can be verified by 
exploiting the local parametrizations of $\partial\Omega$.
\begin{lem}
\label{comin}
Let $\Omega$ be a bounded open Lipschitz subset of ${\mathbb{R}}^{n}$. Then the following statements hold.
\begin{enumerate}
\item[(i)]  Let $\gamma\in ]-\infty, n-1[$. Then 
\[
c'_{\Omega,\gamma}\equiv \sup_{x\in\partial\Omega}
\int_{\partial\Omega}\frac{d\sigma_{y}}{|x-y|^{\gamma}}
<+\infty\,.
\]
\item[(ii)] Let $\gamma\in ]-\infty, n-1[$. Then 
\begin{eqnarray*}
\lefteqn{
c''_{\Omega,\gamma}
\equiv \sup_{x',x''\in\partial\Omega,\ x'\neq x''}
}
\\ \nonumber
&& \qquad\qquad
|x'-x''|^{-(n-1)+\gamma}\int_{{\mathbb{B}}_{n}(x',3|x'-x''|)\cap\partial\Omega}
\frac{d\sigma_{y}}{|x'-y|^{\gamma}}
<+\infty\,.
\end{eqnarray*}

\item[(iii)] Let $\gamma\in ]n-1,+\infty[$. Then 
\[
c'''_{\Omega,\gamma}
\equiv
 \sup_{x',x''\in\partial\Omega,\ x'\neq x''}
 |x'-x''|^{ -(n-1)+\gamma}
 \int_{\partial\Omega\setminus
 {\mathbb{B}}_{n}(x',2|x'-x''|)
 }
 \frac{d\sigma_{y}}{|x'-y|^{\gamma}}\,,
\]
is finite.
\item[(iv)] 
\begin{eqnarray*}
\lefteqn{
c^{iv}_{\Omega}
\equiv
 \sup_{x',x''\in\partial\Omega,\ 0<|x'- x''|<1/e}
 }
\\ \nonumber
&&\qquad\qquad 
 |\ln |x'-x''| |^{-1}
 \int_{\partial\Omega\setminus
 {\mathbb{B}}_{n}(x',2|x'-x''|)
 }
 \frac{d\sigma_{y}}{|x'-y|^{n-1}}<+\infty\,,
\end{eqnarray*}
\end{enumerate}
\end{lem}

\section{Preliminaries on the fundamental solution}
\label{prefundsol}
We first introduce a formula for the fundamental solution of $P[{\mathbf{a}},D]$. To do so, we follow a formulation of Dalla Riva \cite[Thm.~5.2, 5.3]{Da13} and Dalla Riva, Morais and Musolino \cite[Thm.~3.1, 3.2]{DaMoMu13} (see also John~\cite{Jo55}, and Miranda~\cite{Mi65} for homogeneous operators, and Mitrea and Mitrea~\cite[p.~203]{MitMit13}.)
\begin{thm}
\label{fundsol}
Let ${\mathbf{a}}$ be as in (\ref{introd0}), (\ref{ellip}). Let $S_{ {\mathbf{a}} }$ be a fundamental solution of $P[{\mathbf{a}},D]$. Then there exist an analytic function $A_{0}$ from $\partial{\mathbb{B}}_{n}$ to ${\mathbb{C}}$, and an analytic function $A_{1}$ from $\partial{\mathbb{B}}_{n}\times{\mathbb{R}}$ to ${\mathbb{C}}$, and $b_{0}\in {\mathbb{C}}$, and an analytic function $B_{1}$ from ${\mathbb{R}}^{n}$ to ${\mathbb{C}}$ such that $B_{1}(0)=0$, and an analytic function $C$ from ${\mathbb{R}}^{n}$ to ${\mathbb{C}}$ such that
\begin{equation}
\label{fundsol1}
S_{ {\mathbf{a}} }(x)
= |x|^{2-n}A_{0}(\frac{x}{|x|})+|x|^{3-n}A_{1}(\frac{x}{|x|},|x|)
+b_{0}\ln  |x| +B_{1}(x)\ln  |x|+C(x)
\,,
\end{equation}
for all $x\in {\mathbb{R}}^{n}\setminus\{0\}$,
and such that both $b_{0}$ and $B_{1}$   equal zero if $n$ is  odd. Moreover, 
 \[
 |x|^{2-n}A_{0}(\frac{x}{|x|})+\delta_{2,n}b_{0}\ln  |x|
 \]
is a fundamental solution for the principal part
$\sum_{l,j=1}^{n}\partial_{x_{l}}(a_{lj}\partial_{x_{j}} )
$ of $P[{\mathbf{a}},D]$. Here $\delta_{2,n}$ denotes the Kronecker symbol. Namely, 
\[
\delta_{2,n}=1\ {\mathrm{if}}\ n=2,\qquad
\delta_{2,n}=0\ {\mathrm{if}}\ n>2\,.
\]
\end{thm}
Then we have the following.
\begin{corol}
\label{normalizer}
Let ${\mathbf{a}}$ be as in (\ref{introd0}), (\ref{ellip}). Let $S_{ {\mathbf{a}} }$ be a fundamental solution of $P[{\mathbf{a}},D]$. Then the following statements hold.
\begin{enumerate}
\item[(i)] If $n\geq 3$, then there exists one and only one fundamental solution of the principal part
$\sum_{l,j=1}^{n}\partial_{x_{l}}(a_{lj}\partial_{x_{j}} )
$ of $P[{\mathbf{a}},D]$ which is positively homogeneous of degree $2-n$ in 
${\mathbb{R}}^{n}\setminus\{0\}$. 
\item[(ii)] If $n=2$, then there exists one and only one fundamental solution $S(x)$ of the principal part
$\sum_{l,j=1}^{n}\partial_{x_{l}}(a_{lj}\partial_{x_{j}} )
$ of $P[{\mathbf{a}},D]$ such that 
\[
\beta_{0}\equiv \lim_{x\to 0}\frac{S(x)}{\ln  |x|}\in {\mathbb{C}}\,,
\qquad
\int_{\partial{\mathbb{B}}_{n}}S\,d\sigma=0\,,
\]
and such that $S(x)-\beta_{0}\ln |x|$ is  positively homogeneous of degree $0$ in 
${\mathbb{R}}^{n}\setminus\{0\}$. 
\end{enumerate}
\end{corol}
{\bf Proof.} We retain the notation of Theorem \ref{fundsol}. We first consider statement (i). By Theorem \ref{fundsol}, the function 
$|x|^{2-n}A_{0}(\frac{x}{|x|})$ is a fundamental solution of the principal part of $P[{\mathbf{a}},D]$ and is clearly positively homogeneous  of degree $2-n$. Now assume that $u$ is a fundamental solution of the principal part of $P[{\mathbf{a}},D]$ and that $u$  is positively homogeneous  of degree $2-n$ in 
${\mathbb{R}}^{n}\setminus\{0\}$. Then the difference
\[
w(x)\equiv|x|^{2-n}A_{0}(\frac{x}{|x|})-u(x)
\]
defines an entire real analytic function in ${\mathbb{R}}^{n}$ and is positively homogeneous  of degree $2-n$ in 
${\mathbb{R}}^{n}\setminus\{0\}$. In particular,
\[
\lambda^{n-2}w(\lambda x)=w(x)
\qquad\forall (\lambda,x)\in ]0,+\infty[\times({\mathbb{R}}^{n}\setminus\{0\})\,,
\]
and accordingly
\[
\lambda^{(n-2)+|\beta|}D^{\beta}w(\lambda x)=D^{\beta}w(x)
\qquad\forall (\lambda,x)\in ]0,+\infty[\times({\mathbb{R}}^{n}\setminus\{0\})\,,
\]
for all $\beta\in {\mathbb{N}}^{n}$. Then by letting $\lambda$ tend to $0^{+}$, we obtain $D^{\beta}w(0)=0$ for all $\beta\in {\mathbb{N}}^{n}$. Since $w$ is real analytic, we deduce that $w$ is equal to $0$ in ${\mathbb{R}}^{n}$ and thus statement (i) holds. 

We now assume that $n=2$. By Theorem \ref{fundsol}, the function 
$S(x)\equiv A_{0}(\frac{x}{|x|})-\frac{1}{2\pi}\int_{\partial{\mathbb{B}}_{n}}A_{0}\,d\sigma+b_{0}\ln  |x|$ is a fundamental solution of the principal part of $P[{\mathbf{a}},D]$ and satisfies the conditions of statement (ii). We now
 assume that  $u$ is another fundamental solution of the 
principal part as in (ii). Then the difference
\[
w(x)\equiv A_{0}(\frac{x}{|x|})-\frac{1}{2\pi}\int_{\partial{\mathbb{B}}_{n}}A_{0}\,d\sigma+b_{0}\ln |x|-u(x)
\]
defines an entire real analytic function in ${\mathbb{R}}^{n}$ and we have
\[
0=\lim_{x\to 0}\frac{w(x)}{\ln  |x|}=\lim_{x\to 0}\frac{A_{0}(\frac{x}{|x|})-\frac{1}{2\pi}\int_{\partial{\mathbb{B}}_{n}}A_{0}\,d\sigma}{\ln  |x|}
+b_{0}-\lim_{x\to 0}\frac{u(x)}{\ln |x|}\,,
\]
and accordingly
\[
b_{0}=\lim_{x\to 0}\frac{u(x)}{\ln |x|}\equiv\beta_{0}\in {\mathbb{C}}\,.
\]
Then our assumption implies that the analytic function
\[
u(x)-\beta_{0}\ln |x|=u(x)-b_{0}\ln |x|\,,
\]
is positively homogeneous of degree $0$ in ${\mathbb{R}}^{n}\setminus\{0\}$. Hence, there exists a function $g_{0}$ from $\partial{\mathbb{B}}_{n}$ to ${\mathbb{C}}$ such that
\[
u(x)-b_{0}\ln |x|=g_{0}(\frac{x}{|x|})\qquad\forall x\in {\mathbb{R}}^{n}\setminus\{0\}\,.
\]
In particular, $g_{0}$ is real analytic and
\begin{eqnarray*}
\lefteqn{
w(x)=A_{0}(\frac{x}{|x|})-\frac{1}{2\pi}\int_{\partial{\mathbb{B}}_{n}}A_{0}\,d\sigma+b_{0}\ln |x|-(  g_{0}(\frac{x}{|x|})  +b_{0}\ln |x|)
}
\\
&&\qquad\qquad\qquad\qquad\qquad\qquad
=A_{0}(\frac{x}{|x|})-\frac{1}{2\pi}\int_{\partial{\mathbb{B}}_{n}}A_{0}\,d\sigma-g_{0}(\frac{x}{|x|})\,.
\end{eqnarray*}
Moreover, $w$ must be positively homogeneous of degree $0$ in ${\mathbb{R}}^{n}\setminus\{0\}$.   Since  $w$ is  continuous at $0$, $w$ must be  constant in the whole of ${\mathbb{R}}^{n}$. Since 
$\int_{\partial{\mathbb{B}}_{n}}w\,d\sigma=\int_{\partial{\mathbb{B}}_{n}}S\,d\sigma-\int_{\partial{\mathbb{B}}_{n}}u\,d\sigma=0$, such a constant must equal $0$ and thus $A_{0}(\frac{x}{|x|})  -\frac{1}{2\pi}\int_{\partial{\mathbb{B}}_{n}}A_{0}\,d\sigma 
=g_{0}(\frac{x}{|x|})$ for all $x\in 
{\mathbb{R}}^{n}\setminus\{0\}$. Hence, $u(x)=A_{0}(\frac{x}{|x|})-\frac{1}{2\pi}\int_{\partial{\mathbb{B}}_{n}}A_{0}\,d\sigma+b_{0}\ln |x|$ and statement (ii) follows. \hfill  $\Box$ 

\vspace{\baselineskip}

Then we can introduce the following.
\begin{defn}
\label{normalized}
 Let ${\mathbf{a}}$ be as in (\ref{introd0}), (\ref{ellip}). 
 We define as normalized fundamental solution  of the principal part of $P[{\mathbf{a}},D]$, to be the only  fundamental solution of Corollary \ref{normalizer}
\end{defn}

By Theorem \ref{fundsol} and by Corollary \ref{normalizer}, the normalized 
fundamental solution of the principal part of $P[{\mathbf{a}},D]$, equals
\[
|x|^{2-n}A_{0}(\frac{x}{|x|})
\]
if $n\geq 3$ and
\[
A_{0}(\frac{x}{|x|})-\frac{1}{2\pi}\int_{\partial{\mathbb{B}}_{n}}A_{0}\,d\sigma+b_{0}\ln |x|
\]
if $n=2$, where $A_{0}$ and $b_{0}$ are as in Theorem \ref{fundsol}. We now see that if the principal coefficients of $P[{\mathbf{a}},D]$ are real, then the normalized  fundamental solution of the principal part of $P[{\mathbf{a}},D]$ has a very specific form. To do so, we introduce the fundamental solution $S_{n}$ of the Laplace operator. Namely, we set
\[
S_{n}(x)\equiv
\left\{
\begin{array}{lll}
\frac{1}{s_{n}}\ln  |x| \qquad &   \forall x\in 
{\mathbb{R}}^{n}\setminus\{0\},\quad & {\mathrm{if}}\ n=2\,,
\\
\frac{1}{(2-n)s_{n}}|x|^{2-n}\qquad &   \forall x\in 
{\mathbb{R}}^{n}\setminus\{0\},\quad & {\mathrm{if}}\ n>2\,,
\end{array}
\right.
\]
where $s_{n}$ denotes the $(n-1)$ dimensional measure of 
$\partial{\mathbb{B}}_{n}$. Then we have the following elementary statement, which can be verified by the chain rule and by Corollary \ref{normalizer} (cf.~\textit{e.g.},  
Dalla Riva~\cite{Da14}.)
\begin{lem}
\label{exr}
Let ${\mathbf{a}}$ be as in (\ref{introd0}), (\ref{ellip}), (\ref{symr}). Then there exists an invertible matrix $T\in M_{n}({\mathbb{R}})$ such that
\begin{equation}
\label{exr1}
a^{(2)}=TT^{t}
\end{equation}
and the function
\[
S_{a^{(2)}}(x)\equiv \frac{1}{\sqrt{\det a^{(2)} }}S_{n}(T^{-1}x)\qquad\forall x\in {\mathbb{R}}^{n}\setminus\{0\}\,,
\]
coincides with the normalized fundamental solution  of the principal part of  $P[{\mathbf{a}},D]$ if $n\geq 3$, and coincides with the normalized fundamental solution  of the principal part of  $P[{\mathbf{a}},D]$ up to an additive constant if $n=2$. 
\end{lem}
Then Theorem \ref{fundsol}, and Corollary \ref{normalizer},  and Lemma \ref{exr} imply the validity of the following.
\begin{corol}
\label{ourfs} 
Let ${\mathbf{a}}$ be as in (\ref{introd0}), (\ref{ellip}), (\ref{symr}). Let $T\in M_{n}({\mathbb{R}})$  be as in (\ref{exr1}). Let $S_{ {\mathbf{a}} }$ be a fundamental solution of $P[{\mathbf{a}},D]$. 

Then there exist an analytic function $A_{1}$ from $\partial{\mathbb{B}}_{n}\times{\mathbb{R}}$ to ${\mathbb{C}}$, and an analytic function $B_{1}$ from ${\mathbb{R}}^{n}$ to ${\mathbb{C}}$ such that $B_{1}(0)=0$, and an analytic function $C$ from ${\mathbb{R}}^{n}$ to ${\mathbb{C}}$ such that
\begin{eqnarray}
\label{ourfs1}
\lefteqn{
S_{ {\mathbf{a}} }(x)
= 
\frac{1}{\sqrt{\det a^{(2)} }}S_{n}(T^{-1}x)
}
\\ \nonumber
&&\qquad
+|x|^{3-n}A_{1}(\frac{x}{|x|},|x|)
 +(B_{1}(x)+b_{0}(1-\delta_{2,n}))\ln  |x|+C(x)
\,,
\end{eqnarray}
for all $x\in {\mathbb{R}}^{n}\setminus\{0\}$,
 and such that both $b_{0}$ and $B_{1}$   equal zero
if $n$ is odd. Moreover, 
 \[
 \frac{1}{\sqrt{\det a^{(2)} }}S_{n}(T^{-1}x) 
 \]
is a fundamental solution for the principal part
  of $P[{\mathbf{a}},D]$.
\end{corol}
Next we prove the following technical statement. 
\begin{lem}
\label{fun} 
Let ${\mathbf{a}}$ be as in (\ref{introd0}), (\ref{ellip}).  Let $S_{ {\mathbf{a}} }$ be a fundamental solution of $P[{\mathbf{a}},D]$. Let $G$ be a nonempty bounded subset of ${\mathbb{R}}^{n}$. 
\begin{enumerate}
\item[(i)] Let $\gamma\in [0,1[$. Then
\begin{equation}
\label{fun1}
C_{0, S_{ {\mathbf{a}} },G,n-1-\gamma}\equiv
\sup_{ 0<|x|\leq{\mathrm{diam}}\,(G) }|x|^{n-1-\gamma}| S_{ {\mathbf{a}} } (x)|<+\infty\,.
\end{equation}
If $n>2$, then (\ref{fun1}) holds also for $\gamma=1$. 
\item[(ii)]
\begin{eqnarray*}
\lefteqn{
\tilde{C}_{0, S_{ {\mathbf{a}} },G}\equiv
\sup\biggl\{\biggr.
\frac{|x'-y|^{n-1}}{|x'-x''|}
|S_{ {\mathbf{a}} }(x'-y)-
S_{ {\mathbf{a}} }(x''-y)
|
}
\\ \nonumber
&&\qquad\qquad\quad
:\,x',x''\in G, x'\neq x'',  y\in G  
\setminus {\mathbb{B}}_{n}(x',2|x'-x''|)
\biggl.\biggr\}<\infty\,.
\end{eqnarray*}
\end{enumerate}
\end{lem}
{\bf Proof.} Statement (i) is an immediate consequence of formula (\ref{fundsol1}).
We now prove statement (ii). To do so, we resort to formula (\ref{fundsol1}) and we set
\begin{eqnarray*}
A(\theta,r)&\equiv& A_{0}(\theta)+rA_{1}(\theta,r)\qquad
\forall (\theta,r)\in \partial{\mathbb{B}}_{n}\times {\mathbb{R}},
\\
B(x)&\equiv& b_{0}+B_{1}(x)\qquad\forall x\in {\mathbb{R}}^{n}\,.
\end{eqnarray*}
Then Lemmas \ref{rec} and  \ref{fanes} imply that  
\begin{eqnarray*}
\lefteqn{
|S_{ {\mathbf{a}} }(x'-y)-
S_{ {\mathbf{a}} }(x''-y)
|
}
\\ \nonumber
&& 
\leq
|x'-y|^{2-n}\left|
A\left(
\frac{x'-y}{|x'-y|},|x'-y|
\right)
-
A\left(
\frac{x''-y}{|x''-y|},|x''-y|
\right)
\right|
\\ \nonumber
&& \quad
+\left|  
A\left(
\frac{x''-y}{|x''-y|},|x''-y|
\right)
\right|
|\,
|x'-y|^{2-n}-|x''-y|^{2-n}
\,|
\\ \nonumber
&& \quad
+
|\ln  |x'-y||\,|B(x'-y)-B(x''-y)|
\\ \nonumber
&& \quad
+|B(x''-y)|\,
|\ln  |x'-y|-\ln  |x''-y| |
+|C(x'-y)-C(x''-y)|
\\ \nonumber
&& 
\leq |x'-y|^{2-n}M_{A,G}\frac{|x'-x''|}{|x'-y|}
+\left(\sup_{\partial{\mathbb{B}}_{n}\times [0,{\mathrm{diam}}\,(G)]} |A|\right) m_{2-n}\frac{|x'-x''|}{|x'-y|^{n-1}}
\\ \nonumber
&& \quad
+|\ln  |x'-y|\,| \tilde{M}_{B,G}|x'-x''|
+\sup_{G-G}|B|2\frac{|x'-x''|}{|x'-y|}
+\tilde{M}_{C,G}|x'-x''|\,.
\end{eqnarray*}
Since $A$ is continuous on the compact set $\partial{\mathbb{B}}_{n}\times [0,{\mathrm{diam}}\,(G)]$, and $B$ and $C$ are continuous on the compact set ${\mathrm{cl}} (G-G)$, there exists $c>0$ such that
\begin{eqnarray*}
\lefteqn{
|S_{ {\mathbf{a}} }(x'-y)-
S_{ {\mathbf{a}} }(x''-y)
|
}
\\ \nonumber
&&\qquad
\leq c|x'-x''|
 \left\{
|x'-y|^{1-n}+\frac{1}{|x'-y|}+\ln |x'-y|+1
\right\}
\\ \nonumber
&&\qquad
\leq c|x'-x''|\,|x'-y|^{1-n}
\\ \nonumber
&&\qquad\quad
\times
 \{
1+ |x'-y|^{n-2} +|x'-y|^{n-1}\ln |x'-y|+|x'-y|^{n-1}
 \}\,,
\end{eqnarray*}
and thus statement (ii) holds. \hfill  $\Box$ 

\vspace{\baselineskip}

\begin{lem}
\label{grafun}
 Let ${\mathbf{a}}$ be as in (\ref{introd0}), (\ref{ellip}), (\ref{symr}). Let $T\in M_{n}({\mathbb{R}})$  be as in (\ref{exr1}). Let $S_{ {\mathbf{a}} }$ be a fundamental solution of $P[{\mathbf{a}},D]$. Let  $B_{1}$, $C$
 be as in Corollary \ref{ourfs}. 
 Let $G$ be a nonempty bounded subset of ${\mathbb{R}}^{n}$.  Then the following statements hold.
\begin{enumerate}
\item[(i)] There exists a real analytic function $A_{2}$ from $\partial{\mathbb{B}}_{n}\times{\mathbb{R}}$ to ${\mathbb{C}}^{n}$ such that
\begin{eqnarray}
\label{grafun1}
\lefteqn{
DS_{ {\mathbf{a}} }(x)=\frac{1}{ s_{n}\sqrt{\det a^{(2)} } }
|T^{-1}x|^{-n}x^{t}(a^{(2)})^{-1}
}
\\ \nonumber
&& 
+|x|^{2-n}A_{2}(\frac{x}{|x|},|x|)+DB_{1}(x)\ln |x|+DC(x)
\quad\forall x\in {\mathbb{R}}^{n}\setminus\{0\}\,.
\end{eqnarray}
\item[(ii)] 
\[
C_{1, S_{ {\mathbf{a}} },G}\equiv
\sup_{ 0<|x|\leq{\mathrm{diam}}\,(G)}|x|^{n-1}| DS_{ {\mathbf{a}} } (x)|<+\infty\,.
\]
\item[(iii)] 
\begin{eqnarray*}
\lefteqn{
\tilde{C}_{1, S_{ {\mathbf{a}} },G}\equiv
\sup\biggl\{\biggr.
\frac{|x'-y|^{n}}{|x'-x''|}
|DS_{ {\mathbf{a}} }(x'-y)-
DS_{ {\mathbf{a}} }(x''-y)
|
}
\\ \nonumber
&& \qquad\quad
:\,x',x''\in G, x'\neq x'',  y\in G  
\setminus {\mathbb{B}}_{n}(x',2|x'-x''|)
\biggl.\biggr\}<\infty\,.
\end{eqnarray*}
\end{enumerate}
\end{lem}
{\bf Proof.} By formula (\ref{ourfs1}) and by the chain rule, we have  
\begin{eqnarray}
\nonumber
\lefteqn{
DS_{ {\mathbf{a}} }(x)=\frac{1}{ s_{n} \sqrt{\det a^{(2)} } }
|T^{-1}x|^{-n}x^{t}(a^{(2)})^{-1}
+(3-n)|x|^{2-n}\frac{x^{t}}{|x|}A_{1}(\frac{x}{|x|},|x|)
}
\\ \nonumber
&& \quad
+|x|^{3-n}\biggl\{\biggr.
DA_{1}(\frac{x}{|x|},|x|) [|x|I-x\otimes x |x|^{-1}]|x|^{-2}
+\frac{\partial A_{1}}{\partial r}(\frac{x}{|x|},|x|)\frac{x^{t}}{|x|}\biggl.\biggr\}
\\ \label{grafun2}
&& \quad
+DB_{1}(x)\ln  |x|  
+(B_{1}(x)+b_{0}(1-\delta_{2,n}))\frac{x^{t}}{|x|^{2}}+DC(x) \,,
\end{eqnarray}
for all $x\in {\mathbb{R}}^{n}
\setminus\{0\}$, 
where we have still denoted by $A_{1}$ any analytic extension of the function $A_{1}$ of Corollary \ref{ourfs} to an open neighbourhood of $\partial{\mathbb{B}}_{n}\times{\mathbb{R}}$ in 
${\mathbb{R}}^{n+1}$ and where $x\otimes x$ denotes the matrix $(x_{l}x_{j})_{l,j=1,\dots,n}$. Next we consider the term $B_{1}(x)/|x|$. By the Fundamental Theorem of Calculus, we have
\begin{equation}
\label{grafun2a}
B_{1}(x)/|x|=\int_{0}^{1}DB_{1}(t\frac{x}{|x|}|x|)\frac{x}{|x|}\,dt\qquad\forall x\in {\mathbb{R}}^{n}\setminus\{0\}\,.
\end{equation}
Thus if we set
\[
\beta(\theta,r)=\int_{0}^{1}DB_{1}(t\theta r)\theta\,dt\qquad\forall (\theta,r)\in 
{\mathbb{R}}^{n}\times{\mathbb{R}}\,,
\]
the function $\beta$ is analytic and satisfies the equality
\begin{equation}
\label{grafun3}
B_{1}(x)/|x|=\beta(\frac{x}{|x|},|x|)\qquad\forall x\in {\mathbb{R}}^{n}
\setminus\{0\}\,.
\end{equation}
Then we can set
\begin{eqnarray*}
\lefteqn{
A_{2}(\theta,r)\equiv (3-n)\theta^{t} A_{1}(\theta,r)+
DA_{1}(\theta,r) [I-\theta\otimes\theta]
+\frac{\partial A_{1}}{\partial r}(\theta,r)\theta^{t} r
}
\\ \nonumber
&&\qquad\qquad 
+\beta(\theta,r)r^{n-2}\theta^{t}
+r^{n-3}\theta^{t}b_{0}(1-\delta_{2,n})
\qquad
\forall (\theta,r)\in \partial{\mathbb{B}}_{n}\times{\mathbb{R}}\,.
\end{eqnarray*}
 By the analyticity of $A_{1}$ and $\beta$, and by the equality $r^{n-3}\theta^{t}b_{0}(1-\delta_{2,n})=0$ if $n=2$, the function $A_{2}$ is analytic. Hence, equalities (\ref{grafun2}) and (\ref{grafun3}) imply the validity of statement (i). 

Next we turn to the proof of statement (ii). By Lemma \ref{t-1} (ii) and by the Schwartz inequality, we have 
\[
|T^{-1}x|^{-n}\, |x^{t}(a^{(2)})^{-1}|\leq |x|^{1-n}|T|^{n}|(a^{(2)})^{-1}|\,.
\]
Hence, formula (\ref{grafun1}) implies that
\begin{eqnarray*}
\lefteqn{
|x|^{n-1} |DS_{ {\mathbf{a}} }(x)|\leq \frac{1}{ s_{n} \sqrt{\det a^{(2)} } }|T|^{n}|(a^{(2)})^{-1}|  
}
\\ \nonumber
&&\quad
+\biggl\{\biggr.
|x| A_{2}(\frac{x}{|x|},|x|)+(|x|^{n-1}\ln |x|)DB_{1}(x)+|x|^{n-1}DC(x)
\biggl.\biggr\} \,,
\end{eqnarray*}
for all $x\in {\mathbb{R}}^{n}\setminus\{0\}$.  Then the continuity of  $A_{2}$ on the compact set $\partial{\mathbb{B}}_{n}\times [0,{\mathrm{diam}} (G)]$ and the continuity of $DB_{1}$ and $DC$ on the compact set ${\mathrm{cl}}{\mathbb{B}}_{n}(0,{\mathrm{diam}}\,(G))$  imply the validity of statement (ii). 

We now turn to statement (iii). Let  $x',x''\in G$, $x'\neq x''$,  $y\in G   
\setminus {\mathbb{B}}_{n}(x',2|x'-x''|)$. By statement (i), we have
\begin{eqnarray}
\label{grafun4}
\lefteqn{|DS_{ {\mathbf{a}} }(x'-y)-
DS_{ {\mathbf{a}} }(x''-y)
|
}
\\ \nonumber
&&
\leq \frac{1}{ s_{n}\sqrt{\det a^{(2)} } }
\biggl|\biggr.
|T^{-1}(x'-y)|^{-n}(x'-y)^{t}(a^{(2)})^{-1}
\\ \nonumber
&&\qquad\qquad\qquad\qquad
-
|T^{-1}(x''-y)|^{-n}(x''-y)^{t}(a^{(2)})^{-1}
\biggl.\biggr|
\\ \nonumber
&& 
+\biggl|\biggr.
|x'-y|^{2-n}A_{2}(\frac{x'-y}{|x'-y|},|x'-y|)
-
|x''-y|^{2-n}A_{2}(\frac{x''-y}{|x''-y|},|x''-y|)
\biggl.\biggr|
\\ \nonumber
&& 
+\biggl|\biggr.
\ln  |x'-y| DB_{1}(x'-y)-\ln  |x''-y| DB_{1}(x''-y)
\biggl.\biggr|
\\ \nonumber
&& 
+\biggl|\biggr.
DC(x'-y)-DC(x''-y)
\biggl.\biggr|\,.
\end{eqnarray}
We first estimate the first summand in the right hand side of inequality
(\ref{grafun4}).  By the triangular inequality, we have
\begin{eqnarray}
\nonumber
\lefteqn{
\biggl|\biggr.
|T^{-1}(x'-y)|^{-n}(x'-y)^{t}(a^{(2)})^{-1}
-
|T^{-1}(x''-y)|^{-n}(x''-y)^{t}(a^{(2)})^{-1}
\biggl.\biggr|
}
\\ \nonumber
&&\qquad\qquad\qquad
\leq
|x'-y|\,|(a^{(2)})^{-1}|\,
\biggl|\biggr.
|T^{-1}(x'-y)|^{-n}-|T^{-1}(x''-y)|^{-n}
\biggl.\biggr|
\\ \label{grafun5}
&&\qquad\qquad\qquad\quad
+|x'-x''|\,|(a^{(2)})^{-1}|\, |T^{-1}(x''-y)|^{-n}\,.
\end{eqnarray}
Then Lemmas \ref{t-1} (ii),  \ref{rec} (ii), (iv) with $\gamma=-n$, $\Lambda =T^{-1}$ imply that
\begin{eqnarray}
\label{grafun6}
\lefteqn{
\biggl|\biggr.
|T^{-1}(x'-y)|^{-n}-|T^{-1}(x''-y)|^{-n}
\biggl.\biggr|
}
\\ \nonumber
&&\qquad\qquad\qquad 
\leq m_{-n}(T^{-1})|T^{-1}x'-T^{-1}x''|
\,|T^{-1}x'-T^{-1}y|^{-n-1}
\\ \nonumber
&&\qquad\qquad\qquad 
\leq
m_{-n}(T^{-1})|T^{-1}|\,|T|^{n+1} |x'-x''|\,|x'-y|^{-n-1}\,,
\\ \nonumber
\lefteqn{
|T^{-1}(x''-y)|^{-n}
\leq |T|^{n}|x''-y|^{-n}\,,
}
\\ \nonumber
\lefteqn{
|x''-y|^{-n}\leq 2^{n}|x'-y|^{-n}\,.
}
\end{eqnarray}
Next we  estimate the second summand in the right hand side of inequality
(\ref{grafun4}). By Lemmas \ref{rec} (iv) and \ref{fanes} (i), the second summand
 is less or equal to
\begin{eqnarray}
\label{grafun7}
\lefteqn{
\biggl|\biggr.
|x'-y|^{2-n}-|x''-y|^{2-n}
\biggl.\biggr|
\,
\biggl|\biggr.
A_{2}(\frac{x''-y}{|x''-y|},|x''-y|)
\biggl.\biggr|
}
\\ \nonumber
&&\qquad\quad
+|x'-y|^{2-n}
\biggl|\biggr.
A_{2}(\frac{x'-y}{|x'-y|},|x'-y|)
-
A_{2}(\frac{x''-y}{|x''-y|},|x''-y|)
\biggl.\biggr|
\\ \nonumber
&&\qquad
\leq
m_{2-n}|x'-x''||x'-y|^{2-n-1}\sup_{
\partial{\mathbb{B}}_{n}\times [0,{\mathrm{diam}}(G)]
}|A_{2}|
\\ \nonumber
&&\qquad\quad
+|x'-y|^{2-n}
(\sum_{j=1}^{n}
M_{   A_{2,j}   ,G}  ) |x'-x''||x'-y|^{-1}\,.
\end{eqnarray}
 Next we  estimate the third summand in the right hand side of inequality
(\ref{grafun4}). By Lemmas \ref{rec} (v) and \ref{fanes} (ii), the third summand
 is less or equal to
\begin{eqnarray}
\label{grafun8}
\lefteqn{
\biggl|\biggr.
\ln  |x'-y|-\ln  |x''-y| 
\biggl.\biggr|\,
|DB_{1}(x''-y)|
}
\\ \nonumber
&&\quad
+|  \ln  |x'-y|  |\, |DB_{1}(x'-y)-DB_{1}(x''-y)|
\\ \nonumber
&&
\leq 2 |x'-x''|\,|x'-y|^{-1}\sup_{G-G}|DB_{1}|
+
(\sum_{j=1}^{n}\tilde{M}_{\frac{\partial B_{1}}{\partial x_{j}},G})
|x'-x''|\,|\ln  |x'-y| |
\\ \nonumber
&& 
\leq
|x'-x''|\,|x'-y|^{-n}
\biggl\{\biggr.
2 |x'-y|^{n-1}\sup_{G-G}|DB_{1}|
\\ \nonumber
&&\qquad\qquad\qquad\qquad\qquad
+(\sum_{j=1}^{n}\tilde{M}_{\frac{\partial B_{1}}{\partial x_{j}},G})
|x'-y|^{n}|\ln  |x'-y||
\biggl.\biggr\}\,.
\end{eqnarray}
Finally, Lemma \ref{fanes} (ii) implies that
\begin{eqnarray}
\label{grafun9}
\lefteqn{
|DC(x'-y)-DC(x''-y)|\leq 
(\sum_{j=1}^{n}\tilde{M}_{\frac{\partial C}{\partial x_{j}},G})
|x'-x''|
}
\\ \nonumber
&&\qquad
\leq  |x'-x''| |x'-y|^{-n}
(\sum_{j=1}^{n}\tilde{M}_{  \frac{\partial C}{\partial x_{j}},G })
\sup_{(x',y)\in G\times G}|x'-y|^{n}\,.
\end{eqnarray}
Then inequalities (\ref{grafun4})--(\ref{grafun9}) imply the validity of statement (iii).  \hfill  $\Box$ 

\vspace{\baselineskip}

 \section{Preliminary inequalities on the boundary operator}
We now turn to estimate the kernel $\overline{B^{*}_{\Omega,y}}\left(S_{{\mathbf{a}}}(x-y)\right)$ of the double layer potential of (\ref{introd3}). We do so under assumption  (\ref{symr}).
To do so, we introduce some basic inequalities for $\overline{B^{*}_{\Omega,y}}\left(S_{{\mathbf{a}}}(x-y)\right)$ by means of the following. 
\begin{lem}
\label{boest}
 Let ${\mathbf{a}}$ be as in (\ref{introd0}), (\ref{ellip}), (\ref{symr}). Let $T\in M_{n}({\mathbb{R}})$  be as in (\ref{exr1}). Let $S_{ {\mathbf{a}} }$ be a fundamental solution of $P[{\mathbf{a}},D]$. 
 
 Let $\alpha\in]0,1]$. Let $\Omega$ be a bounded open subset of ${\mathbb{R}}^{n}$ of class $C^{1,\alpha}$. Then the following statements hold.
\begin{enumerate}
\item[(i)] If $\alpha\in]0,1[$, then
\begin{equation}
\label{boest1}
b_{\Omega,\alpha}\equiv\sup
\biggl\{\biggr.
|x-y|^{n-1-\alpha}|\overline{B^{*}_{\Omega,y}}\left(S_{{\mathbf{a}}}(x-y)\right)|
:\, x,y\in \partial\Omega, x\neq y
\biggl.\biggr\}<+\infty\,.
\end{equation}
If $n>2$, then (\ref{boest1}) holds also for $\alpha=1$. 
\item[(ii)]
\begin{eqnarray*}
\lefteqn{
\tilde{b}_{\Omega,\alpha}\equiv\sup
\biggl\{\biggr.\frac{|x'-y|^{n-\alpha}}{|x'-x''|}
|\overline{B^{*}_{\Omega,y}}\left(S_{{\mathbf{a}}}(x'-y)\right)
-
\overline{B^{*}_{\Omega,y}}\left(S_{{\mathbf{a}}}(x''-y)\right)|:\,
}
\\ \nonumber
&& 
x',x''\in\partial\Omega, x'\neq x'', y\in\partial\Omega\setminus{\mathbb{B}}_{n}(x',2|x'-x''|)
\biggl.\biggr\}<+\infty\,.
\end{eqnarray*}
\end{enumerate}
\end{lem}
{\bf Proof.} By Lemma \ref{grafun} (i), we have
\begin{eqnarray}
\label{boest3}
\lefteqn{
\overline{B^{*}_{\Omega,y}}\left(S_{{\mathbf{a}}}(x-y)\right)
= -DS_{{\mathbf{a}}}(x-y)    a^{(2)}  \nu  (y)
- 
\nu^{t} (y)
a^{(1)}S_{{\mathbf{a}}}(x-y)
}
\\ \nonumber
&&\qquad\qquad
= 
-\frac{1}{s_{n}\sqrt{\det a^{(2)} }}|T^{-1}(x-y)|^{-n} (x-y)^{t}\nu (y)
\\ \nonumber
&&\qquad\qquad\quad
-
|x-y|^{2-n}
A_{2}(\frac{x-y}{|x-y|},|x-y|)a^{(2)}  \nu  (y)
\\ \nonumber
&&\qquad\qquad\quad
-DB_{1}(x-y)a^{(2)}  \nu  (y)\ln |x-y|
-   DC(x-y)a^{(2)}  \nu  (y)
\\ \nonumber
&&\qquad\qquad\quad
-\nu^{t}(y) a^{(1)}S_{{\mathbf{a}}}(x-y)
\qquad\forall x,y\in\partial\Omega, x\neq y\,.
\end{eqnarray}
Then by Lemmas  \ref{t-1} (ii), \ref{com}, \ref{fun} (i),   and by the   equality in (\ref{boest3}), we have 
\begin{eqnarray*}
\lefteqn{
|x-y|^{n-1-\alpha} |\overline{B^{*}_{\Omega,y}}\left(S_{{\mathbf{a}}}(x-y)\right)|
}
\\ \nonumber
&&\qquad\qquad
\leq \frac{1}{s_{n}\sqrt{\det a^{(2)} }}c_{\Omega,\alpha}|T|^{n}
|x-y|^{-n+1+\alpha+n-1-\alpha}
\\ \nonumber
&&\qquad\qquad\quad
+|x-y|^{2-1-\alpha}|a^{(2)}|\,|A_{2}(\frac{x-y}{|x-y|},|x-y|)|
\\ \nonumber
&&\qquad\qquad\quad
+|x-y|^{n-1-\alpha}|\ln |x-y|\,|\, |a^{(2)}|\,|DB_{1}(x-y)|
\\ \nonumber
&&\qquad\qquad\quad
+|x-y|^{n-1-\alpha}|a^{(2)}|\,|DC(x-y)|
+|  a^{(1)}  | C_{0,S_{ {\mathbf{a}} },
\partial\Omega,n-1-\alpha}\,,
\end{eqnarray*}
for all $x,y\in\partial\Omega$, $x\neq y$. If  either $\alpha\in]0,1[$  
or if $\alpha\in]0,1]$ and $n>2$, then  the right hand side is bounded for $x,y\in\partial\Omega$, $x\neq y$. Hence,  we conclude that statement (i) holds true.

Next we consider statement (ii). 
\begin{eqnarray}
\label{boest5}
\lefteqn{
|\overline{B^{*}_{\Omega,y}}\left(S_{{\mathbf{a}}}(x'-y)\right)
-
\overline{B^{*}_{\Omega,y}}\left(S_{{\mathbf{a}}}(x''-y)\right)|
}
\\ \nonumber
&& 
\leq
\frac{
\biggl|\biggr.
|T^{-1}(x'-y)|^{-n} (x'-y)^{t}\nu(y)
-
|T^{-1}(x''-y)|^{-n} (x''-y)^{t}\nu(y)
\biggl.\biggr|
}
{
s_{n}\sqrt{\det a^{(2)} }
}
\\ \nonumber
&& \quad
+|a^{(2)}|\,
\left|A_{2}(\frac{x'-y}{|x'-y|},|x'-y|)-A_{2}(\frac{x''-y}{|x''-y|},|x''-y|)\right|
 |x'-y|^{2-n} 
\\ \nonumber
&& \quad
+|a^{(2)}|\left|A_{2}(\frac{x''-y}{|x''-y|},|x''-y|)
\right|   |
|x'-y|^{2-n}-|x''-y|^{2-n}
|
\\ \nonumber
&& \quad
+|a^{(2)}|\,|
DB_{1}(x'-y)-DB_{1}(x''-y)|
\,|\ln |x'-y| |  
\\ \nonumber
&& \quad
+|a^{(2)}|\,|DB_{1}(x''-y)|\,
|\ln |x'-y| - \ln |x''-y|  | 
\\ \nonumber
&& \quad
+|a^{(2)}|\,|DC(x'-y)-DC(x''-y)| 
+|a^{(1)}|\,|S_{{\mathbf{a}}}(x'-y)-S_{{\mathbf{a}}}(x''-y)|
\end{eqnarray}
for all $x',x''\in\partial\Omega$, $x'\neq x''$,  $y\in\partial\Omega\setminus{\mathbb{B}}_{n}(x',2|x'-x''|)
$.
Next we denote by $J_{1}$ the first term in the right hand side of 
(\ref{boest5}). By Lemmas  \ref{t-1} (ii),  \ref{rec} (ii) and (iv) with $\gamma=-n$, $\Lambda =T^{-1}$, and by Lemma \ref{com}, we have
\begin{eqnarray}
\label{boest6}
\lefteqn{
J_{1}\leq \frac{ 1  }{  s_{n}\sqrt{\det a^{(2)} } }
}
\\ \nonumber
&&\quad
\times
\biggl\{\biggr.
\biggl|\biggr.
|T^{-1}(x'-y)|^{-n}-|T^{-1}(x''-y)|^{-n} 
\biggl.\biggr| |(x'-y)^{t}\nu(y)|
\\ \nonumber
&&\quad
+|T^{-1}(x''-y)|^{-n} \, |(x'-x'')^{t}\nu(y)|
\biggl.\biggr\}
\\ \nonumber
&&
\leq  \frac{ 1 }{  s_{n}\sqrt{\det a^{(2)} } }
\\ \nonumber
&&\quad
\times
\biggl\{\biggr.
m_{-n}(T^{-1})||T^{-1}x'-T^{-1}x''|\, |
T^{-1}x'-T^{-1}y
|^{-n-1}|x'-y|^{1+\alpha}c_{\Omega,\alpha}
\\ \nonumber
&&\quad
+2^{n}|T|^{n}|x'-y|^{-n}|(x'-x'')^{t}\nu(y)|\biggl.\biggr\}\,,
\end{eqnarray}
for all $x',x''\in\partial\Omega$, $x'\neq x''$,  $y\in\partial\Omega\setminus{\mathbb{B}}_{n}(x',2|x'-x''|)
$. Next we note that
\begin{eqnarray*}
\lefteqn{
|(x'-x'')^{t}\nu(y)|
\leq
|(x'-x'')^{t}(\nu(y)-\nu(x'))|
+|(x'-x'')^{t}\nu(x')|
}
\\ \nonumber
&&\qquad\qquad\qquad\qquad
\leq |x'-x''|\,|\nu|_{\alpha} |x'-y|^{\alpha}
+c_{\Omega,\alpha}|x'-x''|^{1+\alpha}
\\ \nonumber
&&\qquad\qquad\qquad\qquad
\leq |x'-x''|\, |x'-y|^{\alpha} (|\nu|_{\alpha}  +  c_{\Omega,\alpha})\,,
\end{eqnarray*}
and that accordingly
\begin{eqnarray}
\label{boest7}
\lefteqn{
J_{1}\leq 
\frac{   |x'-x''|  }{  s_{n}\sqrt{\det a^{(2)} }  }
\biggl\{\biggr.
m_{-n}(T^{-1})|T^{-1}|\,|T|^{n+1}|x'-y|^{-n-1}|x'-y|^{1+\alpha}c_{\Omega,\alpha}
}
\\ \nonumber
&&\qquad\qquad\qquad\qquad\qquad\qquad
+2^{n}|T|^{n}|x'-y|^{-n}|x'-y|^{\alpha} (|\nu|_{\alpha}+c_{\Omega,\alpha})
\biggl.\biggr\}\,,
\end{eqnarray}
for all $x',x''\in\partial\Omega$, $x'\neq x''$,  $y\in\partial\Omega\setminus{\mathbb{B}}_{n}(x',2|x'-x''|)
$. Next we denote by $J_{2}$ the sum of the  terms different from $J_{1}$ in the right hand side of 
(\ref{boest5}). Then Lemma \ref{rec} (iv), (v) and Lemmas \ref{fanes}, \ref{fun} (ii) imply that 
\begin{eqnarray}
\label{boest8}
\lefteqn{
J_{2}\leq |a^{(2)}|(\sum_{j=1}^{n}M_{A_{2,j},\partial\Omega})\frac{|x'-x''|}{|x'-y|}|x'-y|^{2-n}
}
\\ \nonumber
&&\qquad
+ |a^{(2)}|\sup_{\partial{\mathbb{B}}_{n}\times[0,{\mathrm{diam}}(\partial\Omega)]}|A_{2}|
m_{2-n}|x'-x''|\,|x'-y|^{1-n} 
\\ \nonumber
&&\qquad
+ |a^{(2)}|(\sum_{j=1}^{n}\tilde{M}_{\frac{\partial B_{1}}{\partial x_{j}},\partial\Omega})
|x'-x''|\,|\ln  |x'-y|| 
\\ \nonumber
&&\qquad
+ |a^{(2)}|\sup_{\partial\Omega-\partial\Omega}|DB_{1}|2\frac{|x'-x''| }{|x'-y| } \\ \nonumber
&&\qquad
+ \tilde{M}_{C}|x'-x''| 
+\tilde{C}_{0,S_{ {\mathbf{a}} },
\partial\Omega }
|a^{(1)}|\frac{|x'-x''| }{|x'-y|^{n-1}}
\,,
\end{eqnarray}
for all $x',x''\in\partial\Omega$, $x'\neq x''$,  $y\in\partial\Omega\setminus{\mathbb{B}}_{n}(x',2|x'-x''|)
$. By inequalities (\ref{boest5}), (\ref{boest7}), (\ref{boest8}), we conclude that statement (ii) holds.\hfill  $\Box$ 

\vspace{\baselineskip}

\section{Boundary norms for kernels}
\label{bonoke}

For each subset $A$ of ${\mathbb{R}}^{n}$, we find convenient to set
\[
\Delta_{A}\equiv \{ (x,y)\in A\times A:\, x=y\}\,.
\]
We now introduce a class of functions on $(\partial\Omega)^{2}\setminus\Delta_{\partial\Omega}$ which may carry a singularity as the variable tends to a point of the diagonal, just as in the case of the kernels of integral operators corresponding to layer potentials defined on the boundary of  an open subset $\Omega$ of ${\mathbb{R}}^{n}$. 
\begin{defn}
\label{k0}
Let $G$ be a nonempty bounded subset of ${\mathbb{R}}^{n}$.  Let $\gamma_{1}$, $\gamma_{2}$, $\gamma_{3}\in{\mathbb{R}}$.
We denote by ${\mathcal{K}}_{ \gamma_{1},\gamma_{2},\gamma_{3}   }(G)$
 the set of continuous functions $K$ from $(G\times G)\setminus\Delta_{G}$ to ${\mathbb{C}}$ such that
 \begin{eqnarray*}
\lefteqn{
\|K\|_{  {\mathcal{K}}_{ \gamma_{1},\gamma_{2},\gamma_{3}   }(G)  }
\equiv
\sup\biggl\{\biggr.
|x-y|^{ \gamma_{1} }|K(x,y)|:\,x,y\in G, x\neq y
\biggl.\biggr\}
}
\\ \nonumber
&&\qquad\qquad\qquad
+\sup\biggl\{\biggr.
\frac{|x'-y|^{\gamma_{2}}}{|x'-x''|^{\gamma_{3}}}
|  K(x',y)- K(x'',y)  |:\,
\\ \nonumber
&&\qquad\qquad\qquad 
x',x''\in G, x'\neq x'', y\in G\setminus{\mathbb{B}}_{n}(x',2|x'-x''|)
\biggl.\biggr\}<+\infty\,.
\end{eqnarray*}
\end{defn}
One can easily verify that $({\mathcal{K}}_{ \gamma_{1},\gamma_{2},\gamma_{3}   }(G),\|\cdot\|_{  {\mathcal{K}}_{ \gamma_{1},\gamma_{2},\gamma_{3}   }(G)  })$ is a Banach space.
\begin{rem}
 \label{exk}
 Let ${\mathbf{a}}$ be as in (\ref{introd0}), (\ref{ellip}). Let $S_{ {\mathbf{a}} }$ be a fundamental solution of $P[{\mathbf{a}},D]$. 
\begin{enumerate}
\item[(i)] Let $G$ be a nonempty bounded subset of ${\mathbb{R}}^{n}$. Then Lemma  \ref{fun} implies that $S_{  {\mathbf{a}} }(x-y)\in {\mathcal{K}}_{ n-1-\gamma,n-1,1  }(G)$ for all $\gamma\in[0,1[$ and that the same membership holds also for $\gamma=1$ if $n>2$. 
If we further assume that ${\mathbf{a}}$ satisfies $(\ref{symr})$, then Lemma \ref{grafun} implies that $\frac{\partial}{\partial x_{j}}S_{  {\mathbf{a}} }(x-y)\in
{\mathcal{K}}_{ n-1,n,1  }(G)$ for all $j\in \{1,\dots,n\}$.
\item[(ii)] Let ${\mathbf{a}}$ satisfy $(\ref{symr})$. Let $\alpha\in]0,1[$. Let $\Omega$   be a bounded open  subset of ${\mathbb{R}}^{n}$ of class $C^{1,\alpha}$. Then Lemma \ref{boest} implies that $\overline{B^{*}_{\Omega,y}}\left(S_{{\mathbf{a}}}(x-y)\right)\in {\mathcal{K}}_{ n-1-\alpha,n-\alpha,1  }(\partial\Omega)$.
\end{enumerate}
\end{rem}
For each $\theta\in]0,1]$, we define the function $\omega_{\theta}(\cdot)$ from $]0,+\infty[$ to itself by setting
\[
\omega_{\theta}(r)\equiv
\left\{
\begin{array}{ll}
r^{\theta}|\ln r | &r\in]0,r_{\theta}]\,,
\\
r_{\theta}^{\theta}|\ln r_{\theta} | & r\in ]r_{\theta},+\infty[\,,
\end{array}
\right.
\]
where
\[
r_{\theta} \equiv  e^{-1/\theta} \qquad \forall  \theta\in ]0,1]\,.
\]
Obviously, $\omega_{\theta}(\cdot) $ satisfies (\ref{om}),  (\ref{om1}), and (\ref{om3}) with $\alpha=\theta$. We also note that if ${\mathbb{D}}$ is a   subset of ${\mathbb{R}}^{n}$, then the following continuous imbedding holds
\[
C^{0,\omega_{\theta}(\cdot)}_{b}({\mathbb{D}})
\subseteq
C^{0,\theta'}_{b}({\mathbb{D}})
\]
for all $\theta'\in]0,\theta[$. We now consider the properties of an integral  operator with a kernel in the class
$ {\mathcal{K}}_{ \gamma_{1},\gamma_{2},\gamma_{3}   }(\partial\Omega) $. 
 \begin{prop}
\label{k0a}
Let $\Omega$ be a bounded open Lipschitz subset of ${\mathbb{R}}^{n}$. Let $\gamma_{1}\in ]-\infty,n-1[$,  $\gamma_{2}$, $\gamma_{3}\in{\mathbb{R}}$.
Then the following statements hold.
\begin{enumerate}
\item[(i)]  If $(K,\mu)\in
{\mathcal{K}}_{ \gamma_{1},\gamma_{2},\gamma_{3}   }(\partial\Omega)\times L^{\infty}(\partial\Omega)$, then the function $K(x,\cdot)\mu(\cdot)$ is integrable in $\partial\Omega$  for all $x\in\partial\Omega$, and the function $u[\partial\Omega,K,\mu]$ from $\partial\Omega$ to ${\mathbb{C}}$ defined by 
\begin{equation}
\label{k0a0}
u[\partial\Omega,K,\mu](x)\equiv\int_{\partial\Omega}K(x,y)\mu(y)\,d\sigma_{y} \qquad\forall x\in\partial\Omega\,,
\end{equation}
is continuous. Moreover, the bilinear map from ${\mathcal{K}}_{ \gamma_{1},\gamma_{2},\gamma_{3}   }(\partial\Omega)\times L^{\infty}(\partial\Omega)$ to $C^{0}(\partial\Omega)$, which takes $(K,\mu)$ to $u[\partial\Omega,K,\mu]$  is continuous. 
\item[(ii)] If $\gamma_{1}\in [n-2,n-1[$, $\gamma_{2}\in ]n-1,+\infty[$, 
$\gamma_{3}\in]0,1]$, $(n-1)-\gamma_{2}+\gamma_{3}\in]0,1]$, then the bilinear map from ${\mathcal{K}}_{ \gamma_{1},\gamma_{2},\gamma_{3}   }(\partial\Omega)\times L^{\infty}(\partial\Omega)$ to the space 
$C^{0,\min\{  (n-1)-\gamma_{1}, (n-1)-\gamma_{2}+\gamma_{3} \}}(\partial\Omega)$, which takes $(K,\mu)$ to $u[\partial\Omega,K,\mu]$  is continuous. 
\item[(iii)] If $\gamma_{1}\in [n-2,n-1[$, $\gamma_{2}=n-1$, $\gamma_{3}\in]0,1]$, then the bilinear map  from ${\mathcal{K}}_{ \gamma_{1},\gamma_{2},\gamma_{3}   }(\partial\Omega)\times L^{\infty}(\partial\Omega)$ to the space 
$C^{0,\max\{r^{ (n-1)-\gamma_{1}}, \omega_{\gamma_{3}}(r)\}}(\partial\Omega)$, which takes $(K,\mu)$ to $u[\partial\Omega,K,\mu]$  is continuous. 
\end{enumerate}
\end{prop}
{\bf Proof.} By definition of norm in ${\mathcal{K}}_{ \gamma_{1},\gamma_{2},\gamma_{3}   }(\partial\Omega) $, we have
\[
|K(x,y)\mu(y)|\leq   \| K \|_{  {\mathcal{K}}_{ \gamma_{1},\gamma_{2},\gamma_{3}   }(\partial\Omega)  }\|\mu\|_{  L^{\infty}(\partial\Omega)  }
\frac{1}{ |x-y|^{\gamma_{1}}}
\ \ \forall (x,y)\in (\partial\Omega)^{2}\setminus D_{\partial\Omega}\,.
\]
Then the function $K(x,\cdot)\mu(\cdot)$ is integrable in $\partial\Omega$  for all $x\in\partial\Omega$, and the Vitali Convergence Theorem implies that $u[\partial\Omega,K,\mu]$ is continuous on $\partial\Omega$ (cf.~\textit{e.g.}, Folland~\cite[(2.33) p.~60, p.~180]{Fo84}.) By Lemma \ref{comin} (i), we also have 
\begin{equation}
\label{k0a1}
\left|
\int_{\partial\Omega} K(x,y)\mu(y)\,d\sigma_{y}
\right|
\leq  \| K \|_{  {\mathcal{K}}_{ \gamma_{1},\gamma_{2},\gamma_{3}   }(\partial\Omega)  }\|\mu\|_{  L^{\infty}(\partial\Omega)  }c'_{\Omega,\gamma_{1}}\,.
\end{equation}
Hence, statement (i) follows. Next we turn to estimate the H\"{o}lder coefficient  of $u[\partial\Omega,K,\mu]$ under the assumptions of statements (ii) and (iii). Let $x',x''\in \partial\Omega$, $x'\neq x''$.   By Remark \ref{om4}, there is  no loss of generality in assuming that $0<|x'-x''|\leq r_{\gamma_{3} }$. Then the inclusion
${\mathbb{B}}_{n}(x',2|x'-x''|)\subseteq {\mathbb{B}}_{n}(x'',3|x'-x''|)$ and the triangular inequality imply  that
\begin{eqnarray}
\label{k0a2}
\lefteqn{
|u[\partial\Omega,K,\mu](x')-u[\partial\Omega,K,\mu](x'')|
}
\\ \nonumber
&&\qquad
\leq   \|\mu\|_{  L^{\infty}(\partial\Omega)  }
\biggl\{\biggr.
\int_{{\mathbb{B}}_{n}(x',2|x'-x''|)\cap\partial\Omega}|K(x',y)|\,d\sigma_{y}
 \\
\nonumber
&&\qquad\quad
+
 \int_{{\mathbb{B}}_{n}(x'',3|x'-x''|)\cap\partial\Omega}|K(x'',y)|\,d\sigma_{y}
 \\
\nonumber
&&\qquad\quad
+ 
\int_{\partial\Omega\setminus {\mathbb{B}}_{n}(x',2|x'-x''|)}
 |\, K(x',y)-K(x'',y)  \,|
 \,d\sigma_{y}\biggl.\biggr\} \,.
\end{eqnarray} 
Then Lemma \ref{comin} (ii) implies that
\begin{eqnarray}
\label{k0a3}
\lefteqn{
\int_{{\mathbb{B}}_{n}(x',2|x'-x''|)\cap\partial\Omega}|K(x',y)|\,d\sigma_{y}
 +
 \int_{{\mathbb{B}}_{n}(x'',3|x'-x''|)\cap\partial\Omega}|K(x'',y)|\,d\sigma_{y}
}
\\ \nonumber
&&\qquad\qquad\qquad\qquad
\leq
\| K \|_{  {\mathcal{K}}_{ \gamma_{1},\gamma_{2},\gamma_{3}   }(\partial\Omega)  }
\biggl\{\biggr.
\int_{    {\mathbb{B}}_{n}(x',2|x'-x''|)\cap\partial\Omega     }
\frac{  d\sigma_{y}   }{    |x'-y|^{\gamma_{1}}    }
\\ \nonumber
&&\qquad\qquad\qquad\quad\qquad
+
 \int_{   {\mathbb{B}}_{n}(x'',3|x'-x''|)\cap\partial\Omega    }
 \frac{  d\sigma_{y}  }{   |x''-y|^{\gamma_{1}}  }
 \biggl.\biggr\}
\\ \nonumber
&&\qquad\qquad\qquad\qquad
\leq
\| K \|_{  {\mathcal{K}}_{ \gamma_{1},\gamma_{2},\gamma_{3}   }(\partial\Omega)  }2c''_{\Omega,\gamma_{1}}|x'-x''|^{(n-1)-\gamma_{1}}\,.
\end{eqnarray} 
Moreover, we have 
\begin{eqnarray}
\label{k0a4}
\lefteqn{
\int_{\partial\Omega\setminus {\mathbb{B}}_{n}(x',2|x'-x''|)}
 |\, K(x',y)-K(x'',y)  \,|
 \,d\sigma_{y}
 }
 \\
\nonumber
&&\qquad
\leq
\| K \|_{  {\mathcal{K}}_{ \gamma_{1},\gamma_{2},\gamma_{3}   }(\partial\Omega)  }
\int_{\partial\Omega\setminus {\mathbb{B}}_{n}(x',2|x'-x''|)}
\frac{  |x'-x''|^{\gamma_{3} }   }{  |x'-y|^{\gamma_{2}}   }\,d\sigma_{y}
\end{eqnarray} 
both in case $\gamma_{2}\in ]n-1,+\infty[$ and $\gamma_{2}=n-1$ and for all $\gamma_{3}\in]0,1]$. 

Under the assumptions of statement (ii), Lemma \ref{comin} (iii) implies that 
\begin{equation}
\label{k0a5}
\int_{\partial\Omega\setminus {\mathbb{B}}_{n}(x',2|x'-x''|)}
\frac{  |x'-x''|^{\gamma_{3} }   }{  |x'-y|^{\gamma_{2}}   }\,d\sigma_{y}
\leq c'''_{\Omega,\gamma_{2}}|x'-x''|^{(n-1)-\gamma_{2}+\gamma_{3}}\,.
\end{equation}
Instead, under the assumptions of statement (iii), Lemma \ref{comin} (iv) implies that 
\begin{equation}
\label{k0a6}
\int_{\partial\Omega\setminus {\mathbb{B}}_{n}(x',2|x'-x''|)}
\frac{  |x'-x''|^{\gamma_{3} }   }{  |x'-y|^{\gamma_{2}}   }\,d\sigma_{y}
\leq c^{iv}_{\Omega }
 |x'-x''|^{\gamma_{3} }  
 |\,  \ln  |x'-x''|\,|  
\,.
\end{equation}
Then inequalities (\ref{k0a1})--(\ref{k0a6}) imply the validity of statements (ii), (iii). \hfill  $\Box$ 

\vspace{\baselineskip}

We note that Proposition \ref{k0a} (ii)  for $n=3$, $\gamma_{1}=2-\alpha$, $\gamma_{2}=3-\alpha$, $\gamma_{3}=1$ and for $K$ fixed is known (see Kirsch and Hettlich \cite[\S\ 3.1.3, Thm.~3.17 (a)]{KiHe15}.) Next we introduce two technical lemmas, which we need to define an auxiliary integral operator. 
\begin{lem}
\label{k0b}
Let $\Omega$ be a bounded open Lipschitz subset of ${\mathbb{R}}^{n}$. Let $\alpha,\beta\in]0,1[$. Let $\gamma_{2}\in {\mathbb{R}}$, $\gamma_{3}\in]0,1]$.

If $\gamma_{2}-\beta>n-1$, we further require that $\gamma_{3}+(n-1)-(\gamma_{2}-\beta)>0$.

Then there exists $c>0$ such that the function $u[\partial\Omega,K,\mu]$ defined by (\ref{k0a0}) satisfies the following inequality
\begin{eqnarray}
\label{k0b1}
\lefteqn{
| u[\partial\Omega,K,\mu](x')-u[\partial\Omega,K,\mu](x'')|
}
\\ \nonumber
&&\qquad
\leq
c\| K \|_{  {\mathcal{K}}_{ (n-1)-\alpha,\gamma_{2},\gamma_{3}  }(\partial\Omega)  } \|\mu\|_{  C^{0,\beta}(\partial\Omega)  }\omega(|x'-x''|)
\\ \nonumber
&&\qquad\quad
+\|\mu\|_{  C^{0}(\partial\Omega)  }
| u[\partial\Omega,K,1](x')-u[\partial\Omega,K,1](x'')|
\qquad\forall x',x''\in\partial\Omega\,,
\end{eqnarray}
for all  $(K,\mu)\in {\mathcal{K}}_{(n-1)-\alpha,\gamma_{2},\gamma_{3}  }(\partial\Omega) \times  C^{0,\beta}(\partial\Omega)$, where 
\[
\omega(r)\equiv \left\{
\begin{array}{ll}
r^{\min\{\alpha+\beta,\gamma_{3}\}}  & {\mathrm{if}}\ \gamma_{2}-\beta<n-1\,,
\\
\max\{
r^{\alpha+\beta},\omega_{\gamma_{3}}(r)\}  & {\mathrm{if}}\ \gamma_{2}-\beta=n-1\,,
 \\
r^{\min\{
\alpha+\beta, \gamma_{3}+(n-1)-(\gamma_{2}-\beta) 
\} }& {\mathrm{if}}\ \gamma_{2}-\beta>n-1\,,
\end{array}
\right.
\qquad\forall r\in]0,+\infty[\,.
\]
\end{lem}
{\bf Proof.} By  Remark \ref{om4}  and by Proposition \ref{k0a} (i), it suffices to consider case $0<|x'-x''|<r_{\gamma_{3}}$. By the triangular inequality, we have
\begin{eqnarray}
\label{k0b2}
\lefteqn{
| u[\partial\Omega,K,\mu](x')-u[\partial\Omega,K,\mu](x'')|
}
\\\ \nonumber
&&\qquad
\leq
\left|\int_{\partial\Omega}
[K(x',y)-K(x'',y)]
(\mu(y)-\mu(x'))\,d\sigma_{y}\right|
\\ \nonumber
&&\qquad\quad
+|\mu(x')|\left|\int_{\partial\Omega}
[K(x',y)-K(x'',y)]\,d\sigma_{y}\right|\,.
\end{eqnarray}
By the inclusion ${\mathbb{B}}_{n}(x',2|x'-x''|)\subseteq {\mathbb{B}}_{n}(x'',3|x'-x''|)$, and by the
triangular inequality, and by Lemmas \ref{rec} (i),   \ref{comin} (ii),  and by the inequality
 \[
 |y-x'|^{\beta}\leq  |y-x''|^{\beta}  +  |x'-x''|^{\beta}\,,
 \]
we have
\begin{eqnarray}
\label{k0b3}
\lefteqn{
\left|\int_{\partial\Omega}
[K(x',y)-K(x'',y)]
(\mu(y)-\mu(x'))\,d\sigma_{y}\right|
}
\\ \nonumber
&& 
\leq
\int_{{\mathbb{B}}_{n}(x',2|x'-x''|)\cap\partial\Omega}|K(x',y)|
\,|y-x' |^{\beta}
\,d\sigma_{y}\|\mu\|_{C^{0,\beta}(\partial\Omega)} 
\\ \nonumber
&&\quad
+\int_{{\mathbb{B}}_{n}(x'',3|x'-x''|)\cap\partial\Omega}|K(x'',y)|
\, |y-x'|^{\beta}
\,d\sigma_{y}\|\mu\|_{C^{0,\beta}(\partial\Omega)} 
\\ \nonumber
&&\quad
+ \int_{\partial\Omega\setminus {\mathbb{B}}_{n}(x',2|x'-x''|) }|K(x',y)-K(x'',y)|\,|y-x'|^{\beta}\,d\sigma_{y}\|\mu\|_{C^{0,\beta}(\partial\Omega)} 
\\ \nonumber
&& 
\leq
\| K \|_{  {\mathcal{K}}_{ (n-1)-\alpha,\gamma_{2},\gamma_{3}  }(\partial\Omega)  } \|\mu\|_{  C^{0,\beta}(\partial\Omega)  }
\\ \nonumber
&&\quad
\times
\biggl\{\biggr.
\int_{{\mathbb{B}}_{n}(x',2|x'-x''|)\cap\partial\Omega}
\frac{  d\sigma_{y}      }{
|y-x' |^{(n-1)-(\alpha+\beta)}
}
\\  \nonumber
&&\quad
+\int_{{\mathbb{B}}_{n}(x'',3|x'-x''|)\cap\partial\Omega}
\frac{|x'-x''|^{\beta}\,d\sigma_{y}      }{
|y-x'' |^{(n-1)-\alpha}
}
\\  \nonumber
&&\quad
+\int_{{\mathbb{B}}_{n}(x'',3|x'-x''|)\cap\partial\Omega}
\frac{ d\sigma_{y}      }{
|y-x'' |^{(n-1)-(\alpha+\beta)}
}
\\ \nonumber
&&\quad
+ \int_{\partial\Omega\setminus {\mathbb{B}}_{n}(x',2|x'-x''|) }
\frac{|x'-x''|^{\gamma_{3} }|x'-y|^{\beta}\,d\sigma_{y} 
}{|x'-y|^{ \gamma_{2} }}
\biggl.\biggr\}
\\ \nonumber
&& 
\leq
\| K \|_{  {\mathcal{K}}_{ (n-1)-\alpha,\gamma_{2},\gamma_{3}   }(\partial\Omega)  } \|\mu\|_{  C^{0,\beta}(\partial\Omega)  }
\\ \nonumber
&&\quad \times
\biggl\{\biggr.
2 c''_{\Omega, (n-1)-(\alpha+\beta)}
|x'-x''|^{\alpha+\beta}
+|x'-x''|^{\beta}c''_{\Omega, (n-1)- \alpha }|x'-x''|^{\alpha}
\\ \nonumber
&& \quad 
+|x'-x''|^{ \gamma_{3}}\int_{\partial\Omega\setminus {\mathbb{B}}_{n}(x',2|x'-x''|) }
\frac{ \,d\sigma_{y} 
}{|x'-y|^{ \gamma_{2}-\beta }}
\biggl.\biggr\}\,.
\end{eqnarray}
At this point we distinguish three cases. If $\gamma_{2}-\beta<n-1$, then Lemma \ref{comin} (i) implies that
\[
\int_{\partial\Omega\setminus {\mathbb{B}}_{n}(x',2|x'-x''|) }
\frac{ \,d\sigma_{y} 
}{|x'-y|^{ \gamma_{2}-\beta }}
\leq
\int_{\partial\Omega  }
\frac{ \,d\sigma_{y} 
}{|x'-y|^{ \gamma_{2}-\beta }}\leq c'_{\Omega, \gamma_{2}-\beta}\,,
\]
and thus inequalities (\ref{k0b2}) and (\ref{k0b3}) imply that there exists $c>0$ such that inequality (\ref{k0b1}) holds with $\omega(r)=r^{ \min\{\alpha+\beta,\gamma_{3}\} }$. 
 If $\gamma_{2}-\beta=n-1$, then Lemma \ref{comin} (iv) implies that
\[
\int_{\partial\Omega\setminus {\mathbb{B}}_{n}(x',2|x'-x''|) }
\frac{ \,d\sigma_{y} 
}{|x'-y|^{ \gamma_{2}-\beta }}
\leq c_{\Omega}^{iv}\left| \ln  |x'-x'' | \right|\,,
\]
and thus inequalities (\ref{k0b2}) and (\ref{k0b3}) imply that there exists $c>0$ such that inequality (\ref{k0b1}) holds with $\omega(r)=\max\{
r^{\alpha+\beta},\omega_{\gamma_{3}}(r)\} $. If $\gamma_{2}-\beta>n-1$, then Lemma \ref{comin} (iii) implies that
\[
\int_{\partial\Omega\setminus {\mathbb{B}}_{n}(x',2|x'-x''|) }
\frac{ \,d\sigma_{y} 
}{|x'-y|^{ \gamma_{2}-\beta }}
\leq c'''_{\Omega,\gamma_{2}-\beta}|x'-x'' |^{(n-1)-(\gamma_{2}-\beta)}\,,
\]
and thus inequalities (\ref{k0b2}) and (\ref{k0b3}) imply that there exists $c>0$ such that inequality (\ref{k0b1}) holds with $\omega(r)= 
r^{\min\{ \alpha+\beta ,  \gamma_{3}+(n-1)-(\gamma_{2}-\beta)
\}}$. \hfill  $\Box$ 

\vspace{\baselineskip}

We also point out the validity of the following `folklore' Lemma
\begin{lem}
\label{k0c}
Let $\Omega$ be a bounded open Lipschitz subset of ${\mathbb{R}}^{n}$. Let $\gamma_{1}\in ]-\infty,n-1[$. Let $G$ be a subset of ${\mathbb{R}}^{n}$. Let $K\in C^{0}((G\times\partial\Omega)\setminus\Delta_{\partial\Omega})$ be such that
\[
\kappa_{  \gamma_{1}   }\equiv\sup _{(x,y)\in (G\times\partial\Omega)\setminus\Delta_{\partial\Omega}  }
|x-y|^{\gamma_{1}}|\,K(x,y)\,|
<+\infty\,.
\]
Let $\mu\in L^{\infty}(\partial\Omega) $. Then the function $K(x,\cdot)\mu(\cdot)$ is integrable in $\partial\Omega$ for all $x\in G$ and the function $u^{\sharp}[\partial\Omega,K,\mu]$ from $G$ to ${\mathbb{C}}$ defined by
\[
u^{\sharp}[\partial\Omega,K,\mu](x)\equiv
\int_{\partial\Omega}K(x,y)\mu (y)\,d\sigma_{y}\qquad\forall x\in G\,,
\]
is continuous. If $\sup_{x\in G}\int_{\partial\Omega} \frac{d\sigma_{y}}{|x-y|^{\gamma_{1}}}<\infty$, then $u^{\sharp}[\partial\Omega,K,\mu]$ satisfies the inequality
\begin{equation}
\label{k0c3}
|u^{\sharp}[\partial\Omega,K,\mu](x)|
\leq
\sup_{x\in G}\int_{\partial\Omega} \frac{d\sigma_{y}}{|x-y|^{\gamma_{1}}} \kappa_{  \gamma_{1}   }
\|\mu\|_{   L^{\infty}(\partial\Omega)   }\qquad\forall x\in G\,.
\end{equation}
\end{lem}
{\bf Proof.} The integrability of $K(x,\cdot)\mu(\cdot)$ follows by the inequality
\[
\left|
K(x,y)\mu (y)
\right|\leq \frac{    \kappa_{  \gamma_{1} } \|\mu\|_{L^{\infty}(\partial\Omega) }  }{|x-y|^{\gamma_{1}}}
\qquad {\mathrm{a.a.}}\ y\in\partial\Omega\,.
\]
Since $\sup_{x\in G}\int_{\partial\Omega} \frac{d\sigma_{y}}{|x-y|^{\gamma_{1}}}<\infty$, 
 inequality (\ref{k0c3})  follows 
and the Vitali Convergence Theorem implies that $u^{\sharp}[\partial\Omega,K,\mu]$ is continuous on $G$ (cf.~\textit{e.g.}, Folland~\cite[(2.33) p.~60, p.~180]{Fo84}.) \hfill  $\Box$ 

\vspace{\baselineskip}

We now introduce an auxiliary integral operator, and we deduce some    properties which we need in the sequel  by applying  Proposition \ref{k0a} and Lemma \ref{k0b}. 
\begin{lem}
\label{qk}
 Let $\theta\in]0,1]$.   Let $\Omega$ be a bounded open Lipschitz subset of ${\mathbb{R}}^{n}$. Then the following statements hold.
\begin{enumerate}
\item[(i)] Let $Z\in C^{0}(  ({\mathrm{cl}}\Omega\times\partial\Omega)\setminus\Delta_{\partial\Omega}   )$ satisfy inequality
\begin{equation}
\label{qk1}
\kappa_{n-1}[Z]\equiv \sup_{(x,y)\in  ({\mathrm{cl}}\Omega\times\partial\Omega)\setminus\Delta_{\partial\Omega} }
|x-y|^{n-1}|Z(x,y)|<+\infty\,.
\end{equation}
Let $(f,\mu)\in C^{0,\theta}({\mathrm{cl}}\Omega)\times L^{\infty}(\partial\Omega)$. Let $H^{\sharp}[Z,f]$ be the  function from $({\mathrm{cl}}\Omega\times\partial\Omega)\setminus\Delta_{\partial\Omega}  $ to ${\mathbb{C}}$ defined by 
\[
H^{\sharp}[Z,f](x,y)\equiv (f(x)-f(y))Z(x,y)\qquad\forall (x,y)\in   ({\mathrm{cl}}\Omega\times\partial\Omega)\setminus\Delta_{\partial\Omega}  \,.
\]
 If $x\in {\mathrm{cl}}\Omega$, then the function $H^{\sharp}[Z,f](x,\cdot)$ is Lebesgue integrable in $\partial\Omega$ and the function $Q^{\sharp}[Z,f,\mu]$ 
from ${\mathrm{cl}}\Omega$ to ${\mathbb{C}}$  defined by 
\[
Q^{\sharp}[Z,f,\mu](x)\equiv \int_{\partial\Omega}H^{\sharp}[Z,f](x,y)\mu(y)\,d\sigma_{y}\qquad\forall x\in {\mathrm{cl}}\Omega\,,
\]
is continuous. 
\item[(ii)] The map $H$ from ${\mathcal{K}}_{n-1,n,1}(\partial\Omega)\times C^{0,\theta}(\partial\Omega)$ to ${\mathcal{K}}_{n-1-\theta,n-1,\theta}(\partial\Omega)$, which takes $(Z,g)$ to the function from $(\partial\Omega)^{2}\setminus\Delta_{\partial\Omega}$ to ${\mathbb{C}}$  defined by 
\[
H[Z,g](x,y)\equiv (g(x)-g(y))Z(x,y)\qquad\forall (x,y)\in
(\partial\Omega)^{2}\setminus\Delta_{\partial\Omega}\,,
\]
is bilinear and continuous. 
\item[(iii)] The map $Q$ from ${\mathcal{K}}_{n-1,n,1}(\partial\Omega)\times C^{0,\theta}(\partial\Omega)\times L^{\infty}(\Omega)$ to 
$C^{0,\omega_{\theta}(\cdot)}(\partial\Omega)$, which takes $(Z,g,\mu)$ to the function from $\partial\Omega$ to ${\mathbb{C}}$ defined by 
\[
Q[Z,g,\mu](x)\equiv   \int_{\partial\Omega}
H[Z,g](x,y)\mu(y)\,d\sigma_{y}\qquad\forall x\in\partial\Omega\,,
\]
is trilinear and continuous. 
\item[(iv)] Let $\alpha\in]0,1[$, $\beta\in]0,1]$.  Then there exists $q\in]0,+\infty[$ such that
\begin{eqnarray*}
\lefteqn{
| Q[Z,g,\mu](x')-Q[Z,g,\mu](x'')|
}
\\ \nonumber
&&\quad
\leq
q\|  Z  \|_{  {\mathcal{K}}_{ n-1 ,n,1  }(\partial\Omega)  }\|g\|_{C^{0,\alpha}(\partial\Omega)} \|\mu\|_{  C^{0,\beta}(\partial\Omega)  }|x'-x''|^{\alpha}
\\ \nonumber
&&\quad\quad
+ \|\mu\|_{  C^{0}(\partial\Omega)  }
| Q[Z,g,1](x')-Q[Z,g,1](x'')|
\quad\forall x',x''\in\partial\Omega\,,
\end{eqnarray*}
for all  $(Z,g,\mu)\in {\mathcal{K}}_{n-1,n,1   }(\partial\Omega)\times  C^{0,\alpha}(\partial\Omega) \times  C^{0,\beta}(\partial\Omega)$.
\end{enumerate}
\end{lem}
{\bf Proof.} By assumption (\ref{qk1}), and by the H\"{o}lder continuity of $f$, we have 
\[
\left|    H^{\sharp}[Z,f](x,y)   \right|
\leq
\frac{  |f|_{\theta}  }{|x-y|^{(n-1)-\theta}}
\kappa_{n-1}[Z]\,,
\]
for all $ (x,y) \in ({\mathrm{cl}}\Omega\times\partial\Omega)\setminus\Delta_{\partial\Omega}$. Then Lemma \ref{k0c} implies the validity of statement (i). 

By the H\"{o}lder continuity of $g$, we have 
\begin{equation}
\label{qk8}
\left|    H [Z,g](x,y)   \right|
\leq
\frac{  |g|_{\theta}  }{|x-y|^{(n-1)-\theta}}\|Z\|_{
{\mathcal{K}}_{n-1,n,1   }(\partial\Omega)
}
\qquad
\forall (x,y) \in (\partial\Omega)^{2}\setminus\Delta_{\partial\Omega}\,.
\end{equation}
Now let $x'$, $x''\in\partial\Omega$, $x'\neq x''$, $y\in \partial\Omega\setminus{\mathbb{B}}_{n}(x',2|x'-x''|)$. Then we have 
\begin{eqnarray}
\label{qk9}
\lefteqn{
|H[Z,g](x',y)-H[Z,g](x'',y)|
}
\\ \nonumber
&&\qquad
\leq 
|g(x')-g(y)|\,|Z(x',y)-Z(x'',y)|
+|g(x')-g(x'')|\, |Z(x'',y)|
\\ \nonumber
&&\qquad
\leq
\|g\|_{C^{0,\theta}(\partial\Omega) }\|Z\|_{
{\mathcal{K}}_{n-1,n,1   }(\partial\Omega)}
\left\{
\frac{  |x'-y|^{\theta} |x'-x''|   }{|x'-y|^{n}}
+
\frac{  |x'-x''|^{\theta}   }{   |x''-y|^{n-1}  }
\right\}\,.
\end{eqnarray}
Since $ |x'-x''|\leq  |x'-y|$, we have $ |x'-x''|^{1-\theta}\leq  |x'-y|^{1-\theta}$. Moreover, Lemma \ref{rec} (i) implies that $|x''-y|\geq \frac{1}{2}|x'-y|$ and thus the term in braces in the right hand side of (\ref{qk9}) is less or equal to
\begin{equation}
\label{qk10}
\frac{  |x'-y|\, |x'-x''|^{\theta}  }{|x'-y|^{n}}
+
\frac{2^{n-1}  |x'-x''|^{\theta}   }{   |x'-y|^{n-1}  }
\leq
(1+2^{n-1})\frac{|x'-x''|^{\theta}}{  |x'-y|^{n-1}  }\,.
\end{equation}
Hence, inequalities (\ref{qk8})--(\ref{qk10}) imply that
\begin{equation}
\label{qk11}
\|H[Z,g]  \|_{  {\mathcal{K}}_{n-1-\theta,n-1,\theta}(\partial\Omega)   }
\leq 2^{n}\|Z\|_{ {\mathcal{K}}_{n-1,n,1}(\partial\Omega) }\|g\|_{ C^{0,\theta}(\partial\Omega) }\,.
\end{equation}
Hence statement (ii) holds true. We now turn to prove  (iii). By Proposition \ref{k0a} (iii) with $\gamma_{1}=n-1-\theta$, $\gamma_{2}=n-1$, $\gamma_{3}=\theta$, the map $u[\partial\Omega,\cdot,\cdot]$ is continuous from ${\mathcal{K}}_{n-1-\theta,n-1,\theta}(\partial\Omega)\times L^{\infty}(\partial\Omega)$ to 
$C^{0, \max\{r^{(n-1)-[(n-1)-\theta]},\omega_{\theta}(r)  \} }(\partial\Omega)
=C^{0,\omega_{\theta}(\cdot) }(\partial\Omega)$.    Then statement (ii) implies that $u[\partial\Omega,H[\cdot,\cdot],\cdot]$ is continuous from ${\mathcal{K}}_{n-1,n,1}(\partial\Omega)\times 
C^{0, \theta  }(\partial\Omega)\times L^{\infty}(\partial\Omega)$ to $C^{0,\omega_{\theta}(\cdot) }(\partial\Omega)$. Since 
\begin{equation}
\label{qk12}
u[\partial\Omega,H[Z,g],\mu]
=\int_{\partial\Omega}H[Z,g](x,y)\mu(y)\,d\sigma_{y}\qquad\forall x\in\partial\Omega\,,
\end{equation}
statement (iii) holds true. Since $C^{0,\beta_{1}}(\partial\Omega)$ is continuously imbedded into $C^{0,\beta_{2}}(\partial\Omega)$ whenever $0<\beta_{2}\leq\beta_{1}\leq 1$, then we can
 assume that $\alpha+\beta<1$. Then by equality (\ref{qk12}) and by Lemma \ref{k0b} with   $\gamma_{2}=n-1$, $\gamma_{3}=\alpha$ and by statement (ii) with $\theta=\alpha$, statement (iv)  holds true.  \hfill  $\Box$ 

\vspace{\baselineskip}

 \section{Preliminaries on layer potentials}
 \label{prelat}

Let ${\mathbf{a}}$ be as in (\ref{introd0}), (\ref{ellip}). Let $S_{ {\mathbf{a}} }$ be a fundamental solution of $P[{\mathbf{a}},D]$. Let $\Omega$ be a bounded open  Lipschitz subset of ${\mathbb{R}}^{n}$. If $\mu\in L^{\infty}(\partial\Omega)$, Lemma \ref{fun} (i) ensures the convergence of the integral 
\[
v[\partial\Omega           ,S_{ {\mathbf{a}} },\mu](x)\equiv 
\int_{\partial\Omega}S_{ {\mathbf{a}} }(x-y)\mu(y)\,d\sigma_{y}
\qquad\forall x\in{\mathbb{R}}^{n}\,,
\]
which defines the single layer potential relative to $\mu$, $S_{ {\mathbf{a}} }$.  
We collect in the following statement  some known properties of the single layer potential which we exploit in the sequel (cf. Miranda~\cite{Mi65}, Wiegner~\cite{Wi93}, Dalla Riva \cite{Da13}, Dalla Riva, Morais and Musolino \cite{DaMoMu13} and references therein.)

\vspace{\baselineskip}

\begin{thm}
\label{slay}
Let ${\mathbf{a}}$ be as in (\ref{introd0}), (\ref{ellip}). Let $S_{ {\mathbf{a}} }$ be a fundamental solution of $P[{\mathbf{a}},D]$. 
Let  $\alpha\in]0,1[$, $m\in {\mathbb{N}}\setminus\{0\}$. Let $\Omega$ be a bounded open subset of ${\mathbb{R}}^{n}$ of class $C^{m,\alpha}$. Then the following statements hold. 
\begin{enumerate}
\item[(i)] If $\mu\in C^{m-1,\alpha}(\partial\Omega)$, then the function 
$v^{+}[\partial\Omega        ,S_{{\mathbf{a}}} ,\mu]\equiv v[\partial\Omega, S_{{\mathbf{a}}},\mu]_{|{\mathrm{cl}}\Omega}$ belongs to $C^{m,\alpha}({\mathrm{cl}}\Omega)$ and  the function 
$v^{-}[\partial\Omega ,S_{{\mathbf{a}}} ,\mu]\equiv v[\partial\Omega, S_{{\mathbf{a}}},\mu]_{|{\mathrm{cl}}\Omega^{-}}$ belongs to $C^{m,\alpha}_{{\mathrm{loc}}}({\mathrm{cl}}\Omega^{-})$.
 Moreover the map which takes $\mu$ to the function  $v^{+}[\partial\Omega
 ,S_{{\mathbf{a}}},\mu]$ is  continuous from   $C^{m-1,\alpha}(\partial\Omega)$ to $C^{m,\alpha}({\mathrm{cl}}\Omega)$ and the map from  the space  $C^{m-1,\alpha}(\partial\Omega)$ to $C^{m,\alpha}({\mathrm{cl}}{\mathbb{B}}_{n}(0,R)\setminus \Omega) $ which takes $\mu$ to $v^{-}[\partial\Omega ,S_{{\mathbf{a}}},\mu]_{|{\mathrm{cl}}{\mathbb{B}}_{n}(0,R)\setminus \Omega}$ is  continuous for all 
$R\in]0,+\infty[$ such that ${\mathrm{cl}}\Omega\subseteq {\mathbb{B}}_{n}(0,R)$.
\item[(ii)] Let $l\in\{1,\dots,n\}$. If $\mu\in C^{0,\alpha}(\partial\Omega)$, then 
we have the following jump relation
\begin{eqnarray*}
\lefteqn{
\frac{\partial}{\partial x_{l}}v^{\pm}[\partial\Omega , S_{{\mathbf{a}}},\mu](x)
}
\\
&&
=
\mp\frac{\nu_{l}(x)}{
2\nu(x)^{t}a^{(2)}\nu(x)
}\mu (x)
+
\int_{\partial\Omega} \partial_{x_{l}}S_{{\mathbf{a}}}(x-y)\mu(y)\,d\sigma_{y}\qquad\forall x\in \partial\Omega
\,,
\end{eqnarray*}
where the integral in the right hand side exists in the sense of the principal value. 
\end{enumerate}
\end{thm}

Then we introduce the following refinement of a classical result 
for homogeneous second order elliptic operators
(cf. Miranda \cite{Mi70}.)
\begin{thm}
\label{v0a} 
Let ${\mathbf{a}}$ be as in (\ref{introd0}), (\ref{ellip}). Let $S_{ {\mathbf{a}} }$ be a fundamental solution of $P[{\mathbf{a}},D]$. 
  Let $\Omega$ be a bounded open  Lipschitz subset of ${\mathbb{R}}^{n}$.  Let $\gamma\in]0,1[$. Then the operator $v[\partial\Omega,  S_{{\mathbf{a}}},\cdot]_{|\partial\Omega}$
from $L^{\infty}(\partial\Omega)$ to $C^{0,\gamma  }(\partial\Omega)$ which takes $\mu $ to $v[\partial\Omega,  S_{{\mathbf{a}}},\mu]_{|\partial\Omega}$
is continuous.

If we further assume that $n>2$, then $v[\partial\Omega,  S_{{\mathbf{a}}},\cdot]_{|\partial\Omega}$ is continuous from $L^{\infty}(\partial\Omega)$ to $C^{0,\omega_{1}(\cdot)  }(\partial\Omega)$.
\end{thm}
{\bf Proof.} By Lemma \ref{fun}, we have that $S_{{\mathbf{a}}}(x-y)\in {\mathcal{K}}_{(n-1)-\gamma,n-1,1}(\partial\Omega)$, and also  that   $S_{{\mathbf{a}}}(x-y)\in {\mathcal{K}}_{ n-2,n-1,1}(\partial\Omega)$ if we further assume that $n>2$. 
Since
\[
v[\partial\Omega, S_{{\mathbf{a}}}, \mu]_{|\partial\Omega}=
u[\partial\Omega,S_{{\mathbf{a}}}(x-y),\mu]\,,
\]
Proposition \ref{k0a} (iii) implies that $v[\partial\Omega, S_{{\mathbf{a}}}, \cdot]$ is continuous from $L^{\infty}(\partial\Omega)$ to 
$C^{0,\max\{  r^{\gamma},\omega_{1}(r)   \}   }(\partial\Omega)=C^{0,\gamma  }(\partial\Omega)$, and that  $v[\partial\Omega, S_{{\mathbf{a}}}, \cdot]$ is continuous from $L^{\infty}(\partial\Omega)$ to 
$C^{0,\max\{  r ,\omega_{1}(r)   \}   }(\partial\Omega)=C^{0,\omega_{1}(r)  }(\partial\Omega)$
if we further assume that $n>2$.  \hfill  $\Box$ 

\vspace{\baselineskip}

Next we turn to the double layer potential and we introduce the following technical result  (cf. Miranda~\cite{Mi65}, Wiegner~\cite{Wi93}, Dalla Riva \cite{Da13}, Dalla Riva, Morais and Musolino \cite{DaMoMu13} and references therein.) 
\begin{thm}
\label{dlay}
Let ${\mathbf{a}}$ be as in (\ref{introd0}), (\ref{ellip}). Let $S_{ {\mathbf{a}} }$ be a fundamental solution of $P[{\mathbf{a}},D]$. 
Let  $\alpha\in]0,1[$, $m\in {\mathbb{N}}\setminus\{0\}$. Let $\Omega$ be a bounded open subset of ${\mathbb{R}}^{n}$ of class $C^{m,\alpha}$. Then the following statements hold. 
\begin{enumerate}
\item[(i)] If $\mu\in C^{0,\alpha}(\partial\Omega)$, then the restriction 
$w[\partial\Omega, {\mathbf{a}},S_{{\mathbf{a}}},\mu]_{|\Omega}$ can be extended uniquely to a continuous function $w^{+}[\partial\Omega ,{\mathbf{a}},S_{{\mathbf{a}}},\mu]$ from ${\mathrm{cl}}\Omega$ to ${\mathbb{C}}$, and 
$w[\partial\Omega ,{\mathbf{a}},S_{{\mathbf{a}}},\mu]_{|\Omega^{-}}$ can be extended uniquely to a continuous function $w^{-}[\partial\Omega ,{\mathbf{a}},S_{{\mathbf{a}}},\mu]$ from ${\mathrm{cl}}\Omega^{-}$ to ${\mathbb{C}}$ and we have the following jump relation
\[
w^{\pm}[\partial\Omega ,{\mathbf{a}},S_{{\mathbf{a}}},\mu](x)
=\pm\frac{1}{2}\mu(x)+w[\partial\Omega ,{\mathbf{a}},S_{{\mathbf{a}}},\mu](x)
\qquad\forall x\in\partial\Omega\,.
\]
\item[(ii)] If $\mu\in C^{m,\alpha}(\partial\Omega)$, then $w^{+}[\partial\Omega ,{\mathbf{a}},S_{{\mathbf{a}}},\mu]$ belongs to $ C^{m,\alpha}({\mathrm{cl}}\Omega)$ 
and $w^{-}[\partial\Omega ,{\mathbf{a}},S_{{\mathbf{a}}},\mu] $ belongs to $
C^{m,\alpha}_{ {\mathrm{loc}} }({\mathrm{cl}}\Omega^{-})$. Moreover, the map 
from  the space $C^{m,\alpha}(\partial\Omega)$ to $C^{m,\alpha}({\mathrm{cl}} \Omega)$
 which takes $\mu$ to $
w^{+}[\partial\Omega ,{\mathbf{a}},S_{{\mathbf{a}}},\mu]$ is  continuous and the map from   the space $C^{m,\alpha}(\partial\Omega)$ to $C^{m,\alpha}({\mathrm{cl}}{\mathbb{B}}_{n}(0,R)\setminus \Omega) $ which takes $\mu$ to 
$w^{-}[\partial\Omega ,{\mathbf{a}},S_{{\mathbf{a}}},\mu]_{|{\mathrm{cl}}{\mathbb{B}}_{n}(0,R)\setminus\Omega}$ is continuous for all $R\in]0,+\infty[$ such that ${\mathrm{cl}}\Omega\subseteq {\mathbb{B}}_{n}(0,R)$.
\item[(iii)] Let $r\in\{1,\dots,n\}$. If $\mu\in C^{m,\alpha}(\partial\Omega)$ and $U$ is an open neighborhood of $\partial\Omega$ in ${\mathbb{R}}^{n}$ and $\tilde{\mu}\in C^{m}(U)$, $\tilde{\mu}_{|\partial\Omega}=\mu$, then the following equality holds
\begin{eqnarray}
\label{dlay0}
\lefteqn{\frac{\partial}{\partial x_{r}}w[\partial\Omega ,{\mathbf{a}},S_{{\mathbf{a}}},\mu](x)}
\\   \nonumber
&&
=\sum_{j,l=1}^{n}a_{lj} \frac{\partial}{\partial x_{l}}
\biggl\{
\int_{\partial\Omega}
S_{{\mathbf{a}}}(x-y)
\biggl[
\nu_{r}(y)\frac{\partial\tilde{\mu}}{\partial y_{j}}(y)
-
\nu_{j}(y)\frac{\partial\tilde{\mu}}{\partial y_{r}}(y)
\biggr]\,d\sigma_{y}
\biggr\}
\\   \nonumber
&&\ \ +
\int_{\partial\Omega}
\biggl[
 DS_{{\mathbf{a}}}(x-y)a^{(1)}+
a S_{{\mathbf{a}}}(x-y)
\biggr]   \nu_{r} (y) \mu(y)\,d\sigma_{y}
\\   \nonumber
&\ &\ \ -\int_{\partial\Omega}\partial_{x_{r}}S_{{\mathbf{a}}}(x-y) 
\nu^{t}(y)a^{(1)}\mu(y)\,d\sigma_{y} \qquad\forall x\in{\mathbb{R}}^{n}\setminus\partial\Omega\,.
\end{eqnarray}
\end{enumerate}
\end{thm}
We note that formula (\ref{dlay0}) for the Laplace operator with $n=3$ can be found in 
G\"{u}nter~\cite[Ch.~2, \S\ 10, (42)]{Gu67}. By combining Theorems \ref{slay} and \ref{dlay}, we deduce that under the assumptions of Theorem \ref{dlay} (iii), the following equality holds
\begin{eqnarray}
\label{dlay1}
\lefteqn{\frac{\partial}{\partial x_{r}}w^{+}[\partial\Omega ,{\mathbf{a}},S_{{\mathbf{a}}},\mu] }
\\  \nonumber
&&
=\sum_{j,l=1}^{n}a_{lj} \frac{\partial}{\partial x_{l}}
v^{+}[\partial\Omega ,S_{{\mathbf{a}}},M_{rj}[\mu]] 
+
Dv^{+}[\partial\Omega ,S_{{\mathbf{a}}},\nu_{r}   \mu] a^{(1)}
\\    \nonumber
&&\ \ 
+av^{+}[\partial\Omega ,S_{{\mathbf{a}}},\nu_{r}   \mu] 
 -\frac{\partial}{\partial x_{r}}v^{+}[\partial\Omega ,S_{{\mathbf{a}}},
 (\nu^{t} a^{(1)})\mu]
\qquad{\mathrm{on}}\ {\mathrm{cl}}
\Omega\,.
\end{eqnarray}
Next we introduce  the following result proved by  Schauder \cite[Hilfsatz VII, p.~112]{Sc31} for the Laplace operator and which we extend here  to second order eliptic operators by exploiting Proposition \ref{k0a}. 
\begin{thm}
\label{w0a}
 Let ${\mathbf{a}}$ be as in (\ref{introd0}), (\ref{ellip}), (\ref{symr}). Let $S_{ {\mathbf{a}} }$ be a fundamental solution of $P[{\mathbf{a}},D]$. 
Let  $\alpha\in]0,1[$. Let $\Omega$ be a bounded open subset of ${\mathbb{R}}^{n}$ of class $C^{1,\alpha}$. If $\mu\in L^{\infty}(\partial\Omega)$, then 
$w[\partial\Omega ,{\mathbf{a}},S_{{\mathbf{a}}},\mu]_{|\partial\Omega}\in C^{0,\alpha}(\partial\Omega)$. Moreover the operator from $L^{\infty}(\partial\Omega)$ to $C^{0,\alpha}(\partial\Omega)$ which takes $\mu$ to $w[\partial\Omega ,{\mathbf{a}},S_{{\mathbf{a}}},\mu]_{|\partial\Omega}$ is continuous. 
\end{thm}
{\bf Proof.} By Lemma \ref{boest}, the function $K_{  {\mathbf{a}}  }(x,y)\equiv 
\overline{B^{*}_{\Omega,y}}\left(S_{{\mathbf{a}}}(x-y)\right)$ belongs to ${\mathcal{K}}_{(n-1)-\alpha,n-\alpha,1}(\partial\Omega)$. Since 
\[
w[\partial\Omega, {\mathbf{a}},S_{{\mathbf{a}}},\mu]_{\partial\Omega}
=
u[\partial\Omega,K_{  {\mathbf{a}}  },\mu]\,,
\]
Proposition \ref{k0a} (ii) implies that $w[\partial\Omega, {\mathbf{a}},S_{{\mathbf{a}}},\cdot]_{|\partial\Omega}$ is continuous from $L^{\infty}(\partial\Omega)$ to $C^{0,\min\{\alpha,
(n-1)-(n-\alpha)+1
\} }(\partial\Omega)=C^{0,\alpha}(\partial\Omega)$. \hfill  $\Box$ 

\vspace{\baselineskip}

\section{Auxiliary integral operators}
\label{aio}

In order to compute the tangential derivatives of the double layer potential, we introduce the following two statements which concern two auxiliary integral operators. To shorten our notation, we define the function $\Theta$ from $({\mathbb{R}}^{n}\times {\mathbb{R}}^{n})\setminus\Delta_{{\mathbb{R}}^{n}}$ to $ {\mathbb{R}}^{n}\setminus\{0\}$ by setting
\begin{equation}
\label{theta}
\Theta(x,y)\equiv x-y\qquad\forall (x,y)\in ({\mathbb{R}}^{n}\times {\mathbb{R}}^{n})\setminus\Delta_{{\mathbb{R}}^{n}}\,.
\end{equation}
\begin{thm}
\label{qrssm}
Let ${\mathbf{a}}$ be as in (\ref{introd0}), (\ref{ellip}), (\ref{symr}). Let $S_{ {\mathbf{a}} }$ be a fundamental solution of $P[{\mathbf{a}},D]$. Let $r\in\{1,\dots,n\}$. Then the following statements hold.
\begin{enumerate}
\item[(i)] Let $\Omega$ be a bounded open Lipschitz subset of ${\mathbb{R}}^{n}$. Let $\theta\in ]0,1]$. If $(f,\mu)\in C^{0,\theta}({\mathrm{cl}}\Omega)\times L^{\infty}(\partial\Omega)$, then the function
\[
Q^{\sharp}[\frac{\partial S_{ {\mathbf{a}} }}{\partial x_{r}} \circ\Theta, f,\mu]   (x)
=
\int_{\partial\Omega}(f(x)-f(y))\frac{\partial S_{ {\mathbf{a}} }}{\partial x_{r}}(x-y)\mu(y)\,d\sigma_{y}\quad\forall x\in {\mathrm{cl}}\Omega\,,
\]
is continuous.
\item[(ii)] Let $\alpha\in]0,1[$, $\beta$, $\theta\in ]0,1]$. Let $m\in {\mathbb{N}}\setminus\{0\}$. Let $\Omega$ be a bounded open subset of ${\mathbb{R}}^{n}$ of class $C^{m,\alpha}$. Then the map $Q^{\sharp}[\frac{\partial S_{ {\mathbf{a}} }}{\partial x_{r}}\circ\Theta ,\cdot,\cdot]$ from 
$C^{m-1,\theta}({\mathrm{cl}}\Omega)\times  C^{m-1,\beta}(\partial\Omega)$ to $C^{m-1,\min\{\alpha,\beta,\theta\}}({\mathrm{cl}}\Omega)$, which takes $(f,\mu)$ to $  Q^{\sharp}[\frac{\partial S_{ {\mathbf{a}} }}{\partial x_{r}} \circ\Theta, f,\mu]  $ is bilinear and continuous. 
\end{enumerate}
\end{thm}
{\bf Proof.} By Lemma \ref{grafun} (ii), statement (i)  is an immediate consequence of  
Lemma   \ref{qk} (i). We now consider statement (ii). By treating  separately cases $x\in\partial\Omega$ and $x\in\Omega$, and by exploiting Theorem \ref{slay} (ii), we have
\[
Q^{\sharp}[\frac{\partial S_{ {\mathbf{a}} }}{\partial x_{r}} \circ\Theta, f,\mu]   (x)
=
f(x)\frac{\partial}{\partial x_{r}}v^{+}[\partial\Omega, S_{ {\mathbf{a}} },\mu](x)
-
\frac{\partial}{\partial x_{r}}v^{+}[\partial\Omega, S_{ {\mathbf{a}} },f\mu](x)
\,,
\]
for all $x\in {\mathrm{cl}}\Omega$. Then the statement follows by Theorem \ref{slay} (i) and by continuity of the pointwise product in Schauder spaces. \hfill  $\Box$ 

\vspace{\baselineskip}

Then we have the following. 
 \begin{thm}
\label{qrs}
 Let ${\mathbf{a}}$ be as in (\ref{introd0}), (\ref{ellip}), (\ref{symr}). Let $S_{ {\mathbf{a}} }$ be a fundamental solution of $P[{\mathbf{a}},D]$. Then 
 the following statement holds.
\begin{enumerate}
\item[(i)] Let $\Omega$ be a bounded open Lipschitz subset of ${\mathbb{R}}^{n}$. Let $\theta\in]0,1]$. Then the bilinear map $
Q[\frac{\partial S_{ {\mathbf{a}} }}{\partial x_{r}}\circ\Theta ,\cdot,\cdot]$ from $C^{0,\theta}(\partial\Omega)\times L^{\infty}(\partial\Omega)$ to $C^{0,\omega_{\theta}(\cdot)}(\partial\Omega)$, which takes $(g,\mu)$ to  the function  
\begin{equation}
\label{qrs0}
Q[\frac{\partial S_{ {\mathbf{a}} }}{\partial x_{r}} \circ\Theta,g,\mu](x)
=\int_{\partial\Omega}(g(x)-g(y))\frac{\partial S_{ {\mathbf{a}} }}{\partial x_{r}}(x-y)\mu(y)\,d\sigma_{y}\quad\forall x\in \partial\Omega\,,
\end{equation}
is continuous. 
\item[(ii)]  Let $\alpha\in]0,1[$, $\beta\in]0,1]$. 
Let $\Omega$ be a bounded open subset of ${\mathbb{R}}^{n}$ of class $C^{1,\alpha}$. Then the bilinear map $Q[\frac{\partial S_{ {\mathbf{a}} }}{\partial x_{r}} \circ\Theta,\cdot,\cdot]$ from $C^{0,\alpha}(\partial\Omega)\times C^{0,\beta}(\partial\Omega)$ to  $C^{0,\alpha}(\partial\Omega)$ which takes $(g,\mu)$ to $Q[\frac{\partial S_{ {\mathbf{a}} }}{\partial x_{r}} \circ\Theta,g,\mu]$ is continuous. 
\end{enumerate}
\end{thm} 
{\bf Proof.} By Lemma \ref{grafun}, we have $\frac{\partial S_{ {\mathbf{a}} }}{\partial x_{r}}\in {\mathcal{K}}_{n-1,n,1}(\partial\Omega)$. Then Lemma \ref{qk} (iii) implies the validity of statement (i). 

We now consider statement (ii). By statement (i) and by the continuity of the inclusion of $C^{0,\beta}(\partial\Omega)$ into $L^{\infty}(\partial\Omega)$, we already know that $Q[\frac{\partial S_{ {\mathbf{a}} }}{\partial x_{r}} \circ\Theta,\cdot,\cdot]$ is continuous from $C^{0,\alpha}(\partial\Omega)\times C^{0,\beta}(\partial\Omega)$ to  $C^{0}(\partial\Omega)$. Then it suffices to show that $Q[\frac{\partial S_{ {\mathbf{a}} }}{\partial x_{r}} \circ\Theta,\cdot,\cdot]$ is continuous from $C^{0,\alpha}(\partial\Omega)\times C^{0,\beta}(\partial\Omega)$ to the semi-normed space  $(C^{0,\alpha}(\partial\Omega),  |\cdot:\partial\Omega|_{\alpha})$. By Lemma \ref{qk} (iv), there exists $q\in]0,+\infty[$ such that
  \begin{eqnarray}
 \label{qrs1}
\lefteqn{
\left|Q[\frac{\partial S_{ {\mathbf{a}} }}{\partial x_{r}}\circ\Theta ,g,\mu](x')-Q[\frac{\partial S_{ {\mathbf{a}} }}{\partial x_{r}} \circ\Theta,g,\mu](x'')\right|
}
\\ \nonumber
&&\quad
\leq 
q\left\| \frac{\partial S_{ {\mathbf{a}} }}{\partial x_{r}} \circ\Theta\right\|_{  {\mathcal{K}}_{ n-1 ,n,1  }(\partial\Omega)  }\|g\|_{C^{0,\alpha}(\partial\Omega)} \|\mu\|_{  C^{0,\beta}(\partial\Omega)  }|x'-x''|^{\alpha}
\\ \nonumber
&&\quad\quad
+\|\mu\|_{  C^{0}(\partial\Omega)  }\left| Q[\frac{\partial S_{ {\mathbf{a}} }}{\partial x_{r}} \circ\Theta,g,1](x')-Q[\frac{\partial S_{ {\mathbf{a}} }}{\partial x_{r}}\circ\Theta ,g,1](x'')\right|\,,
\end{eqnarray}
for all $x',x''\in\partial\Omega$. 
Let $R\in]0,+\infty[$ be such that ${\mathrm{cl}}\Omega\subseteq{\mathbb{B}}_{n}(0,R)$. Let  ` $\tilde{\ }$ '  be an extension operator as in Lemma \ref{extsch}, defined  on $C^{0,\alpha}(\partial\Omega)$. Since
\[
Q[\frac{\partial S_{ {\mathbf{a}} }}{\partial x_{r}} \circ\Theta,g,1](x)
=
Q^{\sharp}[\frac{\partial S_{ {\mathbf{a}} }}{\partial x_{r}}\circ\Theta ,\tilde{g},1](x)\qquad\forall x\in\partial\Omega\,,
\]
Theorem \ref{qrssm} (ii) implies that  $Q[\frac{\partial S_{ {\mathbf{a}} }}{\partial x_{r}} \circ\Theta,\cdot,1]$ is continuous from  $C^{0,\alpha}(\partial\Omega)$ to itself, and that accordingly, there exists $q'\in]0,+\infty[$ such that
\begin{equation}
\label{qrs2}
\|Q[\frac{\partial S_{ {\mathbf{a}} }}{\partial x_{r}} \circ\Theta,g,1]\|_{ C^{0,\alpha}(\partial\Omega) }
\leq q'\|g\|_{ C^{0,\alpha}(\partial\Omega) }
\qquad\forall g\in C^{0,\alpha}(\partial\Omega)  \,.
\end{equation}
Then by combining inequalities (\ref{qrs1}) and (\ref{qrs2}), we deduce that $Q[\frac{\partial S_{ {\mathbf{a}} }}{\partial x_{r}}\circ\Theta ,\cdot,\cdot]$ is continuous from $C^{0,\alpha}(\partial\Omega)\times C^{0,\beta}(\partial\Omega)$
to $(C^{0,\alpha}(\partial\Omega),|\cdot:\partial\Omega|_{\alpha})$ and thus the proof is complete.   \hfill  $\Box$ 

\vspace{\baselineskip}

In the next lemma, we introduce a formula for the tangential derivatives of 
$Q[\frac{\partial S_{ {\mathbf{a}} }}{\partial x_{r}}\circ\Theta,g,\mu]$. 
\begin{lem}
\label{formula}
Let ${\mathbf{a}}$ be as in (\ref{introd0}), (\ref{ellip}), (\ref{symr}). Let $S_{ {\mathbf{a}} }$ be a fundamental solution of $P[{\mathbf{a}},D]$. 
 Let $\alpha\in]0,1[$, $\theta\in]0,1]$. 
Let $\Omega$ be a bounded open subset of ${\mathbb{R}}^{n}$ of class $C^{2,\alpha}$. Let $r\in\{1,\dots,n\}$. Let $g\in C^{1,\theta}(\partial\Omega)$, $\mu\in C^{1}(\partial\Omega)$. Then $Q[\frac{\partial S_{ {\mathbf{a}} }}{\partial x_{r}}\circ\Theta,g,\mu]\in C^{1}(\partial\Omega)$ and 
the following formula holds.
\begin{eqnarray}
\label{formula1}
\lefteqn{
M_{lj}\left[Q\left[\frac{\partial S_{ {\mathbf{a}} }}{\partial x_{r}}\circ\Theta,g,\mu\right]\right]
}
\\ \nonumber
&& 
= \nu_l(x)Q\left[\frac{\partial S_{ {\mathbf{a}} }}{\partial x_{r}}\circ\Theta,D_{ {\mathbf{a}} ,j}g,\mu\right](x) - \nu_j(x)Q\left[\frac{\partial S_{ {\mathbf{a}} }}{\partial x_{r}}\circ\Theta,D_{ {\mathbf{a}} ,l}g,\mu\right](x) 
\\ \nonumber
&& 
 + \nu_l(x)Q\left[\frac{\partial S_{ {\mathbf{a}} }}{\partial x_{r}}\circ\Theta, g , 
 \sum_{s=1}^{n} M_{sj}[\sum_{h=1}^{n} \frac{ a_{sh}\nu_h}{
 \nu^{t}a^{(2)}\nu  }\mu]\right](x)
 \\ \nonumber
&& \qquad
- 
\nu_j(x)Q\left[\frac{\partial S_{ {\mathbf{a}} }}{\partial x_{r}}\circ\Theta, g , \sum_{s =1}^{n}M_{sl}[\sum_{h=1}^{n} \frac{a_{sh}\nu_h}{
  \nu^{t}a^{(2)}\nu} \mu ]\right](x)
\\ \nonumber
&&   
  + \sum_{s,h=1}^{n} a_{sh} \nu_l(x) \biggl\{\biggr.
  Q\left[\frac{\partial S_{ {\mathbf{a}} }}{\partial x_{s}}\circ\Theta,\nu_j,\frac{M_{hr}[g]\mu}{ \nu^{t}a^{(2)}\nu}\right](x) 
  \\ \nonumber
&& \qquad
+ 
Q\left[\frac{\partial S_{ {\mathbf{a}} }}{\partial x_{s}}\circ\Theta, g, M_{hr}[\frac{\nu_j\mu}
{\nu^{t}a^{(2)}\nu} ]\right](x) 
\biggl.\biggr\}
\\ \nonumber
&&   
- \sum_{s,h=1}^{n} a_{sh} \nu_j(x) \biggl\{\biggr.
Q\left[\frac{\partial S_{ {\mathbf{a}} }}{\partial x_{s}}\circ\Theta,\nu_l,\frac{M_{hr}[g]\mu}{
\nu^{t}a^{(2)}\nu
}\right](x) 
\\ \nonumber
&& \qquad
+ Q\left[\frac{\partial S_{ {\mathbf{a}} }}{\partial x_{t}}\circ\Theta, g, M_{hr}[\frac{\nu_l\mu}{ \nu^{t}a^{(2)}\nu} ]\right](x) \biggl.\biggr\} 
\\ \nonumber
&& 
 -\sum_{t=1}^{n} a_{s}\biggl\{\biggr.\nu_l(x)
 Q\left[\frac{\partial S_{ {\mathbf{a}} }}{\partial x_{s}}\circ\Theta,g,\frac{\nu_j\nu_r}{\nu^{t}a^{(2)}\nu}\mu\right](x) 
 \\ \nonumber
&& \qquad
 -\nu_j(x)
 Q\left[\frac{\partial S_{ {\mathbf{a}} }}{\partial x_{s}}\circ\Theta,g,\frac{\nu_l\nu_r}{
\nu^{t}a^{(2)}\nu
}\mu\right](x) \biggl.\biggr\} 
\\ \nonumber
&& 
-a\left\{ 
g(x)\left[
\nu_l(x)v[\partial\Omega  ,S_{ {\mathbf{a}} }  ,  
\frac{\nu_j\nu_r}{
\nu^{t}a^{(2)}\nu}\mu](x)
-
\nu_j(x)v[\partial\Omega      ,S_{ {\mathbf{a}} } ,  
\frac{\nu_l\nu_r}{
\nu^{t}a^{(2)}\nu}\mu](x)
\right]\right.
\\ \nonumber
&& 
-
\left.
\left[
\nu_l(x)v[\partial\Omega    ,S_{ {\mathbf{a}} } ,  
g\frac{\nu_j\nu_r}{
\nu^{t}a^{(2)}\nu}\mu](x)
-
\nu_j(x)v[\partial\Omega    ,S_{ {\mathbf{a}} }   ,  
g\frac{\nu_l\nu_r}{
\nu^{t}a^{(2)}\nu}\mu](x)
\right]
\right\} 
\end{eqnarray}
for all $x\in\partial\Omega$ and  $l,j\in\{1,\dots,n\}$. (For  $Q$ see (\ref{qrs0}).) 
\end{lem}
{\bf Proof.}  Let $R\in]0,+\infty[$ be such that ${\mathrm{cl}}\Omega\subseteq{\mathbb{B}}_{n}(0,R)$. Let  ` $\tilde{\ }$ '  be an extension operator as in Lemma \ref{extsch}, defined either on $C^{1,\theta}(\partial\Omega)$ or on $C^{1,\alpha}(\partial\Omega)$ depending on whether it has been  applied to $g\in C^{1,\theta}(\partial\Omega)$ or to $\nu_{l}\in C^{1,\alpha}(\partial\Omega)$ for $l=1,\dots,n$. \par

Now we fix $\beta\in]0,\min\{\theta,\alpha\}[$ and we first prove the formula under the assumption that $\mu\in C^{1,\beta}(\partial\Omega)$. By Theorem 
\ref{qrssm} (ii), we already know that $Q^{\sharp}[\frac{\partial S_{ {\mathbf{a}} }}{\partial x_{r}}\circ\Theta,g,\mu] $ belongs to $  C^{1}({\mathrm{cl}}\Omega)$. Then we find convenient to introduce the notation
\[
M^{\sharp}_{lj}[f](x)\equiv \tilde{\nu}_{l}(x)\frac{\partial f}{\partial x_{j}}(x)
-
\tilde{\nu}_{j}(x)\frac{\partial f}{\partial x_{l}}(x)\qquad\forall x\in {\mathrm{cl}}\Omega\,,
\]
for all $f\in C^{1}({\mathrm{cl}}\Omega)$. If necessary, we write
$M^{\sharp}_{lj,x}$ to emphasize that we are taking $x$ as variable of the differential operator $M^{\sharp}_{lj}$. Next we fix $x\in \Omega$ and we compute
\[
\tilde{\nu}_{l}(x)\frac{\partial}{\partial x_{j}}Q^{\sharp} [\frac{\partial S_{ {\mathbf{a}} }}{\partial x_{r}}\circ\Theta,\tilde{g},\mu](x)
-
\tilde{\nu}_{j}(x)\frac{\partial}{\partial x_{l}}Q^{\sharp}[\frac{\partial S_{ {\mathbf{a}} }}{\partial x_{r}}\circ\Theta,\tilde{g},\mu](x)\,.
\]
Clearly,
\begin{eqnarray*}
\lefteqn{
\frac{\partial}{\partial x_{l}}Q^{\sharp}[\frac{\partial S_{ {\mathbf{a}} }}{\partial x_{r}}\circ\Theta,\tilde{g},\mu](x)
=\int_{\partial\Omega}
\frac{\partial\tilde{g}}{\partial x_{l}}(x)\frac{\partial}{\partial x_{r}}S_{{\mathbf{a}}}(x-y)\mu(y)\,d\sigma_{y}
}
\\
&&
\qquad\qquad\qquad\quad
+
\int_{\partial\Omega} (\tilde{g}(x)-\tilde{g}(y))\frac{\partial^{2}}{\partial x_{l}\partial x_{r}}S_{{\mathbf{a}}}(x-y)
\mu(y)\,d\sigma_{y}\,.
\end{eqnarray*}
 Now to shorten our notation, we set 
 \[
 J_{1}(x)\equiv \int_{\partial\Omega}
\frac{\partial\tilde{g}}{\partial x_{l}}(x)\frac{\partial}{\partial x_{r}}S_{{\mathbf{a}}}(x-y)\mu(y)\,d\sigma_{y}\,.
\]
Then we have 
\begin{eqnarray*}
\lefteqn{
\frac{\partial}{\partial x_{l}}Q^{\sharp}[\frac{\partial S_{ {\mathbf{a}} }}{\partial x_{r}}\circ\Theta,\tilde{g},\mu](x) 
 }
\\  \nonumber
&&
\qquad
= J_1(x) - \int_{\partial\Omega} (\tilde{g}(x)-\tilde{g}(y)) 
\\  \nonumber
&&
\qquad\quad
\times\sum_{s,h=1}^{n}\frac{ \nu_s(y)a_{sh}\nu_h(y)}{    \nu^{t}(y)a^{(2)}\nu (y)  }
\frac{\partial}{\partial y_{l}   }
[\frac{\partial}{\partial x_{r}}S_{   {\mathbf{a}}   }
(x-y)        ]\mu(y)d\sigma_y 
\\  \nonumber
&&
\qquad 
=  J_1(x) - \int_{\partial\Omega} (\tilde{g}(x)-\tilde{g}(y)) 
\\  \nonumber
&&
\qquad\quad
\times  \sum_{s=1}^{n} \left(\nu_{s}(y)
\frac{\partial}{\partial y_{l}}
-
\nu_{l}(y)
\frac{\partial}{\partial y_{s}}
\right)\left[
\frac{\partial}{\partial x_{r}}S_{ {\mathbf{a}} }(x-y)
\right]
\\  \nonumber
&&
\qquad\quad
\times
 \sum_{h=1}^{n} \frac{a_{sh}\nu_h(y)}{     \nu^{t}(y)a^{(2)}\nu (y)    }\mu(y)d\sigma_y
 \\  \nonumber
&&
\qquad\quad
  - \int_{\partial\Omega} (\tilde{g}(x)-\tilde{g}(y))
  \\  \nonumber
&&
\qquad\quad
\times \sum_{s=1}^{n}
 \frac{\partial}{\partial y_{s}}\left[
 \frac{\partial}{\partial x_{r}}S_{ {\mathbf{a}} }(x-y)
 \right]
  \sum_{h=1}^{n} a_{sh}\nu_h(y) \frac{\nu_l(y)}{  \nu^{t}(y)a^{(2)}\nu (y)  }\mu(y)d\sigma_y\,.
 \end{eqnarray*}
By Lemma \ref{gagre}, the second term  in the right hand side  takes the following form  
\begin{eqnarray*}
\lefteqn{
 \int_{\partial\Omega} (\tilde{g}(x)-\tilde{g}(y)) \sum_{s=1}^{n} M_{sl,y}
\left[
 \frac{\partial}{\partial x_{r}}S_{ {\mathbf{a}} }(x-y)
\right]
\frac{  (a^{(2)}\nu(y))_{s}  }{   \nu^{t}(y)a^{(2)}\nu (y)   }\mu(y)\,d\sigma_y
}
\\ \nonumber
&&\qquad
 = - \int_{\partial\Omega} \sum_{s=1}^{n} M_{sl,y} [\tilde{g}(x)-\tilde{g}(y)]
  \frac{\partial}{\partial x_{r}}S_{ {\mathbf{a}} }(x-y)
  \frac{  (a^{(2)}\nu(y))_{s}  }{   \nu^{t}(y)a^{(2)}\nu (y)   }\mu(y)\,d\sigma_y
\\ \nonumber
&&\qquad
   - \int_{\partial\Omega} \sum_{s=1}^{n} (\tilde{g}(x)-\tilde{g}(y)) 
 \frac{\partial}{\partial x_{r}}S_{ {\mathbf{a}} }(x-y)  M_{sl} \left[
  \frac{  (a^{(2)}\nu )_{s}  }{   \nu^{t} a^{(2)}\nu    }\mu 
 \right](y)\,d\sigma_y\,.
\end{eqnarray*}   
Since $M_{sl,y}[  \tilde{g}(x)-\tilde{g}(y)  ]=-M_{sl}[\tilde{g}](y)$, we have 
\begin{eqnarray*}
\lefteqn{
\frac{\partial}{\partial x_{l}}Q^{\sharp}[\frac{\partial S_{ {\mathbf{a}} }}{\partial x_{r}}\circ\Theta,\tilde{g},\mu](x)
} 
\\
&&
=
\frac{\partial\tilde{g}}{\partial x_{l}}(x)\int_{\partial\Omega}
\frac{\partial}{\partial x_{r}}S_{ {\mathbf{a}} }(x-y)\mu(y)\,d\sigma_{y}
\\
&&
\quad
- \int_{\partial\Omega} \sum_{s=1}^{n} M_{sl} [\tilde{g}](y)
\frac{\partial}{\partial x_{r}}S_{ {\mathbf{a}} }(x-y)
 \frac{  (a^{(2)}\nu(y))_{s}  }{   \nu^{t}(y)a^{(2)}\nu (y)   }\mu(y)\,d\sigma_y
 \\
&&
\quad
 + \int_{\partial\Omega} \sum_{s=1}^{n} (\tilde{g}(x)-\tilde{g}(y))
 \frac{\partial}{\partial x_{r}}S_{ {\mathbf{a}} }(x-y)
 M_{sl} \left[
 \frac{  (a^{(2)}\nu )_{s}  }{   \nu^{t} a^{(2)}\nu     } \mu
 \right](y)\,d\sigma_y
 \\
&&
\quad
 - \int_{\partial\Omega} (\tilde{g}(x)-\tilde{g}(y))
 \\  \nonumber
&&
\quad
\times \sum_{s=1}^{n}
 \frac{\partial}{\partial y_{s}} \left[
  \frac{\partial}{\partial x_{r}}S_{ {\mathbf{a}} }(x-y)
 \right](a^{(2)}\nu(y))_{s} \frac{\nu_{l}(y)}{
 \nu^{t}(y)a^{(2)}\nu (y) 
 }\mu(y)\,d\sigma_y
\,.
\end{eqnarray*}
Accordingly, we have
\begin{eqnarray}
\label{formula2}
\lefteqn{M^{\sharp}_{lj}
\left[Q^{\sharp}[\frac{\partial S_{ {\mathbf{a}} }}{\partial x_{r}}\circ\Theta,\tilde{g},\mu]\right](x)
}
\\ \nonumber
&& 
=M^{\sharp}_{lj} [      \tilde{g}]        (x)
\int_{\partial\Omega} \frac{\partial}{\partial x_{r}}S_{ {\mathbf{a}} }(x-y)
\mu(y)\,d\sigma_y
\\ \nonumber
&&\quad
- \int_{\partial\Omega} \sum_{s=1}^{n} \{\tilde{\nu}_l(x) M_{sj}[\tilde{g}](y)
-
\tilde{\nu}_j(x) M_{sl}[\tilde{g}](y)\}\frac{\partial}{\partial x_{r}}S_{ {\mathbf{a}} }(x-y)
\\  \nonumber
&&
\quad
\times
\frac{  (a^{(2)}\nu(y))_{s}  }{   \nu^{t}(y)a^{(2)}\nu (y)   }\mu(y)\,d\sigma_{y}
\\ \nonumber
&&\quad 
+ \int_{\partial\Omega} \sum_{s=1}^{n} (\tilde{g}(x)-\tilde{g}(y)) \frac{\partial}{\partial x_{r}}S_{ {\mathbf{a}} }(x-y)
\\ \nonumber
&&\quad 
\times
\left\{
\tilde{\nu}_{l}(x)M_{sj}\left[  
\frac{  (a^{(2)}\nu )_{s}  }{   \nu^{t} a^{(2)}\nu    }\mu
 \right](y)
-
\tilde{\nu}_{j}(x)M_{sl}\left[  
\frac{  (a^{(2)}\nu )_{s}  }{   \nu^{t} a^{(2)}\nu    }\mu
 \right](y)
\right\}\,d\sigma_{y}
\\ \nonumber
&&\quad 
- \int_{\partial\Omega} (\tilde{g}(x)-\tilde{g}(y)) \sum_{s=1}^{n}
\frac{\partial}{\partial y_{s}}
\left[\frac{\partial}{\partial x_{r}}S_{ {\mathbf{a}} }(x-y)\right]
\\ \nonumber
&&\quad 
\times
 (a^{(2)}\nu )_{s}(y)
 \frac{
 \tilde{\nu}_l(x)\nu_j(y)-\tilde{\nu}_j(x)\nu_l(y)
 }{\nu^{t}(y)a^{(2)}\nu (y) }\mu(y)\,d\sigma_{y}
\end{eqnarray}
We now consider the first two terms in the right hand side of formula 
(\ref{formula2}). By the obvious identity
\[
M_{lj}^{\sharp}[\tilde{g}]=\tilde{\nu}_{l}
\left[
\frac{\partial}{\partial x_{j}}\tilde{g}-\frac{
D\tilde{g}a^{(2)}\tilde{\nu}
}{
\tilde{\nu}^{t} a^{(2)}\tilde{\nu}
}\tilde{\nu}_{j}
\right]
-
\tilde{\nu}_{j}
\left[
\frac{\partial}{\partial x_{l}}\tilde{g}-\frac{
D\tilde{g}a^{(2)}\tilde{\nu}
}{
\tilde{\nu}^{t} a^{(2)}\tilde{\nu}
}\tilde{\nu}_{l}
\right]\qquad{\mathrm{in}}\ {\mathrm{cl}}\Omega\,,
\]
and by the corresponding formula for $M_{lj} [\tilde{g}]$  on $\partial\Omega$, 
and by fomula (\ref{dam}) and by straightforward computations, we obtain
\begin{eqnarray}
\label{formula3}
\lefteqn{
M^{\sharp}_{lj}[\tilde{g}](x)
\int_{\partial\Omega} \frac{\partial}{\partial x_{r}}S_{{\mathbf{a}}}(x-y)\mu(y)\,d\sigma_{y}
}
\\ \nonumber
&&\ 
-\int_{\partial\Omega}\sum_{s=1}^{n}
\left\{
\tilde{\nu}_{l}(x)M_{sj}[\tilde{g}](y)
-
\tilde{\nu}_{j}(x)M_{sl}[\tilde{g}](y)
\right\}
\\ \nonumber
&&\ 
\times
\frac{\partial}{\partial x_{r}}S_{{\mathbf{a}}}(x-y)
\frac{  (a^{(2)}\nu(y) )_{s}  }{   \nu^{t}(y) a^{(2)}\nu  (y)  }
\mu(y)\,d\sigma_{y}
\\ \nonumber
&&
=\tilde{\nu}_{l}(x)
\left[
\frac{\partial}{\partial x_{j}}\tilde{g}(x)
-\frac{D\tilde{g}(x)a^{(2)}\tilde{\nu}(x)}{  \tilde{\nu}^{t}(x) a^{(2)}\tilde{\nu}  (x)   }
\tilde{\nu}_{j}(x)
\right]
\int_{\partial\Omega}\frac{\partial}{\partial x_{r}}S_{{\mathbf{a}}}(x-y)
\mu(y)\,d\sigma_{y}
\\ \nonumber
&&\
-\tilde{\nu}_{j}(x)
\left[
\frac{\partial}{\partial x_{l}}\tilde{g}(x)
-\frac{D\tilde{g}(x)a^{(2)}\tilde{\nu}(x)}{   \tilde{\nu}^{t}(x) a^{(2)}\tilde{\nu}  (x)   }
\tilde{\nu}_{l}(x)
\right]
\int_{\partial\Omega}\frac{\partial}{\partial x_{r}}S_{{\mathbf{a}}}(x-y)
\mu(y)\,d\sigma_{y}
\\ \nonumber
&&\
-\tilde{\nu}_{l}(x)\int_{\partial\Omega}
\left[
\frac{\partial}{\partial y_{j}}\tilde{g}(y)
-\frac{D\tilde{g}(y)a^{(2)}\tilde{\nu}(y)}{   \tilde{\nu}^{t}(y) a^{(2)}\tilde{\nu}  (y)   }
\tilde{\nu}_{j}(y)
\right]
\\  \nonumber
&&\
\times
\left(
\sum_{s, h=1}^{n}\tilde{\nu}_{s}(y)\frac{a^{(2)}_{sh}\nu_{h}(y)}{
\tilde{\nu}^{t}(y) a^{(2)}\tilde{\nu}  (y) 
}
\right)\frac{\partial}{\partial x_{r}}S_{{\mathbf{a}}}(x-y)
\mu(y)\,d\sigma_{y}
\\ \nonumber
&&\
+\tilde{\nu}_{l}(x)\int_{\partial\Omega}\tilde{\nu}_{j}(y)
\biggl\{\biggr.
\sum_{s, h=1}^{n}\frac{\partial}{\partial y_{s}}\tilde{g}(y)
\frac{a_{sh}\nu_{h}(y)}{   \tilde{\nu}^{t}(y) a^{(2)}\tilde{\nu}  (y)  }
-
\frac{D\tilde{g}(y)a^{(2)}\tilde{\nu}(y)}{   \tilde{\nu}^{t}(y) a^{(2)}\tilde{\nu}  (y)   }
\\ \nonumber
&&\
\times\left(
 \tilde{\nu}_{s}(y)\frac{a_{sh}\nu_{h}(y)}{
 \tilde{\nu}^{t}(y) a^{(2)}\tilde{\nu}  (y) 
}
\right)
\biggl.\biggr\} 
 \frac{\partial}{\partial x_{r}}S_{{\mathbf{a}}}(x-y)
\mu(y)\,d\sigma_{y}
\\ \nonumber
&&\
+\tilde{\nu}_{j}(x)\int_{\partial\Omega}
\left[
\frac{\partial}{\partial y_{l}}\tilde{g}(y)
-\frac{D\tilde{g}(y)a^{(2)}\tilde{\nu}(y)}{   \tilde{\nu}^{t}(y) a^{(2)}\tilde{\nu} (y)   }
\tilde{\nu}_{l}(y)
\right]
\\ \nonumber
&&\
\times
\left(
\sum_{s, h=1}^{n}\tilde{\nu}_{s}(y)\frac{a_{sh}\nu_{h}(y)}{
 \tilde{\nu}^{t}(y) a^{(2)}\tilde{\nu}  (y) 
}
\right)\frac{\partial}{\partial x_{r}}S_{{\mathbf{a}}}(x-y)
\mu(y)\,d\sigma_{y}
\\ \nonumber
&&\
-\tilde{\nu}_{j}(x)\int_{\partial\Omega}\tilde{\nu}_{l}(y)
\biggl\{\biggr.
\sum_{s, h=1}^{n}\frac{\partial}{\partial y_{s}}\tilde{g}(y)
\frac{a_{sh}\nu_{h}(y)}{  \tilde{\nu}^{t}(y) a^{(2)}\tilde{\nu}  (y)  }
-
\frac{D\tilde{g}(y)a^{(2)}\tilde{\nu}(y)}{  \tilde{\nu}^{t}(y) a^{(2)}\tilde{\nu}  (y)   }
\\ \nonumber
&&\
\times
\left(
  \tilde{\nu}_{s}(y)\frac{a_{sh}\nu_{h}(y)}{
 \tilde{\nu}^{t}(y) a^{(2)}\tilde{\nu}  (y) 
}
\right)
\biggl.\biggr\} \frac{\partial}{\partial x_{r}}S_{{\mathbf{a}}}(x-y)
\mu(y)\,d\sigma_{y}\,.
\end{eqnarray}
Since
\[
\tilde{\nu}(y)=\nu(y)\,,\qquad
\left(
\sum_{s, h=1}^{n}\tilde{\nu}_{s}(y)\frac{a_{sh}\nu_{h}(y)}{
 \tilde{\nu}^{t}(y) a^{(2)}\tilde{\nu}  (y) 
}
\right)=1\qquad\forall y\in\partial\Omega\,,
\]
we have
\[
\left\{
\sum_{s, h=1}^{n}\frac{\partial}{\partial y_{s}}\tilde{g}(y)
\frac{a_{sh}\nu_{h}(y)}{   \tilde{\nu}^{t}(y) a^{(2)} \tilde{\nu}  (y)  }
-
\frac{D\tilde{g}(y)a^{(2)}\tilde{\nu}(y)}{   \tilde{\nu}^{t}(y) a^{(2)} \tilde{\nu} (y)   }
\left(
 \tilde{\nu}_{s}(y)\frac{a_{sh}\nu_{h}(y)}{
  \tilde{\nu}^{t}(y) a^{(2)} \tilde{\nu} (y) 
}
\right)
\right\}=0\,,
\]
for all $y\in\partial\Omega$ 
and accordingly, the right hand side of (\ref{formula3}) equals
\begin{eqnarray*}
\lefteqn{
\tilde{\nu}_{l}(x)Q^{\sharp}\left[\frac{\partial S_{ {\mathbf{a}} }}{\partial x_{r}}\circ\Theta,
\frac{\partial}{\partial x_{j}}\tilde{g} 
-\frac{D\tilde{g} a^{(2)}\tilde{\nu} }{   \nu^{t}  a^{(2)}\nu    }
\tilde{\nu}_{j},\mu
\right](x)
}
\\
&&\qquad\qquad
-
\tilde{\nu}_{j}(x)Q^{\sharp}\left[\frac{\partial S_{ {\mathbf{a}} }}{\partial x_{r}}\circ\Theta,
\frac{\partial}{\partial x_{l}}\tilde{g} 
-\frac{D\tilde{g} a^{(2)}\tilde{\nu} }{   \nu^{t}  a^{(2)}\nu    }
\tilde{\nu}_{l},\mu
\right](x)\,.
\end{eqnarray*}
Next we consider the third term in the right hand side of  formula (\ref{formula2}), and we note that
\begin{eqnarray}
\label{formula4}
\lefteqn{
\int_{\partial\Omega}\sum_{s=1}^{n}(\tilde{g}(x)-\tilde{g}(y))
\frac{\partial}{\partial x_{r}}S_{ {\mathbf{a}} }(x-y)
}
\\ \nonumber
&&\quad
\times\left\{
\tilde{\nu}_{l}(x)M_{sj}\left[
\frac{
(a^{(2)}\nu)_{s}
}{
 \nu^{t}  a^{(2)}\nu 
}\mu
\right](y)
-
\tilde{\nu}_{j}(x)M_{sl}\left[
\frac{
(a^{(2)}\nu)_{s}
}{
 \nu^{t}  a^{(2)}\nu 
}\mu
\right](y)
\right\}\,d\sigma_{y}
\\ \nonumber
&&
=
\tilde{\nu}_{l}(x)Q^{\sharp}\left[\frac{\partial S_{ {\mathbf{a}} }}{\partial x_{r}}\circ\Theta,\tilde{g},\sum_{s=1}^{n}M_{sj}\left[
\frac{
(a^{(2)}\nu)_{s}
}{
 \nu^{t}  a^{(2)}\nu 
}\mu
\right]
\right](x)
\\ \nonumber
&&\quad
-
\tilde{\nu}_{j}(x)Q^{\sharp}\left[\frac{\partial S_{ {\mathbf{a}} }}{\partial x_{r}}\circ\Theta,\tilde{g},\sum_{s=1}^{n}M_{sl}\left[
\frac{
(a^{(2)}\nu)_{s}
}{
 \nu^{t}  a^{(2)}\nu 
}\mu
\right]
\right](x)\,.
\end{eqnarray}
Next we consider the last integral in the right hand side of  formula (\ref{formula2}) and we note that if $x\in \Omega$ and $y\in\partial\Omega$, we have
\[
\sum_{s,h=1}^{n}\frac{\partial}{\partial x_{h}}
\left[
a_{sh} \frac{\partial}{\partial x_{s}}S_{ {\mathbf{a}} }(x-y) 
\right]
+
\sum_{s=1}^{n}a_{s}\frac{\partial}{\partial x_{s}}S_{ {\mathbf{a}} }(x-y)
+aS_{ {\mathbf{a}} }(x-y)=0\,.
\]
Thus   we obtain
\begin{eqnarray*}
\lefteqn{
\sum_{s,h=1}^{n}a_{sh}\nu_{h}(y)\frac{\partial}{\partial x_{r}}
\left[
\frac{\partial}{\partial y_{s}}S_{ {\mathbf{a}} }(x-y)
\right]
}
\\ \nonumber
&&\qquad
=\sum_{s,h=1}^{n}a_{sh}
\left(
\nu_{h}(y)\frac{\partial}{\partial y_{r}}
-
\nu_{r}(y)\frac{\partial}{\partial y_{h}}
\right)\left[
\frac{\partial}{\partial x_{s}}S_{ {\mathbf{a}} }(x-y)
\right]
\\ \nonumber
&&\qquad\quad
+\nu_{r}(y)\sum_{s=1}^{n}a_{s}
\frac{\partial}{\partial x_{s}}S_{ {\mathbf{a}} }(x-y)
+\nu_{r}(y)aS_{ {\mathbf{a}} }(x-y)\,,
\end{eqnarray*}
and we note that the first parenthesis  in the right hand side 
equals $M_{hr,y}$.
 Then the last integral in the right hand side of formula  (\ref{formula2}) equals 
\begin{eqnarray}
\label{formula5}
\lefteqn{
\int_{\partial\Omega}(\tilde{g}(x)-\tilde{g}(y))
\sum_{s, h=1}^{n}a_{sh}\nu_{h}(y)\frac{\partial}{\partial y_{s}}
\left[ 
\frac{\partial}{\partial x_{r}}S_{ {\mathbf{a}} }(x-y)
\right]
}
\\ \nonumber
&&\qquad\quad
\times\frac{\tilde{\nu}_{l}(x)\nu_{j}(y)-\tilde{\nu}_{j}(x)\nu_{l}(y)}{
 \nu^{t}(y) a^{(2)}\nu  (y)
 }\mu (y)\,d\sigma_{y}
 \\ \nonumber
&&\qquad
=\int_{\partial\Omega}(\tilde{g}(x)-\tilde{g}(y))
\biggl\{\biggr.\sum_{s, h=1}^{n}a_{sh}
M_{hr,y}\left[ 
\frac{\partial}{\partial x_{s}}S_{ {\mathbf{a}} }(x-y)
\right]
 \\ \nonumber
&&\qquad\quad
+\nu_{r}(y)\sum_{s=1}^{n}a_{s}\frac{\partial}{\partial x_{s}}S_{ {\mathbf{a}} }(x-y)
+\nu_{r}(y)aS_{ {\mathbf{a}} }(x-y)\biggl.\biggr\}
\\ \nonumber
&&\qquad\quad
\times \frac{\tilde{\nu}_{l}(x)\nu_{j}(y)-\tilde{\nu}_{j}(x)\nu_{l}(y)}{
 \nu^{t}(y) a^{(2)}\nu  (y)
 }\mu (y)\,d\sigma_{y}
 \\ \nonumber
&&\qquad
 =\sum_{s, h=1}^{n}a_{sh}
\int_{\partial\Omega}(\tilde{g}(x)-\tilde{g}(y))
 M_{hr,y}\left[
 \frac{\partial}{\partial x_{s}}S_{ {\mathbf{a}} }(x-y)
 \right]
  \\ \nonumber
&&\qquad\quad
 \times \frac{\tilde{\nu}_{l}(x)(\tilde{\nu}_{j}(y) -\tilde{\nu}_{j}(x))
 +\tilde{\nu}_{j}(x)(\tilde{\nu}_{l}(x)-\tilde{\nu}_{l}(y))
 }{
 \nu^{t}(y) a^{(2)}\nu  (y)
 }\mu (y)\,d\sigma_{y}
  \\ \nonumber
&&\qquad\quad
+\int_{\partial\Omega}(\tilde{g}(x)-\tilde{g}(y))
\left[
\sum_{s=1}^{n}a_{s}\frac{\partial}{\partial x_{s}}S_{ {\mathbf{a}} }(x-y)
+aS_{ {\mathbf{a}} }(x-y)
\right]
\\ \nonumber
&&\qquad\quad
\times
\frac{\tilde{\nu}_{l}(x)\nu_{j}(y)-\tilde{\nu}_{j}(x)\nu_{l}(y)}{
 \nu^{t}(y) a^{(2)}\nu  (y)
 }\nu_{r}(y)\mu(y)\,d\sigma_{y}\,.
\end{eqnarray}
 We now consider separately each of the terms in the right hand side of formula (\ref{formula5}). By Lemma \ref{gagre} and by equality
 $-M_{hr,y}[\tilde{g}(x)-\tilde{g}(y)]=M_{hr,y}[\tilde{g}(y)]$, the first integral in the right hand side of formula (\ref{formula5}) equals
\begin{eqnarray}
\label{formula6}
\lefteqn{
\int_{\partial\Omega}(\tilde{g}(x)-\tilde{g}(y))
 M_{hr,y}\left[
 \frac{\partial}{\partial x_{s}}S_{ {\mathbf{a}} }(x-y)
 \right]
 }
  \\ \nonumber
&& \quad
 \times \frac{\tilde{\nu}_{l}(x)(\tilde{\nu}_{j}(y) -\tilde{\nu}_{j}(x))
 +\tilde{\nu}_{j}(x)(\tilde{\nu}_{l}(x)-\tilde{\nu}_{l}(y))
 }{
 \nu^{t}(y) a^{(2)}\nu  (y)
 }\mu (y)\,d\sigma_{y}
\\ \nonumber
&& 
=\int_{\partial\Omega}M_{hr}[\tilde{g}] \frac{\partial}{\partial x_{s}}S_{ {\mathbf{a}} }(x-y)
 \\ \nonumber
&& \quad
 \times
 \left(
 -\tilde{\nu}_{l}(x)\frac{\tilde{\nu}_{j}(x)-  \nu_{j}(y)}{
 \nu^{t}(y) a^{(2)}\nu  (y)
 }
 +
 \tilde{\nu}_{j}(x)\frac{\tilde{\nu}_{l}(x)-  \nu_{l}(y)}{
 \nu^{t}(y) a^{(2)}\nu  (y)
 }
 \right)\mu (y)\,d\sigma_{y}
\\ \nonumber
&& \quad
+\int_{\partial\Omega}(\tilde{g}(x)-\tilde{g}(y))
 \frac{\partial}{\partial x_{s}}S_{ {\mathbf{a}} }(x-y)
\\ \nonumber
&& \quad
 \times 
\biggl(
 -\tilde{\nu}_{l}(x)M_{hr}\left[
 \frac{ \nu_{j}\mu     }{    \nu^{t}  a^{(2)}\nu    }
 \right](y)
\biggr.
\biggl.
+\tilde{\nu}_{j}(x)M_{hr}\left[
 \frac{ \nu_{l}\mu     }{    \nu^{t}  a^{(2)}\nu    }
 \right](y)
\biggr)\,d\sigma_{y}
\\ \nonumber
&& 
=-\tilde{\nu}_{l}(x)\int_{\partial\Omega}( \tilde{\nu}_{j}(x)-\nu_{j}(y))
\frac{\partial}{\partial x_{s}}S_{ {\mathbf{a}} }(x-y)
\frac{M_{hr}[\tilde{g}]}{     \nu^{t}(y) a^{(2)}\nu  (y)  }
\mu (y)\,d\sigma_{y}
\\ \nonumber
&& \quad
+\tilde{\nu}_{j}(x)\int_{\partial\Omega}( \tilde{\nu}_{l}(x)-\nu_{l}(y))
\frac{\partial}{\partial x_{s}}S_{ {\mathbf{a}} }(x-y)
\frac{M_{hr}[\tilde{g}]}{     \nu^{t}(y) a^{(2)}\nu  (y)  }
\mu (y)\,d\sigma_{y}
\\ \nonumber
&& \quad
-\tilde{\nu}_{l}(x)\int_{\partial\Omega}(\tilde{g}(x)-\tilde{g}(y))
\frac{\partial}{\partial x_{s}}S_{ {\mathbf{a}} }(x-y)M_{hr}\left[
 \frac{ \nu_{j}\mu     }{    \nu^{t}  a^{(2)}\nu    }
\right](y)\,d\sigma_{y}
\\ \nonumber
&& \quad
+\tilde{\nu}_{j}(x)\int_{\partial\Omega}(\tilde{g}(x)-\tilde{g}(y))
\frac{\partial}{\partial x_{s}}S_{ {\mathbf{a}} }(x-y)M_{hr}\left[
 \frac{ \nu_{l}\mu     }{    \nu^{t}  a^{(2)}\nu    }
\right](y)\,d\sigma_{y}
\\ \nonumber
&& 
=-\tilde{\nu}_{l}(x)
\biggl\{\biggr. Q^{\sharp}\left[\frac{\partial S_{ {\mathbf{a}} }}{\partial x_{s}}\circ\Theta,  \tilde{\nu}_{j}, 
\frac{M_{hr}[g]\mu}{ \nu^{t}  a^{(2)}\nu }\right](x)
\\ \nonumber
&& \qquad
+
 Q^{\sharp}\left[\frac{\partial S_{ {\mathbf{a}} }}{\partial x_{s}}\circ\Theta, \tilde{g}, M_{hr}\left[
\frac{\nu_{j}\mu}{ \nu^{t}  a^{(2)}\nu }\right]\right](x)
\biggl.\biggr\}
\\ \nonumber
&& \quad
+\tilde{\nu}_{j}(x)
\biggl\{\biggr. Q^{\sharp}\left[\frac{\partial S_{ {\mathbf{a}} }}{\partial x_{s}}\circ\Theta, \tilde{\nu}_{l}, 
\frac{M_{hr}[g]\mu}{ \nu^{t}  a^{(2)}\nu }\right](x)
\\ \nonumber
&& \qquad
+
 Q^{\sharp}\left[\frac{\partial S_{ {\mathbf{a}} }}{\partial x_{s}}\circ\Theta, \tilde{g}, M_{hr}\left[
\frac{\nu_{l}\mu}{ \nu^{t}  a^{(2)}\nu }\right]\right](x)
\biggl.\biggr\}\,.
\end{eqnarray}
 Next we note that the second integral   in the right hand side of formula (\ref{formula5}) equals
 \begin{eqnarray*}
\lefteqn{
\sum_{s=1}^{n}a_{s}
\biggl\{\biggr.
\tilde{\nu}_{l}(x)
Q^{\sharp}\left[\frac{\partial S_{ {\mathbf{a}} }}{\partial x_{s}}\circ\Theta,\tilde{g},\frac{\nu_{j}\nu_{r}}{\nu^{t}  a^{(2)}\nu }\mu\right](x)}
\\ \nonumber
&& \qquad
-
\tilde{\nu}_{j}(x)
Q^{\sharp}\left[\frac{\partial S_{ {\mathbf{a}} }}{\partial x_{s}}\circ\Theta,\tilde{g},\frac{\nu_{l}\nu_{r}}{\nu^{t}  a^{(2)}\nu }\mu\right](x)
\biggl.\biggr\}
\\ \nonumber
&&
+a\biggl\{\biggr.
\tilde{g}(x)
\biggl[\biggr.\tilde{\nu}_{l}(x)v\left[\partial\Omega   ,S_{ {\mathbf{a}} }   ,
\frac{\nu_{j}\nu_{r}}{\nu^{t}  a^{(2)}\nu }\mu
\right](x)
\\ \nonumber
&& \qquad
-
\tilde{\nu}_{j}(x)v\left[\partial\Omega   ,S_{ {\mathbf{a}} }   ,
\frac{\nu_{l}\nu_{r}}{\nu^{t}  a^{(2)}\nu }\mu
\right](x)
\biggl.\biggr]
\\ \nonumber
&&
-\left[
\tilde{\nu}_{l}(x)v\left[\partial\Omega    ,S_{ {\mathbf{a}} }   ,g
\frac{\nu_{j}\nu_{r}}{\nu^{t}  a^{(2)}\nu }\mu
\right](x)
-
\tilde{\nu}_{j}(x)v\left[\partial\Omega  ,S_{ {\mathbf{a}} }    ,g
\frac{\nu_{l}\nu_{r}}{\nu^{t}  a^{(2)}\nu }\mu
\right](x)
\right]\biggl.\biggr\}.
\end{eqnarray*}
Then by combining formulas (\ref{formula2})--(\ref{formula6}), we obtain the following formula
\begin{eqnarray}
\label{formula7}
\lefteqn{
M_{lj}^{\sharp}\left[Q^{\sharp}\left[\frac{\partial S_{ {\mathbf{a}} }}{\partial x_{r}}\circ\Theta,\tilde{g},\mu\right]\right](x)
}
\\ \nonumber
&& 
= 
\tilde{\nu}_{l}(x)Q^{\sharp}\left[\frac{\partial S_{ {\mathbf{a}} }}{\partial x_{r}}\circ\Theta,
\frac{\partial}{\partial x_{j}}\tilde{g} 
-\frac{D\tilde{g} a^{(2)}\tilde{\nu} }{   \nu^{t}  a^{(2)}\nu    }
\tilde{\nu}_{j}
,\mu\right](x)
\\ \nonumber
&& 
-
\tilde{\nu}_{j}(x)Q^{\sharp}\left[\frac{\partial S_{ {\mathbf{a}} }}{\partial x_{r}}\circ\Theta,
\frac{\partial}{\partial x_{l}}\tilde{g} 
-\frac{D\tilde{g} a^{(2)}\tilde{\nu} }{   \nu^{t}  a^{(2)}\nu    }
\tilde{\nu}_{l},\mu
\right](x)
\\ \nonumber
&& 
 + \tilde{\nu}_l(x)Q^{\sharp}\left[\frac{\partial S_{ {\mathbf{a}} }}{\partial x_{r}}\circ\Theta, \tilde{g} , 
 \sum_{s=1}^{n} M_{sj}[\sum_{h=1}^{n} \frac{ a_{sh}\nu_h}{
 \nu^{t}a^{(2)}\nu  }\mu]\right](x)
 \\ \nonumber
&& 
  - 
\tilde{\nu}_j(x)Q^{\sharp} \left[\frac{\partial S_{ {\mathbf{a}} }}{\partial x_{r}}\circ\Theta,\tilde{g} , \sum_{s =1}^{n}M_{sl}[\sum_{h=1}^{n} \frac{a_{sh}\nu_h}{
  \nu^{t}a^{(2)}\nu} \mu ]\right](x)
\\ \nonumber
&& 
  + \sum_{s,h=1}^{n} a_{sh} \tilde{\nu}_l(x) \biggl\{\biggr.
  Q^{\sharp}\left[\frac{\partial S_{ {\mathbf{a}} }}{\partial x_{s}}\circ\Theta,\nu_j,\frac{M_{hr}[g]\mu}{ \nu^{t}a^{(2)}\nu}\right](x) 
  \\ \nonumber
&& 
  + 
Q^{\sharp}\left[\frac{\partial S_{ {\mathbf{a}} }}{\partial x_{s}}\circ\Theta, \tilde{g}, M_{hr}[\frac{\nu_j\mu}
{\nu^{t}a^{(2)}\nu} ]\right](x) 
\biggl.\biggr\}
\\ \nonumber
&&  
- \sum_{s,h=1}^{n} a_{sh} \tilde{\nu}_j(x) \biggl\{\biggr.
Q^{\sharp}\left[\frac{\partial S_{ {\mathbf{a}} }}{\partial x_{s}}\circ\Theta,\nu_l,\frac{M_{hr}[g]\mu}{
\nu^{t}a^{(2)}\nu
}\right](x) 
\\ \nonumber
&& 
+ Q^{\sharp}\left[\frac{\partial S_{ {\mathbf{a}} }}{\partial x_{t}}\circ\Theta,\tilde{g}, M_{hr}[\frac{\nu_l\mu}{ \nu^{t}a^{(2)}\nu} ]\right](x) \biggl.\biggr\}
\\ \nonumber
&& 
 -\sum_{s=1}^{n} a_{s}\biggl\{\biggr.     \tilde{\nu}_l(x)
 Q^{\sharp}\left[\frac{\partial S_{ {\mathbf{a}} }}{\partial x_{s}}\circ\Theta,\tilde{g},\frac{\nu_j\nu_r}{\nu^{t}a^{(2)}\nu}\mu\right](x) 
 \\ \nonumber
&& 
 -\tilde{\nu}_j(x)
 Q^{\sharp}\left[\frac{\partial S_{ {\mathbf{a}} }}{\partial x_{s}}\circ\Theta,\tilde{g},\frac{\nu_l\nu_r}{
\nu^{t}a^{(2)}\nu
}\mu\right](x) \biggl.\biggr\}
\\ \nonumber
&& 
-a\left\{ 
g(x)\left[
\tilde{\nu}_l(x)v[\partial\Omega     ,S_{ {\mathbf{a}} } ,  
\frac{\nu_j\nu_r}{
\nu^{t}a^{(2)}\nu}\mu](x)
-
\tilde{\nu}_j(x)v[\partial\Omega   ,S_{ {\mathbf{a}} }  ,  
\frac{\nu_l\nu_r}{
\nu^{t}a^{(2)}\nu}\mu](x)
\right]\right.
\\ \nonumber
&& 
-
\left.
\left[
\tilde{\nu}_l(x)v[\partial\Omega    ,S_{ {\mathbf{a}} }  ,  
g\frac{\nu_j\nu_r}{
\nu^{t}a^{(2)}\nu}\mu](x)
-
\tilde{\nu}_j(x)v[\partial\Omega  ,S_{ {\mathbf{a}} }  ,  
g\frac{\nu_l\nu_r}{
\nu^{t}a^{(2)}\nu}\mu](x)
\right]
\right\}
\end{eqnarray}
Now under our assumptions, the first argument of the maps
  $Q^{\sharp}\left[\frac{\partial S_{ {\mathbf{a}} }}{\partial x_{r}}\circ\Theta,\cdot,\cdot\right]$ and $Q^{\sharp}\left[\frac{\partial S_{ {\mathbf{a}} }}{\partial x_{s}}\circ\Theta,\cdot,\cdot\right]$, which appear in the right hand side of formula (\ref{formula7}) belongs to the space $C^{0,\min\{\alpha,\theta\}}({\mathrm{cl}}\Omega)$ and the second argument of the maps  $Q^{\sharp}\left[\frac{\partial S_{ {\mathbf{a}} }}{\partial x_{r}}\circ\Theta,\cdot,\cdot\right]$, $Q^{\sharp}\left[\frac{\partial S_{ {\mathbf{a}} }}{\partial x_{s}}\circ\Theta,\cdot,\cdot\right] $, which appear in the right hand side of the formula (\ref{formula7}) belongs to $C^{0}(\partial\Omega)$. 
 By Theorem \ref{slay} (i) with $m=1$, the single layer potentials in the right hand side of formula (\ref{formula7}) are continuous in $x\in {\mathrm{cl}}\Omega$. Then Theorem \ref{qrssm} (i) implies that the right hand side of formula (\ref{formula7}) defines a continuous function of the variable $x\in {\mathrm{cl}}\Omega$. Since $\Omega$ is of class $C^{2,\alpha}$ and $\tilde{g}\in C^{1,\theta}({\mathrm{cl}}\Omega)$ and since we are assuming that $\mu\in C^{1,\beta}(\partial\Omega)$, Theorem  \ref{qrssm} (ii) implies that $M_{lj}^{\sharp}[Q^{\sharp}[\frac{\partial S_{ {\mathbf{a}} }}{\partial x_{r}}\circ\Theta,g,\mu]]$
 belongs to $C^{0}({\mathrm{cl}}\Omega)$. Hence, the equation of formula
 (\ref{formula7}) must hold for all $x\in {\mathrm{cl}}\Omega$, and in particular for all $x\in\partial\Omega$. Since     
 $Q^{\sharp}[\frac{\partial S_{ {\mathbf{a}} }}{\partial x_{r}}\circ\Theta,\cdot,\cdot]=Q[\frac{\partial S_{ {\mathbf{a}} }}{\partial x_{r}}\circ\Theta,\cdot,\cdot]$ and $M_{lj}^{\sharp}=M_{lj}$ 
 on $\partial\Omega$, we  conclude that (\ref{formula1}) holds.
\par

Next we assume that $\mu\in C^{1}(\partial\Omega)$. 
We denote by $P_{ljr}[g,\mu]$ the right hand side of (\ref{formula1}).  
 By Theorem \ref{qrs}  (i), the operators  $Q\left[\frac{\partial S_{ {\mathbf{a}} }}{\partial x_{r}}\circ\Theta,g,\cdot\right]$, 
 $Q\left[\frac{\partial S_{ {\mathbf{a}} }}{\partial x_{r}}\circ\Theta,D_{{\mathbf{a}},j}g,\cdot\right]$, $Q\left[\frac{\partial S_{ {\mathbf{a}} }}{\partial x_{r}}\circ\Theta,\nu_{l},\cdot\right]$ 
  are linear and continuous from the space $C^{0}(\partial\Omega)$ to $C^{0}(\partial\Omega)$. Then by Theorem \ref{v0a},
  and by the continuity of the pointwise product in $C^{0}(\partial\Omega)$,    the operator  $P_{ljr}[g,\cdot]$ is    continuous from  $C^{0}(\partial\Omega)$ to $C^{0}(\partial\Omega)$. In particular,  $Q[\frac{\partial S_{ {\mathbf{a}} }}{\partial x_{r}}\circ\Theta,g,\mu]$, $P_{ljr}[g,\mu]\in C^{0}(\partial\Omega)$.\par

We now show that the weak $M_{lj}$-derivative of $Q\left[\frac{\partial S_{ {\mathbf{a}} }}{\partial x_{r}}\circ\Theta,g,\cdot\right]$ in  $\partial\Omega$ coincides with $P_{ljr}[g,\mu]$.\par

By considering an extension of $\mu$ of class $C^{1}$ with compact support in ${\mathbb{R}}^{n}$ and by considering a sequence of mollifiers of such an extension, and by taking the restriction to $\partial\Omega$, we conclude that there exists a sequence of functions $\{\mu_{b}\}_{b\in{\mathbb{N}}  }$ in $C^{2}(\partial\Omega)$ converging to $\mu$ in $C^{1}(\partial\Omega)$. 
Then we note that if $\varphi\in C^{1}(\partial\Omega)$ the validity of (\ref{formula1}) for $\mu_{b}\in C^{2}(\partial\Omega)\subseteq  C^{1,\beta}(\partial\Omega)$, and the membership of $ Q\left[\frac{\partial S_{ {\mathbf{a}} }}{\partial x_{r}}\circ\Theta,g,\mu_{b}\right]$ in $C^{1}(\partial\Omega)$ (see Theorem \ref{qrssm} (ii)), and Lemma \ref{gagre} imply that
\begin{eqnarray*}
\lefteqn{
\int_{\partial\Omega}Q\left[\frac{\partial S_{ {\mathbf{a}} }}{\partial x_{r}}\circ\Theta,g,\mu\right]M_{lj}[\varphi]\,d\sigma=\lim_{b\to\infty}
\int_{\partial\Omega}Q\left[\frac{\partial S_{ {\mathbf{a}} }}{\partial x_{r}}\circ\Theta,g,\mu_{b}\right]M_{lj}[\varphi]\,d\sigma
}
\\
&&\qquad\qquad \qquad
=
-\lim_{b\to\infty}\int_{\partial\Omega}M_{lj}\left[Q\left[\frac{\partial S_{ {\mathbf{a}} }}{\partial x_{r}}\circ\Theta,g,\mu_{b}\right]\right]\varphi\,d\sigma
\\
&& \qquad\qquad\qquad
=
-\lim_{b\to\infty}\int_{\partial\Omega}P_{ljr}[g,\mu_{b}]\varphi\,d\sigma
=-\int_{\partial\Omega}P_{ljr}[g,\mu]\varphi\,d\sigma\,.
\end{eqnarray*}
Hence, $P_{ljr}[g,\mu]$ coincides with the weak $M_{lj}$-derivative of 
$Q\left[\frac{\partial S_{ {\mathbf{a}} }}{\partial x_{r}}\circ\Theta,g,\mu\right]$ for all $l,j\in\{1,\dots,n\}$. Since both $P_{ljr}[g,\mu]$ and 
$Q\left[\frac{\partial S_{ {\mathbf{a}} }}{\partial x_{r}}\circ\Theta,g,\mu\right]$ are continuous functions, it follows that $Q\left[\frac{\partial S_{ {\mathbf{a}} }}{\partial x_{r}}\circ\Theta,g,\mu\right]\in C^{1}(\partial\Omega)$ and that $M_{lj}\left[Q\left[\frac{\partial S_{ {\mathbf{a}} }}{\partial x_{r}}\circ\Theta,g,\mu\right]\right]=P_{ljr}[g,\mu]$ classically. Hence (\ref{formula1}) holds also for $\mu\in C^{1}(\partial\Omega)$.
\hfill  $\Box$ 

\vspace{\baselineskip}

By exploiting  formula (\ref{formula1}), we can prove the following.

\begin{thm}
\label{qrsm}
Let ${\mathbf{a}}$ be as in (\ref{introd0}), (\ref{ellip}), (\ref{symr}). Let $S_{ {\mathbf{a}} }$ be a fundamental solution of $P[{\mathbf{a}},D]$. 
 Let $\alpha\in]0,1[$. Let $m\in{\mathbb{N}}\setminus\{0\}$. 
Let $\Omega$ be a bounded open subset of ${\mathbb{R}}^{n}$ of class $C^{m,\alpha}$. Let $r\in\{1,\dots,n\}$. Then the following statements hold.
\begin{enumerate}
\item[(i)] Let $\theta\in]0,1]$. Then the
 bilinear map $Q\left[\frac{\partial S_{ {\mathbf{a}} }}{\partial x_{r}}\circ\Theta,\cdot,\cdot\right]$  from the space $C^{m-1,\theta}(\partial\Omega)\times C^{m-1}(\partial\Omega)$ to $C^{m-1,\omega_{\theta}(\cdot)}(\partial\Omega)$ which takes  a pair  $(g,\mu)$ to $Q\left[\frac{\partial S_{ {\mathbf{a}} }}{\partial x_{r}}\circ\Theta,g,\mu\right]$ is continuous. 
\item[(ii)] Let $\beta\in]0,1]$. Then the
 bilinear map $Q\left[\frac{\partial S_{ {\mathbf{a}} }}{\partial x_{r}}\circ\Theta,\cdot,\cdot\right]$  from  the space $C^{m-1,\alpha}(\partial\Omega)\times C^{m-1,\beta}(\partial\Omega)$ to $C^{m-1,\alpha}(\partial\Omega)$ which takes a pair   $(g,\mu)$ to $Q\left[\frac{\partial S_{ {\mathbf{a}} }}{\partial x_{r}}\circ\Theta,g,\mu\right]$ is continuous. 
\end{enumerate}
\end{thm}
{\bf Proof.} We first prove statement (i). We proceed by induction on $m$. Case $m=1$ holds by Theorem \ref{qrs} (i). We now prove that if the statement holds for $m$, then it holds for $m+1$. Thus we now assume that $\Omega$ is of class $C^{m+1,\alpha}$ and we turn to prove that $Q[\frac{\partial S_{ {\mathbf{a}} }}{\partial x_{r}}\circ\Theta,\cdot,\cdot]$ is bilinear and continuous from $C^{m,\theta}(\partial\Omega)\times C^{m}(\partial\Omega)$ to $C^{m,\omega_{\theta}(\cdot)}(\partial\Omega)$. By Lemma \ref{tanco} (ii), it suffices to prove that the following two statements hold.
\begin{enumerate}
\item[(j)] $Q[\frac{\partial S_{ {\mathbf{a}} }}{\partial x_{r}}\circ\Theta,\cdot,\cdot]$ is continuous from $C^{m,\theta}(\partial\Omega)\times C^{m}(\partial\Omega)$ to $C^{0}(\partial\Omega)$.
\item[(jj)] $M_{lj}[Q[\frac{\partial S_{ {\mathbf{a}} }}{\partial x_{r}}\circ\Theta,\cdot,\cdot]]$ is continuous from $C^{m,\theta}(\partial\Omega)\times C^{m}(\partial\Omega)$ to the space $C^{m-1,\omega_{\theta}(\cdot)}(\partial\Omega)$ for all $l$, $j\in\{1,\dots,n\}$. 
\end{enumerate}
Statement (j) holds by case $m=1$, and by the imbedding of $C^{m,\theta}(\partial\Omega)\times C^{m}(\partial\Omega)$ into $C^{0,\theta}(\partial\Omega)\times C^{0 }(\partial\Omega)$. We now prove statement (jj).  Since $m+1\geq 2$,  Lemma \ref{formula} and the inductive assumption imply   that we can actually apply $M_{lj}$ to $Q[\frac{\partial S_{ {\mathbf{a}} }}{\partial x_{r}}\circ\Theta,\cdot,\cdot]$. We find convenient to denote by $P_{ljr}[g,\mu]$ the right hand side of formula
(\ref{formula1}). Then we have 
\[
M_{lj}[Q[\frac{\partial S_{ {\mathbf{a}} }}{\partial x_{r}}\circ\Theta,g,\mu]]=P_{ljr}[g,\mu]\qquad\forall (g,\mu)\in
C^{m,\theta}(\partial\Omega)\times C^{m}(\partial\Omega)\,.
\]
By Lemma \ref{thgs}, and by the membership of $\nu$ in $C^{m,\alpha}(\partial\Omega, {\mathbb{R}}^{n})$, which is contained in $ C^{m-1, 1}(\partial\Omega, {\mathbb{R}}^{n})$,   and by the continuity of the pointwise product in Schauder spaces, and by the continuity of the imbedding of $C^{m}(\partial\Omega)$ into 
$C^{m-1}(\partial\Omega)$, and of $C^{m,\alpha}(\partial\Omega)$ into $C^{m-1,\theta}(\partial\Omega)$, and by the inductive assumption on the continuity of $Q\left[\frac{\partial S_{ {\mathbf{a}} }}{\partial x_{r}}\circ\Theta,\cdot,\cdot\right]$, and by the continuity of $v[\partial\Omega ,    S_{ {\mathbf{a}} },\cdot ]_{|\partial\Omega}$ from $C^{m-1,\alpha}(\partial\Omega)$ to $C^{m,\alpha}(\partial\Omega)\subseteq C^{m-1,\theta}(\partial\Omega)$,  
and by the continuity of the imbedding of 
$C^{m}(\partial\Omega)$ into $C^{m-1,\alpha}(\partial\Omega)$  
and of  $C^{m}(\partial\Omega)$ into 
$C^{m-1,\omega_{\theta}(\cdot)}(\partial\Omega)$,  and by the continuity of $D_{  {\mathbf{a}}  }$ from $C^{m,\theta}(\partial\Omega)$ to $C^{m-1,\theta}(\partial\Omega)$, we conclude that $P_{ljr}[\cdot,\cdot]$ is bilinear and continuous from
$C^{m,\theta }(\partial\Omega)\times C^{m}(\partial\Omega)$ to $C^{m-1,\omega_{\theta}(\cdot)}(\partial\Omega)$, and the proof of statement  (jj) and accordingly of statement (i) is complete. The proof of statement (ii) follows the lines of the proof of statement (i), by replacing the use of Theorem \ref{qrs} (i) with that of Theorem  \ref{qrs} (ii). \hfill  $\Box$ 

\vspace{\baselineskip}

\begin{defn}
\label{r}
Let ${\mathbf{a}}$ be as in (\ref{introd0}), (\ref{ellip}), (\ref{symr}). Let $S_{ {\mathbf{a}} }$ be a fundamental solution of $P[{\mathbf{a}},D]$. 
 Let $\alpha\in]0,1[$. Let $\Omega$ be a bounded open subset of ${\mathbb{R}}^{n}$ of class $C^{1,\alpha}$.  
 Then we set
\begin{eqnarray*}
\lefteqn{
R[g,h,\mu]
}
\\ \nonumber
&&
\equiv \sum_{r=1}a_{r} 
\left\{
Q[\frac{\partial S_{ {\mathbf{a}} }}{\partial x_{r}}\circ\Theta,gh,\mu]-g
Q[\frac{\partial S_{ {\mathbf{a}} }}{\partial x_{r}}\circ\Theta,h,\mu]
-Q[\frac{\partial S_{ {\mathbf{a}} }}{\partial x_{r}}\circ\Theta,h,g\mu]
\right\}
\\ \nonumber
&&\
+a\left\{
gv[\partial\Omega     , S_{ {\mathbf{a}} }   ,h\mu  ]
- h v[\partial\Omega     , S_{ {\mathbf{a}} }  ,g\mu  ]
\right\}\,
\end{eqnarray*} 
for all $ (g,h,\mu)\in (C^{0,\alpha}(\partial\Omega))^{2}\times C^{0}(\partial\Omega)$.
\end{defn}
 Since
\begin{eqnarray*}
\lefteqn{
g(x)h(y)-g(y)h(x)=[g(x)h(x)-g(y)h(y)]
}
\\ \nonumber
&&\qquad
-g(x)[h(x)-h(y)]-g(y)[h(x)-h(y)]\qquad\forall x,y\in\partial\Omega\,,
\end{eqnarray*}
we have 
\begin{eqnarray*}
\lefteqn{ 
 R[g,h,\mu]
 =\int_{\partial\Omega}
\left\{
\sum_{r=1}^{n}a_{r}
 \frac{\partial}{\partial x_{r}}S_{{\mathbf{a}}}(x-y)
 +aS_{{\mathbf{a}}}(x-y)
\right\}
}
\\ \nonumber
&&\qquad\qquad\qquad\qquad
\times [g(x)h(y)-g(y)h(x)]\mu(y)\,d\sigma_{y}\qquad\forall x\in\partial\Omega\,.
\end{eqnarray*} 
Since $R$ is a composition of the operator $Q[\frac{\partial S_{ {\mathbf{a}} }}{\partial x_{r}}\circ\Theta,\cdot,\cdot]$ and of a single layer potential, Theorems \ref{slay}, \ref{v0a} and Theorem \ref{qrsm} and the continuity of the product in Schauder spaces 
and of the imbedding of $C^{m-1}(\partial\Omega)$ into $C^{m-2,\alpha}(\partial\Omega)$
for $m\geq 2$ and of the imbedding of $C^{m-1,\alpha}(\partial\Omega)$ into 
$C^{m-1,\omega_{\alpha}(\cdot)}(\partial\Omega)$ and of 
$C^{m,\beta}(\partial\Omega)$ into 
$C^{m-1,\alpha }(\partial\Omega)$,  imply the validity of the  following. 
\begin{thm}
\label{rrsm} 
Let ${\mathbf{a}}$ be as in (\ref{introd0}), (\ref{ellip}), (\ref{symr}). Let $S_{ {\mathbf{a}} }$ be a fundamental solution of $P[{\mathbf{a}},D]$. 
 Let $\alpha\in]0,1[$. Let $m\in{\mathbb{N}}\setminus\{0\}$. Let $\Omega$ be a bounded open subset of ${\mathbb{R}}^{n}$ of class $C^{m,\alpha}$. Then the following statements hold. 
\begin{enumerate}
\item[(i)]   The
 trilinear map $R$  from the space $\left(C^{m-1,\alpha}(\partial\Omega)\right)^{2}\times C^{m-1}(\partial\Omega)$ to $C^{m-1,\omega_{\alpha}(\cdot)}(\partial\Omega)$ which takes a pair $(g,h,\mu)$ to $R[g,h,\mu]$ is continuous. 
\item[(ii)] Let $\beta\in]0,1]$. Then the
 trilinear map $R$  from the space  $\left(C^{m-1,\alpha}(\partial\Omega)\right)^{2}\times C^{m-1,\beta}(\partial\Omega)$ to $C^{m-1,\alpha}(\partial\Omega)$ which takes a pair $(g,h,\mu)$ to $R[g,h,\mu]$ is continuous. 
\end{enumerate}
\end{thm}

\section{Tangential derivatives and regularizing properties of the double layer potential}

We now exploit Theorems \ref{dlay}, \ref{w0a}  and Lemma \ref{formula}, and Theorems \ref{qrsm}, \ref{rrsm} in order to prove a formula for the tangential derivatives of the double layer potential, which generalizes the corresponding formula of 
 Hofmann, Mitrea and Taylor~\cite[(6.2.6)]{HoMitTa10} for homogeneous operators.  We do so by means of the following. 
\begin{thm}
\label{wtg} 
Let ${\mathbf{a}}$ be as in (\ref{introd0}), (\ref{ellip}), (\ref{symr}). Let $S_{ {\mathbf{a}} }$ be a fundamental solution of $P[{\mathbf{a}},D]$. 
 Let $\alpha\in]0,1[$.  Let $\Omega$ be a bounded open subset of ${\mathbb{R}}^{n}$ of class $C^{1,\alpha}$. If $\mu\in C^{1}(\partial\Omega)$, then $w[\partial\Omega     ,{\mathbf{a}}, S_{ {\mathbf{a}} },\mu]_{|\partial\Omega}\in 
 C^{1 }(\partial\Omega)$ and
 \begin{eqnarray}
\label{wtg1}
\lefteqn{
M_{lj}[w[\partial\Omega    ,{\mathbf{a}}, S_{ {\mathbf{a}} } ,\mu]_{|\partial\Omega}]
=w[\partial\Omega     ,{\mathbf{a}}, S_{ {\mathbf{a}} }     ,M_{lj}[\mu]    ]_{|\partial\Omega}
}
\\ \nonumber
&& \qquad\quad 
+\sum_{b,r=1}^{n}a_{br}
\left\{
Q\left[\frac{\partial S_{ {\mathbf{a}} }}{\partial x_{b}}\circ\Theta,\nu_{l},M_{jr}[\mu]\right]-Q\left[\frac{\partial S_{ {\mathbf{a}} }}{\partial x_{b}}\circ\Theta,\nu_{j}, M_{lr}[\mu]\right]
\right\}
\\ \nonumber
&&  \qquad\quad 
+\nu_{l} Q\left[\frac{\partial S_{ {\mathbf{a}} }}{\partial x_{j}}\circ\Theta,\nu\cdot a^{(1)},\mu\right] 
-\nu_{j} Q\left[\frac{\partial S_{ {\mathbf{a}} }}{\partial x_{l}}\circ\Theta,\nu\cdot a^{(1)},\mu\right] 
\\ \nonumber
&& \qquad\quad 
+\nu \cdot a^{(1)}
\left\{
Q\left[\frac{\partial S_{ {\mathbf{a}} }}{\partial x_{l}}\circ\Theta,\nu_{j},\mu\right] 
-
Q\left[\frac{\partial S_{ {\mathbf{a}} }}{\partial x_{j}}\circ\Theta,\nu_{l},\mu\right] 
\right\}
\\ \nonumber
&& \qquad\quad 
-\nu \cdot a^{(1)}
v[\partial\Omega      , S_{ {\mathbf{a}} }   , M_{lj}[\mu]]
+v[\partial\Omega   , S_{ {\mathbf{a}} }         , \nu \cdot a^{(1)}M_{lj}[\mu]]
\\ \nonumber
&&  \qquad\qquad\qquad\qquad
+R[\nu_{l},\nu_{j},\mu]
\qquad{\mathrm{on}}\ \partial\Omega\,,
\end{eqnarray}
for all $l,j\in\{1,\dots,n\}$. (For  $Q$ see (\ref{qrs0}).) 
\end{thm}
{\bf Proof.} We fix $\beta\in]0,\alpha[$ and we  first consider the specific case in which $\mu\in C^{1,\beta}(\partial\Omega)$. Let $R\in]0,+\infty[$ be such that ${\mathrm{cl}}\Omega\subseteq {\mathbb{B}}_{n}(0,R)$. Let  ` $\tilde{\ }$ '  be an extension operator of $C^{1,\beta}(\partial\Omega)$ to $C^{1,\beta}({\mathrm{cl}} {\mathbb{B}}_{n}(0,R))$  as in Lemma \ref{extsch}. By Theorem \ref{dlay} (i), (ii), we have
$w^{+}[\partial\Omega   ,{\mathbf{a}}, S_{ {\mathbf{a}} }  ,\mu ]\in 
C^{1,\beta}({\mathrm{cl}} \Omega)$  and
\begin{equation}
\label{wtg1a}
M_{lj}[w^{+}[\partial\Omega     ,{\mathbf{a}}, S_{ {\mathbf{a}} },\mu]_{|\partial\Omega}]
=\frac{1}{2}M_{lj}[\mu]+M_{lj}[w [\partial\Omega ,{\mathbf{a}}, S_{ {\mathbf{a}} },\mu]_{|\partial\Omega}]\,.
\end{equation}
By the definition of $M_{lj}$ and by  equality (\ref{dlay1}), we have 
\begin{eqnarray}
\label{wtg2}
\lefteqn{
M_{lj}[w^{+}[\partial\Omega ,{\mathbf{a}}, S_{ {\mathbf{a}} },\mu]_{|\partial\Omega}]
}
\\ \nonumber
&& 
=\nu_{l}\frac{\partial}{\partial x_{j}}w^{+}[\partial\Omega ,{\mathbf{a}}, S_{ {\mathbf{a}} },\mu]
-\nu_{j}\frac{\partial}{\partial x_{l}}w^{+}[\partial\Omega ,{\mathbf{a}}, S_{ {\mathbf{a}} },\mu]
\\ \nonumber
&& 
=\nu_{l}\biggl[
\sum_{b,r=1}^{n}a_{br}\frac{\partial}{\partial x_{b}}
v^{+}[\partial\Omega        ,S_{ {\mathbf{a}} }, M_{jr}[\mu]]
+\sum_{b=1}^{n}a_{b}\frac{\partial}{\partial x_{b}}
v^{+}[\partial\Omega , S_{ {\mathbf{a}} }, \nu_{j}\mu]
\biggr.
\\ \nonumber
&& \quad
\biggl.
-\frac{\partial}{\partial x_{j} }
v^{+}[\partial\Omega     ,S_{ {\mathbf{a}} }  ,(\nu^{t}\cdot a^{(1)})\mu]
+av^{+}[\partial\Omega ,S_{ {\mathbf{a}} },\nu_{j}\mu]
\biggr]
\\ \nonumber
&& \quad
-\nu_{j}\biggl[
\sum_{b,r=1}^{n}a_{br}\frac{\partial}{\partial x_{b}}
v^{+}[\partial\Omega ,S_{ {\mathbf{a}} }, M_{lr}[\mu]]
+\sum_{b=1}^{n}a_{b}\frac{\partial}{\partial x_{b}}
v^{+}[\partial\Omega , S_{ {\mathbf{a}} }, \nu_{l}\mu]
\biggr.
\\ \nonumber
&& \quad
\biggl.
-\frac{\partial}{\partial x_{l} }
v^{+}[\partial\Omega ,S_{ {\mathbf{a}} } ,(\nu^{t}\cdot a^{(1)})\mu]
+av^{+}[\partial\Omega ,S_{ {\mathbf{a}} },\nu_{l}\mu]
\biggr]
\\ \nonumber
&&\qquad
=\sum_{b,r=1}^{n}a_{br}
\left\{
\nu_{l}\frac{\partial}{\partial x_{b}}
v^{+}[\partial\Omega ,S_{ {\mathbf{a}} } , M_{jr}[\mu]]
-
\nu_{j}\frac{\partial}{\partial x_{b}}
v^{+}[\partial\Omega ,S_{ {\mathbf{a}} } , M_{lr}[\mu]]
\right\}
\\ \nonumber
&& \quad
+\sum_{b=1}^{n}a_{b}
\left\{
\nu_{l}\frac{\partial}{\partial x_{b}}
v^{+}[\partial\Omega , S_{ {\mathbf{a}} }, \nu_{j}\mu]
-
\nu_{j}\frac{\partial}{\partial x_{b}}
v^{+}[\partial\Omega , S_{ {\mathbf{a}} } , \nu_{l}\mu]
\right\}
\\ \nonumber
&& \quad
-
\left\{
\nu_{l}\frac{\partial}{\partial x_{j} }
v^{+}[\partial\Omega ,S_{ {\mathbf{a}} } ,(\nu^{t}\cdot a^{(1)})\mu]
-
\nu_{j}\frac{\partial}{\partial x_{l} }
v^{+}[\partial\Omega ,S_{ {\mathbf{a}} } ,(\nu^{t}\cdot a^{(1)})\mu]
\right\}
\\ \nonumber
&& \quad
+a\left\{
\nu_{l}v[\partial\Omega          ,S_{ {\mathbf{a}} }  , \nu_{j}\mu]
-
\nu_{j}v[\partial\Omega  ,S_{ {\mathbf{a}} }  ,\nu_{l}\mu   ]
\right\}\qquad{\mathrm{on}}\ \partial\Omega\,.
\end{eqnarray}
We now consider the first term  in braces in the right hand side of (\ref{wtg2}), and we note that
\begin{eqnarray}
\label{wtg3}
\lefteqn{
\left\{
\nu_{l}(x)\frac{\partial}{\partial x_{b}}
v^{+}[\partial\Omega ,S_{ {\mathbf{a}} }, M_{jr}[\mu]](x)
-
\nu_{j}\frac{\partial}{\partial x_{b}}
v^{+}[\partial\Omega ,S_{ {\mathbf{a}} }, M_{lr}[\mu]](x)
\right\}
}
\\ \nonumber
&&\quad
=
-\frac{\nu_{l}(x)\nu_{b}(x)}{ 2  \nu^{t}(x)  a^{(2)}\nu (x)     }M_{jr}[\mu](x)
+\nu_{l}(x)\int_{\partial\Omega}\frac{\partial}{\partial x_{b}}S_{{\mathbf{a}}}(x-y)M_{jr}[\mu](y)\,d\sigma_{y}
\\ \nonumber
&&\quad\quad
+
\frac{\nu_{j}(x)\nu_{b}(x)}{ 2  \nu^{t}(x)  a^{(2)}\nu (x)     }M_{lr}[\mu](x)
-\nu_{j}(x)\int_{\partial\Omega}\frac{\partial}{\partial x_{b}}S_{{\mathbf{a}}}(x-y)M_{lr}[\mu](y)\,d\sigma_{y}
\\ \nonumber
&&\quad
=\nu_{b}(x)\frac{
-\nu_{l}(x)M_{jr}[\mu](x)+\nu_{j}(x)M_{lr}[\mu](x)
}{   2  \nu^{t}(x)  a^{(2)}\nu (x)   }
\\ \nonumber
&&\quad\quad
+\int_{\partial\Omega}\frac{\partial}{\partial x_{b}}
S_{{\mathbf{a}}}(x-y)
\{
\nu_{l}(x)M_{jr}[\mu](y)
-
\nu_{j}(x)M_{lr}[\mu](y)
\}\,d\sigma_{y}\,.
\end{eqnarray}
Next we note that
\begin{eqnarray}
\label{wtg4}
\lefteqn{
[
\nu_{l}M_{jr}[\mu]-\nu_{j}M_{lr}[\mu]
]
}
\\
\nonumber
&&
\quad
=\nu_{l}\nu_{j}\frac{\partial\mu}{\partial x_{r}}
-\nu_{l}\nu_{r}\frac{\partial\mu}{\partial x_{j}}
-\nu_{j}\nu_{l}\frac{\partial\mu}{\partial x_{r}}
+\nu_{j}\nu_{r}\frac{\partial\mu}{\partial x_{l}}
=-\nu_{r}M_{lj}[\mu]\qquad{\mathrm{on}}\ \partial\Omega\,.
\end{eqnarray}
Then we obtain
\begin{eqnarray}
\label{wtg5}
\lefteqn{
\sum_{b,r=1}^{n}a_{br}\nu_{b} \frac{
-\nu_{l} M_{jr}[\mu] +\nu_{j} M_{lr}[\mu] 
}{   2  \nu^{t}   a^{(2)}\nu     }
}
\\
\nonumber
&&\quad
=\sum_{b,r=1}^{n}a_{br}\nu_{b}\frac{  \nu_{r}M_{lj}[\mu]  }{2    \nu^{t}  a^{(2)}\nu   }
=\frac{
\sum_{b,r=1}^{n}\nu_{b}a_{br}\nu_{r}
}{2    \nu^{t}  a^{(2)}\nu}M_{lj}[\mu]=\frac{1}{2}M_{lj}[\mu]
\qquad{\mathrm{on}}\ \partial\Omega\,.
\end{eqnarray}
Next we consider the term in braces  in the argument of the integral in the right hand side of (\ref{wtg3}), and we note that equality (\ref{wtg4}) implies that
\begin{eqnarray}
\label{wtg6}
\lefteqn{
 \nu_{l}(x)M_{jr}[\mu](y)-\nu_{j}(x)M_{lr}[\mu](y)
 }
 \\   
 \nonumber
 &&\qquad\qquad
=
[\nu_{l}(x) -\nu_{l}(y) ]M_{jr}[\mu](y)
+
[
\nu_{l}(y)M_{jr}[\mu](y)-\nu_{j}(y)M_{lr}[\mu](y)
]
\\   
 \nonumber
 && \qquad\qquad\quad 
-[\nu_{j}(x)-\nu_{j}(y)  ] M_{lr}[\mu](y)
\\   
 \nonumber
 && \qquad\qquad
=
[\nu_{l}(x) -\nu_{l}(y) ]M_{jr}[\mu](y)
-\nu_{r}(y)M_{lj}[\mu](y)
\\   
 \nonumber
 && \qquad\qquad\quad 
-[\nu_{j}(x)-\nu_{j}(y)  ] M_{lr}[\mu](y)
\qquad \forall x,y\in\partial\Omega\,.
\end{eqnarray}
Next we consider the term in the second braces in the right hand side of equality 
(\ref{wtg2}) and we note that
 \begin{eqnarray}
\label{wtg7}
\lefteqn{
\nu_{l}(x)\frac{\partial}{\partial x_{b}}
v^{+}[\partial\Omega , S_{ {\mathbf{a}} }, \nu_{j}\mu](x)
-
\nu_{j}(x)\frac{\partial}{\partial x_{b}}
v^{+}[\partial\Omega , S_{ {\mathbf{a}} }, \nu_{l}\mu](x)
}
\\ \nonumber
&& 
=-\nu_{l}(x)\frac{\nu_{b}(x)}{   2  \nu^{t}(x)  a^{(2)}\nu (x)   }\nu_{j}(x)\mu(x)
+\nu_{l}(x)\int_{\partial\Omega}\frac{\partial}{\partial x_{b}}S_{{\mathbf{a}}}(x-y)\nu_{j}(y)\mu(y)\,d\sigma_{y}
\\ \nonumber
&& \quad
+ \nu_{j}(x)\frac{\nu_{b}(x)}{   2  \nu^{t}(x)  a^{(2)}\nu (x)   }\nu_{l}(x)\mu(x)
-\nu_{j}(x)\int_{\partial\Omega}\frac{\partial}{\partial x_{b}}S_{{\mathbf{a}}}(x-y)\nu_{l}(y)\mu(y)\,d\sigma_{y}
\\ \nonumber
&& 
=
\int_{\partial\Omega}\frac{\partial}{\partial x_{b}}S_{{\mathbf{a}}}(x-y)
[
\nu_{l}(x)\nu_{j}(y)-\nu_{j}(x)\nu_{l}(y)
]
\mu(y)\,d\sigma_{y}
\qquad\forall x\in \partial\Omega\,.
\end{eqnarray}
Next we consider the term in the third braces in the right hand side of equality 
(\ref{wtg2}) and we note that
\begin{eqnarray}
\nonumber
\lefteqn{
\nu_{l}(x)\frac{\partial}{\partial x_{j} } 
v^{+}[\partial\Omega ,S_{ {\mathbf{a}} },(\nu^{t}\cdot a^{(1)})\mu](x)
-
\nu_{j}(x)\frac{\partial}{\partial x_{l} } 
v^{+}[\partial\Omega ,S_{ {\mathbf{a}} },(\nu^{t}\cdot a^{(1)})\mu](x)
}
\\ \label{wtg8}
&& 
=-\nu_{l}(x)\frac{\nu_{j}(x)}{   2  \nu^{t}(x)  a^{(2)}\nu (x)   }
(\nu^{t}(x)\cdot a^{(1)})\mu(x)
\\ \nonumber
&&\quad
+\nu_{l}(x)
\int_{\partial\Omega}\frac{\partial}{\partial x_{j}}S_{{\mathbf{a}}}(x-y)\nu^{t}(y)\cdot a^{(1)}\mu(y)\,d\sigma_{y}
\\ \nonumber
&&\quad
+\nu_{j}(x)\frac{\nu_{l}(x)}{   2  \nu^{t}(x)  a^{(2)}\nu (x)   }
(\nu^{t}(x)\cdot a^{(1)})\mu(x)
\\ \nonumber
&&\quad
-\nu_{j}(x)
\int_{\partial\Omega}\frac{\partial}{\partial x_{l}}S_{{\mathbf{a}}}(x-y)\nu^{t}(y)\cdot a^{(1)}\mu(y)\,d\sigma_{y}
\\ \nonumber
&& 
=-\nu_{l}(x)\int_{\partial\Omega}\left[
(\nu^{t}(x)\cdot a^{(1)})-(\nu^{t}(y)\cdot a^{(1)})
\right]\frac{\partial}{\partial x_{j}}S_{{\mathbf{a}}}(x-y)\mu(y)\,d\sigma_{y}
\\ \nonumber
&&\quad
+\nu_{l}(x)\int_{\partial\Omega}(\nu^{t}(x)\cdot a^{(1)})
\frac{\partial}{\partial x_{j}}S_{{\mathbf{a}}}(x-y)\mu(y)\,d\sigma_{y}
\\ \nonumber
&&\quad
+\nu_{j}(x)\int_{\partial\Omega}\left[
(\nu^{t}(x)\cdot a^{(1)})-(\nu^{t}(y)\cdot a^{(1)})
\right]\frac{\partial}{\partial x_{l}}S_{{\mathbf{a}}}(x-y)\mu(y)\,d\sigma_{y}
\\ \nonumber
&&\quad
-\nu_{j}(x)\int_{\partial\Omega}(\nu^{t}(x)\cdot a^{(1)})
\frac{\partial}{\partial x_{l}}S_{{\mathbf{a}}}(x-y)\mu(y)\,d\sigma_{y}
\\ \nonumber
&& 
=-\nu_{l}(x)\int_{\partial\Omega}\left[
(\nu^{t}(x)\cdot a^{(1)})-(\nu^{t}(y)\cdot a^{(1)})
\right]\frac{\partial}{\partial x_{j}}S_{{\mathbf{a}}}(x-y)\mu(y)\,d\sigma_{y}
\\ \nonumber
&&\quad 
+\nu_{j}(x)\int_{\partial\Omega}\left[
(\nu^{t}(x)\cdot a^{(1)})-(\nu^{t}(y)\cdot a^{(1)})
\right]\frac{\partial}{\partial x_{l}}S_{{\mathbf{a}}}(x-y)\mu(y)\,d\sigma_{y}
\\ \nonumber
&&\quad 
+(\nu^{t}(x)\cdot a^{(1)})\int_{\partial\Omega}
\left(
\nu_{l}(x)\frac{\partial}{\partial x_{j}}
-
\nu_{j}(x)\frac{\partial}{\partial x_{l}}
\right)S_{{\mathbf{a}}}(x-y)\mu(y)\,d\sigma_{y}
\\ \nonumber
&& 
=-\nu_{l}(x)\int_{\partial\Omega}\left[
(\nu^{t}(x)\cdot a^{(1)})-(\nu^{t}(y)\cdot a^{(1)})
\right]\frac{\partial}{\partial x_{j}}S_{{\mathbf{a}}}(x-y)\mu(y)\,d\sigma_{y}
\\ \nonumber
&&\quad 
+\nu_{j}(x)\int_{\partial\Omega}\left[
(\nu^{t}(x)\cdot a^{(1)})-(\nu^{t}(y)\cdot a^{(1)})
\right]\frac{\partial}{\partial x_{l}}S_{{\mathbf{a}}}(x-y)\mu(y)\,d\sigma_{y}
\\ \nonumber
&&\quad 
+(\nu^{t}(x)\cdot a^{(1)})\biggl\{\biggr.\int_{\partial\Omega}
(\nu_{l}(x)-\nu_{l}(y))\frac{\partial}{\partial x_{j}}S_{{\mathbf{a}}}(x-y)\mu(y)\,d\sigma_{y}
\\ \nonumber
&&\quad 
-\int_{\partial\Omega}
(\nu_{j}(x)-\nu_{j}(y))\frac{\partial}{\partial x_{l}}S_{{\mathbf{a}}}(x-y)\mu(y)\,d\sigma_{y}
\biggl.\biggr\}
\\ \nonumber
&&\quad
+(\nu^{t}(x)\cdot a^{(1)})\int_{\partial\Omega}
\left(
\nu_{l}(y)\frac{\partial}{\partial x_{j}}
-
\nu_{j}(y)\frac{\partial}{\partial x_{l}}
\right)S_{{\mathbf{a}}}(x-y)\mu(y)\,d\sigma_{y}\,,
\end{eqnarray}
for all $x\in \partial\Omega$. 
By Lemma \ref{gagre}, the last integral in the right hand side of (\ref{wtg8}) equals
\begin{equation}
\label{wtg9}
-\int_{\partial\Omega}M_{lj,y}[S_{{\mathbf{a}}}(x-y)]\mu(y)\,d\sigma_{y}
=\int_{\partial\Omega}S_{{\mathbf{a}}}(x-y)M_{lj }[\mu](y)\,d\sigma_{y}
\qquad\forall x\in \partial\Omega\,.
\end{equation}
Thus the last term in the right hand side of (\ref{wtg8}) equals
\begin{eqnarray}
\label{wtg10}
\lefteqn{
(\nu^{t}(x)\cdot a^{(1)})\int_{\partial\Omega}
S_{{\mathbf{a}}}(x-y)M_{lj }[\mu](y)\,d\sigma_{y}
}
\\ \nonumber
&&\qquad
=\int_{\partial\Omega}\left[
(\nu^{t}(x)\cdot a^{(1)})-(\nu^{t}(y)\cdot a^{(1)})
\right]S_{{\mathbf{a}}}(x-y)M_{lj }[\mu](y)\,d\sigma_{y}
\\ \nonumber
&&\qquad\quad
+\int_{\partial\Omega}(\nu^{t}(y)\cdot a^{(1)})S_{{\mathbf{a}}}(x-y)M_{lj }[\mu](y)\,d\sigma_{y}
\qquad\forall x\in \partial\Omega\,.
\end{eqnarray}
In the fourth and last term in braces of equation (\ref{wtg2}), we have
\begin{equation}
\label{wtg11}
\int_{\partial\Omega}S_{{\mathbf{a}}}(x-y)\left[
\nu_{l}(x)\nu_{j}(y)-\nu_{j}(x)\nu_{l}(y)
\right]\mu(y)\,d\sigma_{y}
\qquad\forall x\in \partial\Omega\,.
\end{equation}
Then by combining (\ref{wtg1a})--(\ref{wtg3}),    (\ref{wtg5})--(\ref{wtg11}), we obtain
\begin{eqnarray*}
\lefteqn{
M_{lj}[w[\partial\Omega ,{\mathbf{a}}, S_{ {\mathbf{a}} },\mu](x)
}
\\ \nonumber
&&\qquad
=\sum_{b,r=1}^{n}a_{br}\biggl\{\biggr.
\int_{\partial\Omega}(\nu_{l}(x)-\nu_{l}(y))
\frac{\partial}{\partial x_{b}}S_{ {\mathbf{a}} }(x-y)
M_{jr}[\mu](y)\,d\sigma_{y}
\\ \nonumber
&&\qquad\quad
-\int_{\partial\Omega}(\nu_{j}(x)-\nu_{j}(y))
\frac{\partial}{\partial x_{b}}S_{ {\mathbf{a}} }(x-y)
M_{lr}[\mu](y)\,d\sigma_{y}
\\ \nonumber
&&\qquad\quad
-\int_{\partial\Omega}\nu_{r}(y)
\frac{\partial}{\partial x_{b}}S_{ {\mathbf{a}} }(x-y)
M_{lj}[\mu](y)\,d\sigma_{y}\biggl.\biggr\}
\\ \nonumber
&&\qquad\quad
+\sum_{b=1}^{n}a_{b}\int_{\partial\Omega}\frac{\partial}{\partial x_{b}}S_{ {\mathbf{a}} }(x-y)
[ \nu_{l}(x)\nu_{j}(y)-\nu_{j}(x)\nu_{l}(y) ]\mu(y)\,d\sigma_{y}
\\ \nonumber
&&\qquad\quad
+\nu_{l}(x)\int_{\partial\Omega}
[(\nu^{t}(x)\cdot a^{(1)})-(\nu^{t}(y)\cdot a^{(1)})]\frac{\partial}{\partial x_{j}}S_{ {\mathbf{a}} }(x-y)
\mu(y)\,d\sigma_{y}
\\ \nonumber
&&\qquad\quad
-\nu_{j}(x)\int_{\partial\Omega}
[(\nu^{t}(x)\cdot a^{(1)})-(\nu^{t}(y)\cdot a^{(1)})]\frac{\partial}{\partial x_{l}}S_{ {\mathbf{a}} }(x-y)
\mu(y)\,d\sigma_{y}
\\ \nonumber
&&\qquad\quad
-(\nu^{t}(x)\cdot a^{(1)})\biggl\{\biggr.\int_{\partial\Omega}
(\nu_{l}(x)-\nu_{l}(y))\frac{\partial}{\partial x_{j}}S_{ {\mathbf{a}} }(x-y)\mu(y)\,d\sigma_{y}
\\ \nonumber
&&\qquad\quad
-
\int_{\partial\Omega}
(\nu_{j}(x)-\nu_{j}(y))\frac{\partial}{\partial x_{l}}S_{ {\mathbf{a}} }(x-y)\mu(y)\,d\sigma_{y}
\biggl.\biggr\}
\\ \nonumber
&&\qquad\quad
-
\int_{\partial\Omega}[(\nu^{t}(x)\cdot a^{(1)})-(\nu^{t}(y)\cdot a^{(1)})]S_{ {\mathbf{a}} }(x-y)
M_{lj}[\mu](y)\,d\sigma_{y}
\\ \nonumber
&&\qquad\quad
-
\int_{\partial\Omega}(\nu^{t}(y)\cdot a^{(1)})S_{ {\mathbf{a}} }(x-y)
M_{lj}[\mu](y)\,d\sigma_{y}
\\ \nonumber
&&\qquad\quad
+a\int_{\partial\Omega}S_{ {\mathbf{a}} }(x-y)
[ \nu_{l}(x)\nu_{j}(y)-\nu_{j}(x)\nu_{l}(y) ]\mu(y)\,d\sigma_{y}
\qquad\forall x\in \partial\Omega\,,
\end{eqnarray*}
which we rewrite as
\begin{eqnarray*}
\lefteqn{
M_{lj}[    w[\partial\Omega ,{\mathbf{a}}, S_{ {\mathbf{a}} },\mu]](x)
}
\\ \nonumber
&&
=\sum_{b,r=1}^{n}a_{br}\biggl\{\biggr.
Q\left[\frac{\partial S_{ {\mathbf{a}} }}{\partial x_{b}}\circ\Theta,\nu_{l},M_{jr}[\mu]\right] (x)
-
Q\left[\frac{\partial S_{ {\mathbf{a}} }}{\partial x_{b}}\circ\Theta,\nu_{j},M_{lr}[\mu]\right] (x)\biggl.\biggr\}
\\ \nonumber
&& \ 
+\nu_{l}(x)Q\left[\frac{\partial S_{ {\mathbf{a}} }}{\partial x_{j}}\circ\Theta, \nu^{t}\cdot a^{(1)} ,\mu\right](x)
-\nu_{j}(x)Q\left[\frac{\partial S_{ {\mathbf{a}} }}{\partial x_{l}}\circ\Theta,  \nu^{t}\cdot a^{(1)} ,\mu\right](x)
\\ \nonumber
&& \ 
+w[\partial\Omega ,{\mathbf{a}}, S_{ {\mathbf{a}} }  ,M_{lj}[\mu] ](x)
\\ \nonumber
&& \ 
+(\nu^{t}(x)\cdot a^{(1)})\left\{
Q\left[\frac{\partial S_{ {\mathbf{a}} }}{\partial x_{l}}\circ\Theta, \nu_{j},\mu\right](x)-Q\left[\frac{\partial S_{ {\mathbf{a}} }}{\partial x_{j}}\circ\Theta,  \nu_{l},\mu\right](x)
\right\}
\\ \nonumber
&& \
-(\nu^{t}(x)\cdot a^{(1)})v[\partial\Omega       ,S_{ {\mathbf{a}} }  ,M_{lj}[\mu]](x)
+ v[\partial\Omega         ,S_{ {\mathbf{a}} }    ,(\nu^{t} \cdot a^{(1)})M_{lj}[\mu]](x)
\\ \nonumber
&& \
+R[\nu_{l},\nu_{j},\mu](x)
\qquad\forall x\in \partial\Omega\,.
\end{eqnarray*}
Thus we have proved formula (\ref{wtg1}) for $\mu\in C^{1,\beta}(\partial\Omega)$. 

Next we assume that $\mu\in C^{1}(\partial\Omega)$.
We denote by $T_{lj}[\mu]$ the right hand side of (\ref{wtg1}). By 
the continuity of $M_{lj}$ from $ C^{1}(\partial\Omega)$ to $ C^{0}(\partial\Omega)$, and by the continuity of $w[\partial\Omega ,{\mathbf{a}}, S_{ {\mathbf{a}} },\cdot]_{|\partial\Omega}$ 
and of $v[\partial\Omega   , S_{ {\mathbf{a}} }     ,\cdot]_{|\partial\Omega}$ from  $ C^{0}(\partial\Omega)$ to  $ C^{0,\alpha}(\partial\Omega)$, and by the continuity of 
$Q[\frac{\partial S_{ {\mathbf{a}} }}{\partial x_{r}}\circ\Theta,\cdot,\cdot]$ from $C^{0,\alpha}(\partial\Omega)\times
C^{0 }(\partial\Omega)$ to $C^{0,\omega_{\alpha}}(\partial\Omega)$, and by the continuity of $R$ from $\left(C^{0,\alpha}(\partial\Omega)\right)^{2}\times
C^{0 }(\partial\Omega)$ to $C^{0,\omega_{\alpha}}(\partial\Omega)$, and by the continuity of the pointwise product in Schauder spaces,  we conclude that the operators $w[\partial\Omega ,{\mathbf{a}}, S_{ {\mathbf{a}} },\cdot]_{|\partial\Omega}$ and $T_{lj}[\cdot]$ are
continuous from  $C^{1}(\partial\Omega)$ to $C^{0,\alpha}(\partial\Omega)$
and from $C^{1}(\partial\Omega)$ to $C^{0,\omega_{\alpha}(\cdot)}(\partial\Omega)$, respectively. 
In particular, $T_{lj}[\mu]$ and $w[\partial\Omega ,{\mathbf{a}}, S_{ {\mathbf{a}} },\mu]_{|\partial\Omega}$
belong to $C^{0}(\partial\Omega)$.
We now show that the weak $M_{lj}$-derivative of 
$w[\partial\Omega ,{\mathbf{a}}, S_{ {\mathbf{a}} },\mu]_{|\partial\Omega}$   coincides with $T_{lj}[\mu]$.\par

By arguing so as  at the end of the proof of Lemma \ref{formula}, 
there exists a sequence of functions $\{\mu_{b}\}_{b\in{\mathbb{N}}  }$ in $C^{1,\alpha}(\partial\Omega)$, which converges  to $\mu$ in $C^{1}(\partial\Omega)$. 
Then we note that if $\varphi\in C^{1}(\partial\Omega)$ the validity of (\ref{wtg1}) for $\mu_{b}\in C^{1,\alpha}(\partial\Omega)$,
and the membership of $w[\partial\Omega ,{\mathbf{a}}, S_{ {\mathbf{a}} },\mu_{b}]_{|\partial\Omega}$ in $C^{1,\alpha}(\partial\Omega)$, and   the above mentioned continuity of $w[\partial\Omega ,{\mathbf{a}}, S_{ {\mathbf{a}} },\cdot]_{|\partial\Omega}$, and Lemma \ref{gagre}  imply that
\begin{eqnarray*}
\lefteqn{
\int_{\partial\Omega}w[\partial\Omega ,{\mathbf{a}}, S_{ {\mathbf{a}} },\mu]_{|\partial\Omega}
M_{lj}[\varphi]\,d\sigma=\lim_{b\to\infty}
\int_{\partial\Omega}w[\partial\Omega ,{\mathbf{a}}, S_{ {\mathbf{a}} },\mu_{b}]_{|\partial\Omega}M_{lj}[\varphi]\,d \sigma
}
\\
&&\qquad\qquad\qquad\qquad
=
-\lim_{b\to\infty}\int_{\partial\Omega}M_{lj}[w[\partial\Omega ,{\mathbf{a}}, S_{ {\mathbf{a}} },\mu_{b}]_{|\partial\Omega}]\varphi\,d\sigma
\\
&&\qquad\qquad\qquad \qquad
=
-\lim_{b\to\infty}\int_{\partial\Omega}T_{lj}[\mu_{b}]\varphi\,dx
=-\int_{\partial\Omega}T_{lj}[\mu]\varphi\,dx\,.
\end{eqnarray*}
Hence, $T_{lj}[\mu]$ coincides with the weak $M_{lj}$-derivative of 
$w[\partial\Omega ,{\mathbf{a}}, S_{ {\mathbf{a}} },\mu]_{|\partial\Omega}$ for all $l,j$ in $\{1,\dots,n\}$. Since both $T_{lj}[\mu]$ and 
$w[\partial\Omega ,{\mathbf{a}}, S_{ {\mathbf{a}} },\mu]_{|\partial\Omega}$ are continuous functions, it follows that $w[\partial\Omega ,{\mathbf{a}}, S_{ {\mathbf{a}} },\mu]_{|\partial\Omega}\in C^{1}(\partial\Omega)$ and that $M_{lj}[w[\partial\Omega ,{\mathbf{a}}, S_{ {\mathbf{a}} },\mu]_{|\partial\Omega}]=T_{lj}[\mu]$ classically. Hence (\ref{wtg1}) holds also for $\mu\in C^{1}(\partial\Omega)$.\hfill  $\Box$ 

\vspace{\baselineskip}

By exploiting formula (\ref{wtg1}), we now prove the following result, which says that the double layer potential on $\partial\Omega$ has a regularizing effect. 
\begin{thm}
\label{wreg} 
Let ${\mathbf{a}}$ be as in (\ref{introd0}), (\ref{ellip}), (\ref{symr}). Let $S_{ {\mathbf{a}} }$ be a fundamental solution of $P[{\mathbf{a}},D]$. 
 Let $\alpha\in]0,1[$. Let $m\in{\mathbb{N}}\setminus\{0\}$. 
Let $\Omega$ be a bounded open subset of ${\mathbb{R}}^{n}$ of class $C^{m,\alpha}$. Then the following statements hold.
\begin{enumerate}
\item[(i)] The operator $w[\partial\Omega ,{\mathbf{a}}, S_{ {\mathbf{a}} },\cdot]_{|\partial\Omega}$ is linear and continuous from $C^{m}(\partial\Omega)$ to $C^{m,\omega_{\alpha}(\cdot)}(\partial\Omega)$.
\item[(ii)] Let $\beta\in]0,\alpha]$. Then the operator $w[\partial\Omega ,{\mathbf{a}}, S_{ {\mathbf{a}} },\cdot]_{|\partial\Omega}$ is linear and continuous from $C^{m,\beta}(\partial\Omega)$ to $C^{m,\alpha}(\partial\Omega)$.
\end{enumerate}
\end{thm}
{\bf Proof.} We  prove statement (i) by induction on $m$. As in the previous proof, we denote by $T_{lj}[\mu]$ the right hand side of formula (\ref{wtg1}).  We first consider case $m=1$. By Lemma \ref{tanco} (ii) and formula (\ref{wtg1}), it suffices to prove that the following two statements hold.
\begin{enumerate}
\item[(j)] $w[\partial\Omega ,{\mathbf{a}}, S_{ {\mathbf{a}} },\cdot]_{|\partial\Omega}$ is continuous from $C^{1}(\partial\Omega)$ to $C^{0}(\partial\Omega)$.
\item[(jj)] $T_{lj}[\cdot]$ is continuous from $C^{1}(\partial\Omega)$ to $C^{0,\omega_{\alpha}(\cdot)}(\partial\Omega)$ for all $l$, $j\in\{1,\dots,n\}$. 
\end{enumerate}
Theorem \ref{w0a}   implies the validity of (j). Statement (jj) follows by the continuity of the pointwise product in Schauder spaces and 
 by the continuity of $M_{lj}$ from $C^{1}(\partial\Omega)$ to $C^{0}(\partial\Omega)$, and by the continuity of 
$v[\partial\Omega, S_{ {\mathbf{a}} },\cdot]_{|\partial\Omega} $ and of 
$w[\partial\Omega ,{\mathbf{a}}, S_{ {\mathbf{a}} },\cdot]_{|\partial\Omega}$ from $C^{0}(\partial\Omega)$ to $C^{0,\alpha}(\partial\Omega)$ (cf. Theorems \ref{v0a}, \ref{w0a}), and by the continuity of   $Q[\frac{\partial S_{ {\mathbf{a}} }}{\partial x_{r}}\circ\Theta,\cdot,\cdot]$  from $C^{0,\alpha}(\partial\Omega)\times C^{0}(\partial\Omega)$ to $C^{0,\omega_{\alpha}(\cdot)}(\partial\Omega)$ (cf.  Theorem \ref{qrs} (i)) and by the continuity of $R$ from $\left(C^{0,\alpha}(\partial\Omega)\right)^{2}\times C^{0}(\partial\Omega)$ to $C^{0,
\omega_{\alpha}(\cdot)
}(\partial\Omega)$ (cf. Theorem \ref{rrsm} (i).)  \par 

 Next we assume
that  $\Omega$ is of class $C^{m+1,\alpha}$ and we turn to prove that $w[\partial\Omega ,{\mathbf{a}}, S_{ {\mathbf{a}} },\cdot]_{|\partial\Omega}$ is continuous from $C^{m+1}(\partial\Omega)$ to $C^{m+1,\omega_{\alpha}(\cdot)}(\partial\Omega)$.
By Lemma \ref{tanco} (ii) and formula (\ref{wtg1}), it suffices to prove that the following two statements hold.
\begin{enumerate}
\item[(a)] $w[\partial\Omega ,{\mathbf{a}}, S_{ {\mathbf{a}} },\cdot]_{|\partial\Omega}$ is continuous from $C^{m+1}(\partial\Omega)$ to $C^{0}(\partial\Omega)$.
\item[(b)] $T_{lj}[\cdot]$ is continuous from $C^{m+1}(\partial\Omega)$ to $C^{m,\omega_{\alpha}(\cdot)}(\partial\Omega)$. for all $l$, $j\in\{1,\dots,n\}$. 
\end{enumerate}
Statement (a) holds by the inductive assumption. We now prove statement (b).  Since $\Omega$ is of class $C^{m+1,\alpha}$, then $\nu$ is of class $C^{m,\alpha}(\partial\Omega)$.  Then 
Theorem \ref{qrsm} (i) ensures that $Q\left[\frac{\partial S_{ {\mathbf{a}} }}{\partial x_{r}}\circ\Theta,
\nu \cdot a^{(1)},\cdot\right]$ and $Q\left[\frac{\partial S_{ {\mathbf{a}} }}{\partial x_{r}}\circ\Theta,
\nu_{j},\cdot\right]$ are continuous from $C^{m}(\partial\Omega)$ to $C^{m,\omega_{\alpha}}(\partial\Omega)$ for all $l$, $j$, $r$ in  $\{1,\dots,n\}$.  Since $M_{lj}$ is continuous from $C^{m+1}(\partial\Omega)$ to 
$C^{m}(\partial\Omega)$, the inductive assumption implies that $w[\partial\Omega ,{\mathbf{a}}, S_{ {\mathbf{a}} },M_{lj}[\cdot]]_{|\partial\Omega} 
 $ is   continuous from $C^{m+1}(\partial\Omega)$ to 
$C^{m, \omega_{\alpha}(\cdot)}(\partial\Omega)$ for all $l$, $j$ in  $\{1,\dots,n\}$. 

Since $M_{lj}$ is continuous from $C^{m+1}(\partial\Omega)$ to 
 $C^{m-1,\alpha}(\partial\Omega)$ and $v[\partial\Omega ,  S_{ {\mathbf{a}} }        ,\cdot  ]_{|\partial\Omega}$ is continuous from $C^{m-1,\alpha}(\partial\Omega)$ to $C^{m,\alpha}(\partial\Omega)$
and $\nu\in (C^{m,\alpha}(\partial\Omega))^{n}$ and $C^{m,\alpha}(\partial\Omega)$ is continuously imbedded into $C^{m,\omega_{\alpha}(\cdot)}(\partial\Omega)$, we conclude that 
 $v[\partial\Omega      , S_{ {\mathbf{a}} }     ,M_{lj}[\cdot] ]_{|\partial\Omega}$ and 
 $v[\partial\Omega    , S_{ {\mathbf{a}} }      ,\nu\cdot a^{(1)}M_{lj}[\cdot] ]_{|\partial\Omega}$  are continuous from the space  $C^{m+1}(\partial\Omega)$ to $C^{m,\omega_{\alpha}(\cdot)}(\partial\Omega)$ for all $l$, $j$ in $\{1,\dots,n\}$. Moreover, $R$ is continuous  from $\left(C^{m,\alpha}(\partial\Omega)\right)^{2}\times C^{m}(\partial\Omega)$ to $C^{m,
\omega_{\alpha}(\cdot)
}(\partial\Omega)$ (cf. Theorem \ref{rrsm} (i).) Then   statement (b) holds true. \par

Statement (iii) can be proved by the same argument of the proof of statement (i) by exploiting   Theorem \ref{qrsm} (ii) instead of Theorem \ref{qrsm} (i) and Theorem \ref{rrsm} (ii) instead of Theorem \ref{rrsm} (i). \hfill  $\Box$ 

\vspace{\baselineskip}

Since $C^{m,\omega_{\alpha}(\cdot)}(\partial\Omega)$ is compactly imbedded into $C^{m}(\partial\Omega)$, and $C^{m,\alpha}(\partial\Omega)$ is compactly imbedded into
$C^{m,\beta}(\partial\Omega)$ for all $\beta\in]0,\alpha[$, we have the followng immediate consequence of Theorem \ref{wreg}. 

\begin{corol}
Under the assumptions of Theorem \ref{wreg}, the linear operator  $w[\partial\Omega ,{\mathbf{a}}, S_{ {\mathbf{a}} },\cdot]_{|\partial\Omega}$ is compact from $C^{m }(\partial\Omega)$ to itself, and from $C^{m,\omega_{\alpha}(\cdot)}(\partial\Omega)$ to itself, and from $C^{m,\alpha}(\partial\Omega)$ to itself. 
\end{corol}

\section{Other layer potentials  associated to $P[{\mathbf{a}},D]$}

Another relevant layer potential operator associated to the analysis of boundary value problems for the operator $P[{\mathbf{a}},D]$ is the following 
\[
w_{*}[\partial\Omega ,{\mathbf{a}}, S_{ {\mathbf{a}} },\mu](x)\equiv
\int_{\partial\Omega}\mu(y)DS_{ {\mathbf{a}} }(x-y)a^{(2)}\nu(x)\,d\sigma_{y}\qquad\forall x\in\partial\Omega\,,
\]
which we now turn to consider.
\begin{thm}
\label{v*reg}
Let ${\mathbf{a}}$ be as in (\ref{introd0}), (\ref{ellip}), (\ref{symr}). Let $S_{ {\mathbf{a}} }$ be a fundamental solution of $P[{\mathbf{a}},D]$. 
 Let $\alpha\in]0,1[$. Let $m\in{\mathbb{N}}\setminus\{0\}$. 
Let $\Omega$ be a bounded open subset of ${\mathbb{R}}^{n}$ of class $C^{m,\alpha}$. Then the following statements hold.
\begin{enumerate}
\item[(i)] The operator $w_{*}[\partial\Omega ,{\mathbf{a}}, S_{ {\mathbf{a}} },\cdot]_{|\partial\Omega}$ is linear and continuous from $C^{m-1}(\partial\Omega)$ to $C^{m-1,\omega_{\alpha}(\cdot)}(\partial\Omega)$.
\item[(ii)] Let $\beta\in]0,\alpha]$. Then the operator $w_{*}[\partial\Omega ,{\mathbf{a}}, S_{ {\mathbf{a}} },\cdot]_{|\partial\Omega}$ is linear and continuous from $C^{m-1,\beta}(\partial\Omega)$ to $C^{m-1,\alpha}(\partial\Omega)$.
\end{enumerate}
\end{thm}
{\bf Proof.} We first note that
\begin{eqnarray}
\label{v*reg1}
\lefteqn{
w_{*}[\partial\Omega ,{\mathbf{a}}, S_{ {\mathbf{a}} },\mu](x)
=\sum_{b,r=1}^{n}a_{br}\int_{\partial\Omega}\nu_{r}(x)\frac{\partial}{\partial x_{b}}S_{{\mathbf{a}}}(x-y)\mu(y)\,d\sigma_{y}
}
\\ \nonumber
&&\qquad\qquad
=\sum_{b,r=1}^{n}a_{br}
Q[\frac{\partial S_{ {\mathbf{a}} }}{\partial x_{b}}\circ\Theta,\nu_{r},\mu](x)
\\ \nonumber
&&\qquad\qquad\qquad
+\sum_{b,r=1}^{n}a_{br}\int_{\partial\Omega}\nu_{r}(y)\frac{\partial}{\partial x_{b}}S_{{\mathbf{a}}}(x-y)\mu(y)\,d\sigma_{y}
\\ \nonumber
&& \qquad\qquad
=\sum_{b,r=1}^{n}a_{br}
Q[\frac{\partial S_{ {\mathbf{a}} }}{\partial x_{b}}\circ\Theta,\nu_{r},\mu](x)
\\ \nonumber
&&\qquad\qquad\qquad
-
 \int_{\partial\Omega}\mu(y)\sum_{b,r=1}^{n}a_{br}\nu_{r}(y)\frac{\partial}{\partial y_{b}}S_{{\mathbf{a}}}(x-y) \,d\sigma_{y}
\\ \nonumber
&& \qquad\qquad
=\sum_{b,r=1}^{n}a_{br}
Q[\frac{\partial S_{ {\mathbf{a}} }}{\partial x_{b}}\circ\Theta,\nu_{r},\mu](x)
\\ \nonumber
&&\qquad\qquad\qquad
-
w[\partial\Omega ,{\mathbf{a}}, S_{{\mathbf{a}}} ,\mu ](x)
-
v[\partial\Omega       , S_{{\mathbf{a}}}       ,(a^{(1)}\nu) \mu](x)\,,
\end{eqnarray}
for all $x\in \partial\Omega$ and for all $\mu\in C^{0}(\partial\Omega)$.   

If $m=1$, then Theorem \ref{v0a} implies that $v[\partial\Omega  , S_{{\mathbf{a}}}      
,\cdot]_{\partial\Omega}$ is linear and continuous from $C^{m-1}(\partial\Omega)$ to $C^{m-1,\alpha}(\partial\Omega)$.

If $m> 1$, then $C^{m-1}(\partial\Omega)$ is continuously imbedded into 
$C^{m-2,\alpha}(\partial\Omega)$ and Theorem \ref{slay} implies that 
$v[\partial\Omega  , S_{{\mathbf{a}}}      ,\cdot]_{|\partial\Omega}$ is linear and continuous from $C^{m-2,\alpha}(\partial\Omega)$ to $C^{m-1,\alpha}(\partial\Omega)$. Hence, 
$v[\partial\Omega    , S_{{\mathbf{a}}}     ,\cdot]_{|\partial\Omega}$ is continuous  from the space  $C^{m-1}(\partial\Omega)$ to $C^{m-1,\alpha}(\partial\Omega)$ for all $m\geq 1$. Then formula (\ref{v*reg1}), and the continuity of the imbedding of $C^{m-1,\alpha}(\partial\Omega)$ into $C^{m-1,\omega_{\alpha}}(\partial\Omega)$, and Theorems \ref{qrsm} (i), \ref{wreg} (i)  imply the validity of statement (i). 

We now consider statement (ii). Since 
$v[\partial\Omega    , S_{{\mathbf{a}}}     ,\cdot]_{|\partial\Omega}$ is continuous from 
$C^{m-1,\beta}(\partial\Omega)$ to $C^{m,\beta}(\partial\Omega)$
and 
$C^{m,\beta}(\partial\Omega)$ is continuously imbedded into $C^{m-1,\alpha}(\partial\Omega)$, then the operator 
$v[\partial\Omega   , S_{{\mathbf{a}}}       ,\cdot]_{|\partial\Omega}$ is continuous from 
$C^{m-1,\beta}(\partial\Omega)$ into $C^{m-1,\alpha}(\partial\Omega)$. Then 
formula (\ref{v*reg1}), and  Theorems \ref{qrsm} (ii), \ref{wreg} (ii)  imply the validity of statement (ii). \hfill  $\Box$ 

\vspace{\baselineskip}

Since the space $C^{m-1,\omega_{\alpha}(\cdot)}(\partial\Omega)$ is compactly imbedded into $C^{m-1}(\partial\Omega)$, and $C^{m-1,\alpha}(\partial\Omega)$  is compactly imbedded into $C^{m-1,\beta}(\partial\Omega)$ for all $\beta\in]0,\alpha[$, we have the followng immediate consequence ot Theorem \ref{v*reg} (ii).

\begin{corol}
Under the assumptions of Theorem \ref{v*reg},  
$w_{*}[\partial\Omega ,  S_{ {\mathbf{a}} } ,\cdot]_{|\partial\Omega}$ is compact from $C^{m-1 }(\partial\Omega)$ to itself, and from $C^{m-1,\omega_{\alpha}(\cdot)}(\partial\Omega)$ to itself, and from $C^{m-1,\alpha}(\partial\Omega)$ to itself. 
\end{corol}

\vspace{\baselineskip}

 \noindent
{\bf Acknowledgement} This paper represents an extension of   the work performed 
by F.~Dondi in his `Laurea Magistrale' dissertation under the guidance of 
M.~Lanza de Cristoforis, and contains the results of  \cite{Do12}, \cite{Do14}, \cite{La08}.
The authors are indebted to M.~Dalla Riva for a help in the formulation of Lemma \ref{exr}, and of Theorem
\ref{fundsol}, and of Corollary \ref{ourfs} on the fundamental solution and of Theorems \ref{slay}, \ref{dlay} on layer potentials. The authors are indebted to P. Luzzini for a comment which has improved the statement
of Lemma \ref{k0b}.  The authors also wish to thank D.~Natroshvili for pointing out a number of references.

\noindent
Francesco Dondi\\
Ascent software, Malta Office\\
92/3, Alpha Center,\\
Tarxien Road, 
Luqa, LQA 1815, Malta\\
francesco314@gmail.com\\

\noindent
Massimo Lanza de Cristoforis\\
Dipartimento di Matematica\\
Universit\`{a} degli Studi di Padova\\
Via Trieste 63, Padova 35121, Italy\\
mldc@math.unipd.it\\


\begin{thebibliography}{11}
 

\bibitem{Ci00}
A.~Cialdea, {\em Appunti di teoria del potenziale}.
Corso di Analisi Superiore, A.A. 2001--02, Dipartimento di 
Matematica, Universit\`{a} degli Studi della Basilicata, 2000.

\bibitem{CoKr83}
D.~Colton and R.~Kress, {\em Integral equation methods in scattering 
theory},  Wiley, New York, 1983. 

\bibitem{Da13}
M.~Dalla Riva, {\em A family of fundamental solutions of elliptic partial differential operators with real constant coefficients}. Integral Equations Operator Theory, {\bf 76} (2013),  1--23.

\bibitem{Da14}
M.~Dalla Riva, {\em  Anisotropic heat transmission}. Typewritten manuscript, 2014. 

\bibitem{DaMoMu13} 
 M.~Dalla Riva, J.~Morais, and P.~Musolino, {\em A family of fundamental solutions of elliptic partial differential operators with quaternion constant coefficients}. Math.~Methods Appl.~Sci.,  {\bf 36} (2013),  1569--1582.
 

\bibitem{De85}
 K.~Deimling, {\em Nonlinear Functional Analysis}. 
Springer-Verlag, Berlin, {\em etc.},  1985.

\bibitem{DiMi04}
M.~Dindo\v{s} and  M.~Mitrea, {\em  The stationary Navier-Stokes system in nonsmooth manifolds: the Poisson problem in Lipschitz and $C^1$ domains}. Arch. Ration. Mech. Anal.,  {\bf 174} (2004), 1--47.



\bibitem{Do12}
F.~Dondi, {\em  Comportamento al contorno dei potenziali di strato}. Tesi di laurea triennale, relatore M.~Lanza de Cristoforis, Universit\`{a} degli studi di Padova, pp.~1--57, 2012. 

\bibitem{Do14}
F.~Dondi, {\em  Derivate tangenziali del potenziale di doppio strato per problemi ellittici del secondo ordine a coefficienti costanti}. Tesi di laurea magistrale, relatore M.~Lanza de Cristoforis, Universit\`{a} degli studi di Padova, pp.~1--86, 2014. 

 
\bibitem{Du10}
R.~Duduchava, {\em Lions' lemma, Korn's inequalities and the Lam\'{e}  operator on hypersurfaces}.   Recent Trends in Toeplitz and Pseudodifferential Operators,  43--77. Springer, 2010.


\bibitem{DuMiMi06}
R.~Duduchava, D.~Mitrea, and M.~Mitrea, {\em Differential operators and boundary value problems on hypersurfaces}. Math. Nachr.,  {\bf 279} (2006),  996--1023.

\bibitem{FaJoRi78}
E.B.~Fabes, M.~Jodeit, Jr., N.M.~Rivi\`{e}re, {\em Potential techniques for boundary value problems on C1-domains}, Acta Math. {\bf 141} (1978),  165--186.


 
\bibitem{Fo84}
G.~Folland, {\em Real Analysis, Modern techniques and their 
applications}, John Wiley \& Sons, 1984.


\bibitem{Gu67}
N.M.~G\"{u}nter, {\em
Potential theory and its applications to basic problems of mathematical physics}, 
translated from the Russian by John R. Schulenberger, Frederick Ungar Publishing Co., New York, 1967.


\bibitem{He92}
U.~Heinemann, {\em Die regularisierende Wirkung der Randintegraloperatoren der klassischen Potentialtheorie in den R\"{a}umen h\"{o}lderstetiger Funktionen}, Diplomarbeit, Universit\"{a}t Bayreuth, 1992.


\bibitem{HoMitTa10}
S.~Hofmann, M.~Mitrea and M.~Taylor, {\em 
Singular integrals and elliptic boundary problems on regular Semmes-Kenig-Toro domains}.
Int. Math. Res. Not. IMRN 2010, no. 14, 2567--2865. 

\bibitem{Jo55}
F.~John,  {\em Plane waves and spherical means applied to partial differential
  equations}.  Interscience Publishers, New York-London, 1955.


\bibitem{Ki89}
A.~Kirsch, {\em Surface gradients and continuity properties for some 
integral operators in classical scattering theory}, Math.~Methods  
Appl.~Sciences, {\bf 11} (1989), 789--804. 

\bibitem{KiHe15}
A.~Kirsch and  F.~Hettlich, {\em  The Mathematical Theory of Time-Harmonic Maxwell's Equations; Expansion-, Integral-, and Variational Methods},  Springer International Publishing Switzerland, 2015.  

\bibitem{KuGeBaBu79}
V.D.~Kupradze, T.G.~Gegelia, M.O.~Basheleishvili and 
T.V.~Burchuladze, {\em Three-dimensional Problems of the Mathematical 
Theory of Elasticity and Thermoelasticity}, North-Holland Publ.~Co., Amsterdam, 
1979. 

\bibitem{La08}
M.~Lanza~de~Cristoforis, {\em Properties of the integral operators associated to the layer potentials}. Typewritten handout for students, 2008.

\bibitem{MaSh05}
V.~Maz'ya and T.~Shaposhnikova, {\em Higher regularity in the classical layer potential theory for Lipschitz domains}. Indiana University Mathematics Journal, {\bf 54} (2005),   99--142.

\bibitem{Mik70}
S.G.~Mikhlin, {\em Mathematical physics, an advanced course},  translated from the Russian, North-Holland Publishing Co., Amsterdam-London,  1970.

	
\bibitem{Mi65}
C.~Miranda, {\em Sulle propriet\`{a} di regolarit\`{a} di certe 
trasformazioni integrali}, Atti Accad. Naz.   
Lincei Mem. Cl. Sci. Fis. Mat. Natur. Sez I, {\bf 7} (1965), 303--336. 

\bibitem{Mi70} 
C.~Miranda, {\em Partial differential equations of elliptic type, 
Second revised edition},   Springer-Verlag, Berlin,  1970.


\bibitem{Mit14}
M.~Mitrea, {\em The almighty double layer: recent perspectives}. Invited presentation at the 13$th$ International Conference on Integral methods in Science and Engineering, July 21--25, 2014. 

\bibitem{MitMit13}
I.~Mitrea and M.~Mitrea, {\em Multi-Layer Potentials and Boundary Problems, for Higher-Order Elliptic Systems in Lipschitz Domains}.
Lecture Notes in Mathematics, Springer, Berlin, {\it etc.} 2013.


\bibitem{MitMitVe14}
D.~Mitrea, M.~Mitrea  and J.~Verdera. {\em Characterizing Lyapunov Domains via Riesz Transforms on H\"{o}lder Spaces}. Preprint 18. October 2014. 


\bibitem{Ne67}
J.~Ne\v cas, {\em Les m\'{e}thodes directes en th\'{e}orie des \'{e}quations elliptiques.} Masson et Cie,  Paris 1967.

\bibitem{Sc31}
J.~Schauder,  {\em Potentialtheoretische Untersuchungen}, Math. 
Z., {\bf 33}  (1931),  602--640.

\bibitem{Sc32}
J.~Schauder,  {\em Bemerkung zu meiner Arbeit ``Potentialtheoretische Untersuchungen I (Anhang)''}, Math.  Z., {\bf 35} (1932), 536--538.

\bibitem{Schi82}
H.~Schippers, {\em  On the regularity of the principal value of the double-layer potential}. Journal of Engineering Mathematics, {\bf 16} (1982), 59--76.


\bibitem{ScWe99}
C.~Schwab and W.~Wendland, {\em  On the extraction technique in boundary integral
equations}. Math. Comp. {\bf 68} (1999), 92--122.

\bibitem{Tr87}
G.M.~Troianiello, {\em Elliptic Differential Equations 
and Obstacle Problems}, Plenum Press, New York and London, 1987.

\bibitem{vo90}
W.~von Wahl, {\em Absch\"{a}tzungen f\"{u}r das Neumann-\-Problem
und die Helm\-holtz-\-Zerlegung von $L^{p}$},  Nachr. Akad. Wiss. G\"{o}ttingen Math.-Phys. Kl. II, {\bf 2},   1990.   

\bibitem{Wi93}
M.~Wiegner, {\em Schauder estimates for boundary layer potentials}, Math. Methods Appl.~Sci., {\bf 16}   (1993),   877--894.
 






\end{thebibliography}
\end{document}